\documentclass[matann,final,a4paper]{svjour}[28/09/2006]

\usepackage[T1]{fontenc}

\usepackage{amssymb}
\usepackage{mathrsfs}
\usepackage{amsmath}

\input{macroperso}

\spnewtheorem*{notation}{Notation}{\bf}{\normalsize}
\spnewtheorem{theorem-section}{Theorem}[subsection]{\bf}{\normalsize}
\spnewtheorem*{theorem-intro}{Theorem.}{\bf}{\normalsize}

\newcommand{\adj}[1]{\hspace{0.02in}^{\mathbf{t}}\overline{#1}}
\newcommand{\tr}{\mathrm{tr}}

\newcommand{\normop}[1]{\vert \vert \vert #1 \vert \vert \vert}
\renewcommand{\Im}{\mathrm{Im}}

\smartqed 

\date{}

\begin{document}
\setcounter{tocdepth}{3}

\bibliographystyle{spmpsci}

\title{Vortex type equations and canonical metrics}
\author{Julien Keller}
\institute{\textsc{Imperial College - London} \\ \email{j.keller@imperial.ac.uk, keller@cict.fr}}
\maketitle
\begin{abstract}
We introduce a notion of Gieseker stability for a filtered holomorphic vector bundle $\mathcal{F}$ over a projective manifold. We relate it to an analytic condition in terms of hermitian metrics on $\mathcal{F}$ coming from a construction of the Geometric Invariant Theory (G.I.T). We prove that if there is a $\boldsymbol{\tau}$-Hermite-Einstein metric $h_{HE}$ on $\mathcal{F}$,
then there exists a sequence of such balanced metrics that converges and its limit is $h_{HE}$. As a corollary, we obtain an approximation theorem for quiver Vortex equations and other classical equations.
\end{abstract}

\tableofcontents

\section{Introduction}

Intrinsic global methods have appeared in K\"ahler geometry with the fundamental work of E.\ Calabi and S-T.\ Yau in the study of Einstein type metrics. Recently, S.K. Donaldson in \cite{Do3} has made a crucial advance in view of an algebro-geometric interpretation of the existence of K\"ahler-Einstein metrics by studying a notion of stability introduced by H. Luo for the couples $\left( M,L^{k}\right) $ where $M$ is projective and $L$ is a polarization over $M$. The main result of \cite{Do3} (also see \cite{Bi1} for a survey) is the convergence of a sequence of \lq balanced' metrics, constructed algebraically via the embeddings in the projective spaces $\mathbb{P}H^{0}\left(M,L^{k}\right)$ for $k$ large, towards a constant scalar metric (when its existence is assumed \textit{a priori}). \\
 In the case of bundles (with or without some extra decorated structures), such correspondences have been proved in most of the known cases, at least when $M$ is compact. Therefore, one can ask if the method of \cite{Do3,Do2} can be applied to give an approximation of the equations that appear for these objects. In particular, this means that the gobal analysis used in Hitchin-Kobayashi-Donaldson-Uhlenbeck-Yau (we will say HKDUY in short) type correspondences 
contains some key technics for the construction of algebraic objects. By \lq algebraic objects', we mean some canonical metrics coming from an underlying G.I.T construction in a finite dimensional framework.
 A preliminary work has been done by C. Drouet \cite{Dr} for the case of Hermite-Einstein equations over a curve, and later X. Wang \cite{W1,W2} has given a complete solution of the problem in any dimension using Gieseker's results and Donaldson's breakthrough.
  
  In this paper, we are interested in a more general type of equation, so called Vortex equations, that appear in bosonic theories for non linear $\sigma$-models under the form of Bogomol'nyi-Prasad-Sommerfeld states. In the K\"ahler setting, and considering only a single bundle, these equations have been studied essentially by C.H.\ Taubes, S.\ Bradlow \cite{Bra2}, O. Garc\'{i}a-Prada \cite{GP0,GP1} and later D.\ Banfield \cite{Ba} who proposed a generic frame. In \cite{GP} (see also \cite{B-GP2}), O. Garc\'{i}a-Prada introduced a notion of coupled Vortex equations relative to a triple (i.e two holomorphic bundles and a bundle map between them). Later, this idea of considering Vortex equations for some data involving more than one bundle or sheaf has been developped in the work of  L.\ {\'A}lvarez-C{\'o}nsul and O.\ Garc\'{i}a-Prada   \cite{AC-GP,AC-GP2,AC-GP3} who also investigated the related notion of quivers. Nevertheless, a G.I.T construction for the objects studied by Banfield or {\'A}lvarez-C{\'o}nsul and Garc\'{i}a-Prada is still to be achieved in full generality, even if some recent progress has been made in \cite{G-S,Sc1}. From a general point of view, the study of Vortex equations and their moduli space of solutions has had some important consequences in algebraic geometry (for instance the Velinde formula and its generalizations) or in Gauge theory with the computation of Gromov-Witten invariants \cite{OT0}.

\medskip
In our work, we will focus on the $\boldsymbol{\tau}$-Hermite-Einstein equations introduced by \'{A}lvarez-C\'{o}nsul and Garc\'{i}a-Prada \cite{AC-GP} for a holomorphic filtration $\mathscr{F}$  (or filtered holomorphic bundle) over a smooth projective manifold. In a first part, we introduce a weaker notion of stability for holomorphic filtrations and we relate it to  certain Gieseker spaces using a G.I.T construction (Theorem \ref{theo1}). We use a Kempf-Ness type argument to show that the stability condition can be written as an analytic condition involving the Bergman kernel associated to the filtration (Theorem \ref{Br403a}). Thus, a sequence of canonical metrics will be defined, called balanced metrics as foreseen in \cite{Do2} (Definition \ref{equili}). Then we use an idea of Donaldson to find zeros of moment maps and we study the combined action of the Gauge group and the special unitary group $SU(N)$. Finally, we will prove in the fourth part (Theorem \ref{final}),
\begin{theorem-intro} \label{central}
Let $M$ be a smooth projective manifold. If $\mathscr{F}$ is an irreducible holomorphic filtration of a holomorphic vector bundle $\mathcal{F}$ over $M$ equipped with a metric $h_{HE}$ solution of the $\boldsymbol{\tau}$-Hermite-Einstein equation
$$\sqrt{-1}\Lambda F_{h_{HE}} = \sum_i \widetilde{\tau}_i\pi^{i}_{h_{HE}}(\mathscr{F}),$$
then there exists a sequence of balanced  metrics $h_k$ on $\mathcal{F}$ which converges, up to conformal change, towards the metric $h_{HE}$ in $C^{\infty}$ topology.
\end{theorem-intro}
As a corollary, we give a procedure to get an approximation of solutions of quiver Vortex equations (Theorem \ref{final2}) using a dimensional reduction argument based on \cite{AC-GP,AC-GP3}. In particular, these equations overlap with Hermite-Einstein equations, special Witten triples (non abelian monopoles) and
critical Higgs equations over a curve (see Section \ref{lastpar} for details).

\section{Preliminaries}

\subsection{Notions of stability}

The aim of this part is to introduce different notions of stability for a holomorphic filtration on a projective manifold. We will follow here the ideas of D. Gieseker and the Mumford theory in order to introduce a Gieseker space for which the G.I.T stable points correspond to stable holomorphic filtrations.   \newline

Let $M$ be a projective manifold of complex dimension $n$, and $L$ an ample line bundle on $M$.

\begin{definition}
A filtration of sheaves of length $m$ is a finite sequence of coherent subsheaves $$\mathscr{F}:0 = \mathcal{F}_0 \hookrightarrow \mathcal{F}_1 \hookrightarrow ...
\hookrightarrow \mathcal{F}_m=\mathcal{F}$$ and we shall say that $\mathscr{F}$ is
a holomorphic filtration if the sheaves  $\mathcal{F}_i$ and $\mathcal{F}$ are subbundles.
\end{definition}

\begin{definition}
A subfiltration of the filtration $\mathscr{F}$ is a filtration of sheaves of length
 $m$
$$\mathscr{F}':0 \hookrightarrow \mathcal{F}_1' \hookrightarrow ...
\hookrightarrow \mathcal{F}_m'=\mathcal{F}',$$
where $\mathcal{F}'$  is a subsheaf of $\mathcal{F}$ and such that $\mathcal{F}_i'\subseteq\mathcal{F}_i \cap \mathcal{F}'$ for any $1\leq i \leq m$. A
subfiltration is said to be proper if $r(\mathcal{F}') < r(\mathcal{F})$.
\end{definition}

\begin{definition}
A holomorphic filtration $\mathscr{F}$ is irreducible if it cannot be written as
$$\mathscr{F}=\mathscr{F}_1 \oplus \mathscr{F}_2, $$
where the $\mathscr{F}_i\neq \mathscr{F}$ are holomorphic subfiltrations. 
\end{definition}

\begin{definition}
We will say that the filtration $\mathscr{F}$ is simple if any endomorphism $f \in End(\mathcal{F})$
which preserves the filtration (i.e. $f(\mathcal{F}_i)\subset \mathcal{F}_i)$ is a scalar
multiple of $Id$.
\end{definition}

We recall the notion of (slope) stability for a filtration (we refer to \cite{AC-GP,B-GP4}
for details).
\begin{definition} \label{Defpaire}
Let $\mathscr{F}$ be a filtration of length $m$ and $\boldsymbol{\tau}=(\tau_1,...,\tau_{m-1})$
a $(m-1)$-tuple of real numbers. Define the $\boldsymbol{\tau}$-degree of $\mathscr{F}$
as
$$\deg_{\boldsymbol{\tau}}(\mathscr{F})=\deg(\mathcal{F})+\tsum_{i=1}^{m-1}\tau_i r(\mathcal{F}_i),$$
and the $\boldsymbol{\tau}$-slope of $\mathscr{F}$ as
$$\mu_{\boldsymbol{\tau}}(\mathscr{F})=\frac{\deg_{\boldsymbol{\tau}}(\mathscr{F})}{r(\mathcal{F})}.$$
We shall say that the filtration $\mathscr{F}$ is $\boldsymbol{\tau}$-stable (resp. semi stable)
if for all proper subfiltration $\mathscr{F}'\hookrightarrow \mathscr{F}$, we have $$\mu_{\boldsymbol{\tau}}(\mathscr{F}') < \mu_{\boldsymbol{\tau}}(\mathscr{F}) \hspace{1cm}(\text{resp. }\leq).$$
A filtration is said to be polystable if it is a direct sum of  $\boldsymbol{\tau}$-stable filtrations with same slope $\mu_{\boldsymbol{\tau}}$.
\end{definition}

\begin{remark}
If the filtration reduces to a bundle $\mathcal{F}$ (i.e. $m=1$) we recover the case of Mumford stability for bundles with $\tau_i=0$ and in this case, we will denote by $\mu(\mathcal{F})$ the slope of $\mathcal{F}$. Notice that the $\boldsymbol{\tau}$-stability of a holomorphic filtration does not imply the Mumford stability of the subbundles.
\end{remark}

\begin{lemma} \label{lemme1}
Let $\mathscr{F}^1$ and $\mathscr{F}^2$ be two filtrations of torsion free sheaves with same length,  $\boldsymbol{\tau}$-stable and of same slope. Let $\varrho:\mathscr{F}^1 \rightarrow \mathscr{F}^2$ a non zero homomorphism such that for all $i$, $\varrho(\mathcal{F}^1_i) \subset \mathcal{F}^2_i$. Then $\varrho$ is injective. In particular, if a holomorphic filtration is stable then it is simple.
\end{lemma}
\begin{proof}
If $\varrho$ is non injective then $\mathcal{F}^3=\Im(\varrho)$ is a (proper) torsion free quotient of $\mathcal{F}^1$ and we get  $\mathscr{F}^3$ subfiltration of 
$\mathscr{F}^2$ such that $$\mu_{\boldsymbol{\tau}}(\mathscr{F}^3)>\mu_{\boldsymbol{\tau}}(\mathscr{F}^1)=\mu_{\boldsymbol{\tau}}(\mathscr{F}^2)$$ by Whitney product formula. Since $\mathscr{F}^2$ is stable
and $\mathscr{F}^3 \subset \mathscr{F}^2$, one has necessarily that $r(\mathcal{F}^3)=r(\mathcal{F}^2)$. \\
 Now it is also clear that for a holomorphic subfiltration $\mathscr{F}'\subset \mathscr{F}$ such that $r(\mathcal{F})=r(\mathcal{F}')$, the following inequality always holds
\begin{equation}
\mu_{\boldsymbol{\tau}}(\mathscr{F}') \leq \mu_{\boldsymbol{\tau}}(\mathscr{F}) \label{stabil}.
\end{equation}
This comes from the fact that if $\mathcal{F}' \ncong \mathcal{F}$ there exists an effective divisor $D$ such that $\det(\mathcal{F}) \cong \det(\mathcal{F}'\otimes \mathcal{O}_M(D))$ and so $\mu(\mathcal{F})=\mu(\mathcal{F}')+\frac{\deg(D)}{r(\mathcal{F}')}>\mu(\mathcal{F}')$.
In the case of non holomorphic filtrations (i.e. given by torsion free sheaves), inequality (\ref{stabil}) remains true. Then we get a contradiction: $\varrho$ is injective. \\
\- \quad If $\mathscr{F}^1$ and $\mathscr{F}^2$ are two holomorphic filtrations, then $\mathscr{F}^3$ is a holomorphic subfiltration of $\mathscr{F}^2$, with same slope and $r(\mathcal{F}^3)=r(\mathcal{F}^2)$. We notice that in the case of holomorphic  filtrations, the only case of equality in (\ref{stabil}) is $\mathscr{F}'= \mathscr{F}$. Thus, $\varrho$ is an isomorphism. Therefore, if $\mathscr{F}$ is a  holomorphic stable filtration, any non zero endomorphism $\varrho$ of $\mathscr{F}$ such that $\varrho(\mathcal{F}_i) \subset \mathcal{F}_i$ is an isomorphism and by Schur lemma, $\{\varrho \in End(\mathcal{F}): \varrho(\mathcal{F}_i) \subset \mathcal{F}_i \}$ is a divison algebra of finite dimension and isomorphic to $\mathbb{C}$. \qed
\end{proof} 

\begin{notation}
For a holomorphic filtration $\mathscr{F}$ of length $m$ and  a hermitian metric  $h$ on $\mathcal{F}$ correspond $h$-orthogonal smooth projections
on the subbundle $\mathcal{F}_{i}$ of $\mathcal{F}$ for all $1\leq i\leq m$, 
that we shall  note $$\pi_{h,i}^{\mathscr{F}}:\mathcal{F}\rightarrow \mathcal{F}_i$$
with the convention $\pi_{h,m}^{\mathscr{F}}=Id_{\mathcal{F}}.$ 
\end{notation}

The main result of  \cite{AC-GP} is the existence of a HKDUY correspondence for holomorphic filtrations in terms of metrics of
 $\boldsymbol{\tau}$-Hermite-Einstein type.

\begin{theorem-section} Fix $\omega $ a  K\"{a}hler metric on the compact manifold $%
M$. Let $\boldsymbol{\tau}\in \mathbb{R}^{m-1}_{+}$ and $\mathscr{F}$ be a holomorphic filtration of length $m$.  A holomorphic filtration $\mathscr{F}$ is $\boldsymbol{\tau}$-polystable
if and only if there exists a smooth hermitian metric $h$ solution of the equation
\begin{equation}
\label{Br904} 
{\sqrt{-1}}\Lambda _{\omega }F_{h}+ \sum_{i=1}^{m-1} {\tau}_i \pi_{h,i}^{\mathscr{F}} = \mu_{\boldsymbol{\tau}}(\mathscr{F})Id_{\mathcal{F}}.
\end{equation}
This can be also written as
\begin{equation}
{\sqrt{-1}}\Lambda _{\omega }F_{h}= \sum_{i=1}^{m} \widetilde{\tau}_i \pi_{h}^{i}(\mathscr{F}), \label{neweqn}
\end{equation}
where $$\widetilde{\tau}_i=\mu_{\boldsymbol{\tau}}(\mathscr{F})-\tsum_{j=i}^{m-1} \tau_j, \hspace{1cm}
\widetilde{\tau}_{m}=\mu_{\boldsymbol{\tau}}(\mathscr{F}),$$
and $\pi_{h}^{i}(\mathscr{F})=\pi_{h,i}^{\mathscr{F}}-\pi_{h,i-1}^{\mathscr{F}}$ is the  projection (with respect to $h$) on the orthogonal of the subbundle $\mathcal{F}_{i-1}$ of $\mathcal{F}_i$ with the convention $\pi_{h}^{1}(\mathscr{F})=\pi_{h,1}^{\mathscr{F}}$.
We shall say under these conditions that  $h$ is $\boldsymbol{\tau}$-Hermite-Einstein. The filtration $\mathscr{F}$ will be said $\boldsymbol{\tau}$-Hermite-Einstein.
\end{theorem-section}

\begin{remark}
The assumption of non negativity of the real numbers $\tau_i$ is crucial in the proof of
 \cite[Theorem 2.1]{AC-GP}. Henceforth, when we consider the stability of a filtration  of length $m$, it is relatively to a  $(m-1)$-tuple of non negative real numbers.
If we take the trace of  (\ref{neweqn}), we notice that the parameters 
$\widetilde{\tau}_i$ satisfy the relation
$$\tsum_{i=1}^{m} \widetilde{\tau}_i r(\mathcal{F}_i/\mathcal{F}_{i-1})=\deg(\mathcal{F}),$$
which means that we only have $m-1$ degrees of freedom as expected.
\end{remark}

We introduce now a notion of Gieseker stability for filtrations.
\begin{definition}
Let $\mathbf{R}=(R_1,..,R_{m-1})$ be a collection of $(m-1)$ polynomials with rational coefficients of degrees $d_i<n$ and positive for $k$ large. Define $$\mathscr{P}_{\mathbf{R},\mathscr{F}}(k)=\chi(\mathcal{F}\otimes L^k)+\tsum_{i=1}^{m-1} r(\mathcal{F}_i)R_i(k)$$ the $\mathbf{R}$-Hilbert polynomial of the filtration $\mathscr{F}$ of length $m$. Then $\mathscr{F}$ is said to be Gieseker $\mathbf{R}$-stable (resp. semi-stable) if for $k$ large, one has for all proper subfiltration $\mathscr{F}'$ of $\mathscr{F}$,
\begin{eqnarray*}
 \frac{\mathscr{P}_{\mathbf{R},\mathscr{F}'}(k)}{r(\mathcal{F}')} &<& \frac{\mathscr{P}_{\mathbf{R},\mathscr{F}}(k)}{r(\mathcal{F})}  \hspace{1cm} (\text{resp.} \leq).
\end{eqnarray*}
\end{definition}

We get immediately by applying Riemann-Roch formula,
\begin{proposition}
If the filtration $\mathscr{F}$ is $\boldsymbol{\tau}$-Mumford stable, then it is also  $\mathbf{R}$-Gieseker stable for
\begin{equation*}
R_i(k)=\tau_i k^{n-1} +O(k^{n-2}).
\end{equation*}
\end{proposition}

\subsection{G.I.T construction and extended Gieseker space for filtrations}

We present a G.I.T frame for the  holomorphic filtrations on a projective manifold, inspired from the work of \cite{H&L3,H&L2,Sc1}. We introduce a \lq Gieseker space' which parametrizes the Gieseker stable holomorphic filtrations. Firstly, we notice that the considered Gieseker semi-stable objects live in bounded families, i.e. are parametrized by a scheme of finite type over $\mathbb{C}$.

\begin{proposition}
The set of isomorphy classes of Gieseker semi-stable filtrations given by torsion free coherent sheaves with fixed $\mathbf{R}$-Hilbert polynomial is bounded.
\end{proposition}
\begin{proof}
Let us consider a filtration $\mathscr{F}$. By considering the leading coefficients of $\mathscr{P}_{\mathbf{R},\mathcal{F}}(k)$ and using the Gieseker semi-stability condition, one reaches the conclusion that the slopes $\mu(\mathcal{F}_i)$ are all bounded.
Now, we can apply the boundedness theorem of \cite[Section 3]{M2}. \qed 
\end{proof}
 
Fix $\mathbf{R}$ a collection of $(m-1)$ polynomials as before and let $\mathscr{F}$ be a holomorphic filtration of length $m$ and note $r_i=r(\mathcal{F}_i)$ with $r_0=0$, $r=r(\mathcal{F})$ and $\mathfrak{p}$ its $\mathbf{R}$-Hilbert polynomial. By Kodaira's embedding  theorem, there exists an integer $k_0$ such that for $k\geq k_0$,
the bundles $\mathcal{F}_i\otimes L^k$ are globally generated and the cohomology groups of higher dimension of $\mathcal{F}_i\otimes L^k$ are trivial, i.e.  $$H^{j}(M,\mathcal{F}_i\otimes L^k)=0, \hspace{1cm}\forall j\geq 1. $$ For such a $k$, consider a vector space $V$ isomorphic to $H^0(M,\mathcal{F}\otimes L^k)$, and let $(v_i)$ be a basis of $V$. One can construct a quasi-projective quot scheme  $\mathfrak{Q}'$ parametrizing equivalence classes of quotients $\{q:V\otimes L^{-k}\rightarrow \mathfrak{F}\}$ where $\mathfrak{F}$ is a filtration of torsion free coherent sheaves of length $m$ with $\mathbf{R}$-Hilbert polynomial equal to $\mathfrak{p}$ and $H^0(q\otimes id_{L^k})$ is an isomorphism. This gives us  universal quotients $\tilde{q_i}:V\otimes \pi_M^*(L^{-k})\rightarrow \tilde{\mathfrak{F}_i}$ over $\mathfrak{Q}\times M$ where $\mathfrak{Q}$ is the union of components of $\mathfrak{Q}'$ that contain $\mathbf{R}$-semi-stable elements. The line bundles $\det(\tilde{\mathfrak{F}_i})$ induce morphisms $\upsilon_i:\mathfrak{Q} \rightarrow Pic(M)$ and we denote $\mathfrak{A}_i$ the union of the finitely many components of $Pic(M)$ hit by the morphism $\upsilon_i$. By previous proposition, this does not depend on the choice of $k$ and we can  assume that $k_0$ is large enough so that for all $[\mathcal{L}]\in \mathfrak{A}_i$, $\mathcal{L}\otimes L^{(r-r_i)k}$ is globally generated and without higher cohomology. Fix a Poincar\'{e} line bundle $\tilde{\mathcal{L}}$ on $Pic(M)\times M$ and note $\tilde{\mathcal{L}}_{\mathfrak{A}_i}$ its restriction to  $\mathfrak{A}_i\times M$.
Again by previous proposition, we can assume that $\tilde{\mathcal{L}}_{\mathfrak{A}_i}\otimes \pi_M^*L^{(r-r_i) k}$
is globablly generated and without higher cohomology for $k \geq k_0$. Introduce the extended Gieseker space $${\mathfrak{G}}=\prod_{i=0}^{m-1} \mathbb{P} \left(Hom \left(\wedge^{r-r_i} V\otimes \mathcal{O}_{\mathfrak{A}_i},(\pi_{\mathfrak{A}_i})_*(\tilde{\mathcal{L}}_{\mathfrak{A}_i}\otimes \pi_M^*L^{(r-r_i) k})\right)^\vee \right).$$ 
The  morphism $(\pi_{\mathfrak{Q}})_*(\wedge^{r-r_i}( \tilde{q_i} \otimes id_{\pi_M^*L^{k}}):\wedge^{r-r_i} V\otimes \mathcal{O}_{\mathfrak{Q}} \rightarrow (\pi_{\mathfrak{Q}})_*(\det(\tilde{\mathfrak{F}_i})\otimes \pi_M^*L^{(r-r_i)k})$ gives rise to an injective and $SL(V)$-equivariant morphism $$Gies:\mathfrak{Q}\rightarrow {\mathfrak{G}}.$$ Moreover, the Gieseker space $\mathfrak{G}$ maps $SL(V)$-invariantly to $\prod_i \mathfrak{A}_i$
and the fibers are close $SL(V)$-invariant subschemes; i.e., over the point $(\mathcal{L}_1,..,\mathcal{L}_m)\in \prod_i \mathfrak{A}_i$, sits the space $$\widetilde{\mathfrak{G}}_{(\mathcal{L}_1,..,\mathcal{L}_m)}=\prod_{i} \mathbb{P} \left(Hom \left(\wedge^{r-r_i} V,H^0(\mathcal{L}_i\otimes L^{(r-r_i) k})\right)^\vee \right).$$
We now focus on determining the (semi-)stable points of this space. For the holomorphic filtration $\mathscr{F}$ of lentgh $m$, let $\pi_i$ be the natural onto map  from $\mathcal{F}$ to the quotient $\mathcal{F}\otimes L^k / \mathcal{F}_i \otimes L^k$ with the convention $\pi_0=Id$. To the filtration $\mathscr{F}$, we associate by the following morphisms of evaluation,
\begin{eqnarray*}
T_i&:& (v_{j_1} \wedge ... \wedge v_{j_{r-r_i}}) \mapsto \left( p\mapsto \pi_i ev_p(v_{j_1}) \wedge ... \wedge \pi_i  ev_p(v_{j_{r-r_i}})\right),
\end{eqnarray*}
a point in the space
\begin{equation*}
\widetilde{\mathfrak{G}}_k=\prod_{i=0}^{m-1} \mathbb{P}Hom \left(\wedge^{r-r_i} V, H^0\left(M,\det(\mathcal{F}\otimes L^k/ \mathcal{F}_i \otimes L^k)\right)\right).
\end{equation*}
The action of $SL(V)$ is given by $$g \star (T_0,...,T_{m-1})=\left(T_0 \circ \wedge^r g^{-1}, ..., T_{m-1}\circ\wedge^{r-r_{m-1}} g^{-1} \right),$$ and we can consider 
for a choice of parameters $\varepsilon_i>0$ the G.I.T stability of a point of $\widetilde{\mathfrak{G}}_k$ relatively to a  $SL(V)$-linearization of the very ample bundle 
 $\mathcal{O}_{\widetilde{\mathfrak{G}}_k}(\varepsilon_0,...,\varepsilon_{m-1})$.

 Let $\underline{\lambda}:G_t\rightarrow SL(V)$ be a 1-parameter subgroup and ${v}_i$ a basis of ${V}$ such that $G_t$ acts on  ${V}$ by $\underline{\lambda}$ with the weights $\gamma_i \in \mathbb{Z}$, i.e. for all $t\in G_t$,
\begin{equation*}
\underline{\lambda}(t)\cdot {v}_i=t^{\gamma_i}{v}_i.
\end{equation*}
Of course, it can be assumed that $\gamma_i \leq \gamma_{i+1}$ and $\tsum_i 
\gamma_i=0$. For any multi-index $I=(i_1,...,i_{r-r_j})$ of length $|I|=r-r_j$ with  $1\leq i_1 <...<i_{r-r_j}\leq \dim({V})$, let ${v}_I={v}_{i_1} \wedge ... \wedge {v}_{i_{r-r_j}}$ and 
$\gamma_I=\gamma_{i_1}+...+\gamma_{i_{r-r_j}}.$ The group $SL({V})$ acts on $\bigwedge^{r-r_j} V$ with weight $\gamma_I$ relatively to the basis $v_I$.
The classical Hilbert-Mumford criterion asserts that the point of $\widetilde{\mathfrak{G}}_k$ is G.I.T-(semi-)stable relatively to the linearization that we have fixed and the action of $SL({V})$ if and only if, for all 1-parameter subgroup,
\begin{equation*}
\tsum_{j=0}^{m-1}\varepsilon_j\min_{\{I:|I|=r-r_j\}}\lbrace\gamma_I:T_j({v}_I)\neq 0\rbrace
 <0
\hspace{1cm} \text{(resp. }\leq).
\end{equation*}

\begin{remark}
 With Riemann-Roch theorem, we compute the dimension of $V$,
 \begin{eqnarray*}
 \dim(V) &=& \int_M Ch(\mathcal{F}\otimes L^k)Todd(M) \\
 &=& r(\mathcal{F})k^n \int_M \frac{c_1(L)^n}{n!} + k^{n-1}\int_M \left( \frac{r}{2}c_1(M)+c_1(\mathcal{F})\right)\frac{c_1(L)^{n-1}}{(n-1)!}+...
 \end{eqnarray*}
\end{remark}  

Consider now the case of a subspace $V'\subset V$ and the action of the subgroup  associated is given by
\begin{equation*}
\lambda(t)=\left(
\begin{array}{ccc}
t^{-\mathrm{codim}(V')} & & 0\\
 & \ddots &  \\
0 & & t^{\dim(V')} \\
\end{array}\right),
\end{equation*}
with $\gamma_1\hspace{-0.05cm}=\hspace{-0.05cm}..\hspace{-0.05cm}=\hspace{-0.05cm}\gamma_{\dim(V')}\hspace{-0.05cm}=\hspace{-0.05cm}-\mathrm{codim}(V')$ and $\gamma_{\dim(V')+1}\hspace{-0.05cm}=\hspace{-0.05cm}..\hspace{-0.05cm}=\gamma_{\dim(V)}\hspace{-0.05cm}=\hspace{-0.05cm}\dim(V')$.  Via the morphism ${V}\otimes \mathcal{O}_M\rightarrow \mathcal{F}\otimes L^k$ obtained for $k$ sufficiently large, we have a holomorphic filtration
\begin{eqnarray*}
0\subsetneq \mathcal{F}_{(V')} \subset \mathcal{F}_{(V)}=\mathcal{F}, \\
\mathcal{F}_{(V'),i}=\mathcal{F}_{i}\cap \mathcal{F}_{(V')}
\end{eqnarray*}
and under these conditions,
\begin{eqnarray*}
\hspace{-0.1cm}\min_{\{|I|=r-r_j\}} \lbrace\gamma_I:T_j({v}_I)\neq 0\rbrace
       &=&\dim(V')\left(r(\mathcal{F})-r(\mathcal{F}_i)\right) \\ &&-\dim(V)\left(r(\mathcal{F}_{(V')})-r(\mathcal{F}_{(V'),i})\right).
\end{eqnarray*}
Thus, we have shown that if the point of $\widetilde{\mathfrak{G}}_k$, defined by the filtration $\mathscr{F}$, is G.I.T (semi-)stable, then we have
\begin{eqnarray*}
\varepsilon\left(\dim(V')r(\mathcal{F})-\dim(V)r(\mathcal{F}_{(V')})\right)\\
& \hspace{-5cm}+ \sum_{i=0}^{m-1} \varepsilon_i \left(\dim(V)r(\mathcal{F}_{(V'),i})-\dim(V')r(\mathcal{F}_{i})\right)<
0,
\end{eqnarray*}
where we have set $$\varepsilon=\tsum_{i=0}^{m-1}\varepsilon_i.$$
In fact, we also have the converse:
\begin{lemma}
\label{Br901b} The point of $\widetilde{\mathfrak{G}}_k$ defined by the  holomorphic filtration $\mathcal{F}$ is G.I.T (semi-)stable if and only if for any subspace $V' \subset V$,
\begin{eqnarray*}
\varepsilon \left(\dim(V')r(\mathcal{F})-\dim(V)r(\mathcal{F}_{(V')}) \right)\\
& \hspace{-5cm}+ \sum_{i=0}^{m-1}  \varepsilon_i \left(\dim(V)r(\mathcal{F}_{(V'),i})-\dim(V')r(\mathcal{F}_{i})\right)<
0,
\end{eqnarray*}
where $\mathscr{F}_{(V')}$ is the filtration generated by $V'\otimes \mathcal{O}_M$ and $\mathcal{F}/\mathcal{F}_{(V')}$ 
is torsion free.
\end{lemma}

\begin{proof}
Let $({v}_1,...,{v}_{\dim(V)})$ be a  basis of ${V}$. If one denotes $\mathcal{F}_{(i)}=\mathcal{F}_{\langle v_1,...,v_i \rangle}$, we have a filtration
\begin{equation*}
\mathcal{F}_{(0)} \subset ... \subset \mathcal{F}_{\dim({V})}=\mathcal{F}
\end{equation*}
and we get $\mathcal{F}_{(i)}=\mathcal{F}_{(i-1)}$ or $r(\mathcal{F}_{(i)})>r(\mathcal{F}_{(i-1)})$. Consequently,
there exist $r$ integers between $1$ and $\dim({V})$ that mark the jumps of
ranks. We shall denote them $(k_1,..,k_r)$. If one considers the action associated to  $(v_i,\gamma_i)$, we get by \cite[Lemma 1.23]{H&L2}, $\min_I\lbrace\gamma_I:T_0({v}_I)\neq 0\rbrace=\gamma_{k_1}+...+\gamma_{k_r}$. From a similar way, there exist $r-r_j$ integers between  $1$ and $\dim({V})$ that mark the jumps of ranks for $\mathcal{F}\otimes L^k/\mathcal{F}_{(i),j}\otimes L^k$. We denote them  $(k_1^j,..,k_{r-r_j}^j)$. Therefore, we get  $$\min_I\lbrace\gamma_I:T_j({v}_I)\neq 0\rbrace=\gamma_{k_1^j}+...+\gamma_{k_{r-r_j}^j}.$$
In order to apply the Hilbert-Mumford criterion, we shall consider the set of all the vectors with weight $\gamma_i$. Of course, these vectors are  generated by the following vectors
\begin{equation*}
\gamma^{(i)}=(\underbrace{i-\dim({V}),\ldots,i-\dim(V)}_{i},\underbrace{i,\ldots,i}_{\dim(V)-i})
\end{equation*}
for $i=1,...,\dim(V)$. All weighted vectors  $\gamma=\left(\gamma_1,...,\gamma_{\dim({V})}\right)$ can be expressed as
$\gamma= \sum_{i=1}^{\dim({V})}c_i \gamma^{(i)}$ where $c_i= \frac{\gamma_{i+1}-\gamma_i}{\dim({V})}$ are non negative rational coefficients.
Let us apply the Hilbert-Mumford criterion to $\gamma^{(i)}$; we get
\begin{eqnarray*}
\mu^{(i)}&=&{\tsum}_{j=0}^{m-1}\varepsilon_j\min_{\{I:|I|=r-r_j\}}\lbrace\gamma_I:T_j({v}_I)\neq 0\rbrace \\
&=&
  -\dim({V})\left(\displaystyle{\tsum}_j
  \varepsilon_j\max_l\{k_l^j \leq i \}\right)   +i\left(\displaystyle{\tsum}_j (r-r_j) \varepsilon_j\right).
\end{eqnarray*}
If $i$ grows,  $\mu^{(i)}$ goes down except for $k_j$ or one $k_j^l$. One has to evaluate  $\mu^{(i)}$ at the values $k_j-1$ or $k_j^l-1$. This leads to
\begin{equation*}
\mu^{(i)}=-\dim(V)\displaystyle{\tsum}_l \varepsilon_l(j^l-1)
+\displaystyle{\tsum}_l \varepsilon_l (r-r_l)(k_j^l-1).
\end{equation*}
Eventually, we can forget the choice of the basis ${v}_i$ and we get that the point of $\widetilde{\mathfrak{G}}_k$ is G.I.T stable (resp. semi-stable) if and only if
\begin{eqnarray*}
\varepsilon\dim(V')r(\mathcal{F})-\dim(V')\displaystyle{\tsum}_{i}\varepsilon_i r(\mathcal{F}_{i}) &<&
\varepsilon\dim(V)r(\mathcal{F}_{(V')}) \\
&& \hspace{0.5cm} -\dim(V)\displaystyle{\tsum}_{i}\varepsilon_i r(\mathcal{F}_{(V'),i}) 
\end{eqnarray*}
(resp. $\leq$) for any subspace $0 \neq V'\subset {V},$ with $r \left(\mathcal{F}_{(V')}\right)\leq r(\mathcal{F})$.\qed
\end{proof}

\- \newline By \cite[Lemma 1.26]{H&L2}, the G.I.T stability (resp. semi-stability) can be written as a condition on the subsheaves of $\mathcal{F}$ instead of the subspaces of $V$,
\begin{eqnarray*}
\dim({V}\cap H^0(\mathcal{F}'\otimes L^k))\left(\varepsilon r(\mathcal{F})-
\textstyle\sum_i \varepsilon_i r(\mathcal{F}_i)\right) &\\
& \hspace{-5cm} <\dim({V})\left(\varepsilon r(\mathcal{F}')- \textstyle\sum_i \varepsilon_i r \left(\mathcal{F}_{
V \cap H^0(\mathcal{F}'\otimes L^k),i}\right)\right) 
\end{eqnarray*}
(resp. $\leq$) for any proper holomorphic filtration $\mathscr{F}' \subset \mathscr{F}$. Indeed, if $\mathcal{F}' \subset \mathcal{F}$, choose $V'=H^0(\mathcal{F}'\otimes L^k) \cap {V}.$ Then for the morphism $q:V\otimes \mathcal{O}_M \rightarrow \mathcal{F} \otimes L^k$, $q(V'\otimes \mathcal{O}_M)\subset \mathcal{F}'\otimes L^k$ and $r(\mathcal{F}')=r(\mathcal{F}_{(V')})$.  For the converse, we consider $V'\subset {V}$ and set $\mathcal{F}'=q(V'\otimes \mathcal{O}_M)$. Thus, we have $V'\subset {V} \cap H^0(\mathcal{F}'\otimes L^k)$ and $r(\mathcal{F}_{(V')})=r(\mathcal{F}')$.

Now, it is clear that the Gieseker stability condition for the holomorphic filtration $\mathcal{F}$ implies the previous condition for a convenient choice of $\{\varepsilon,\varepsilon_1,...,\varepsilon_{m-1}\}$, i.e.  we have proved the
\begin{theorem} \label{theo1}
Let $\mathbf{R}=(R_1,...,R_{m-1})$ be a collection of $(m-1)$ polynomials with rational  coefficients of degree $d_i<n$ and with positive leading coefficient.
The holomorphic filtration $\mathscr{F}$ of length $m$ is $\mathbf{R}$-stable (resp. semi-stable) if for $k$ large, the associated point of $\widetilde{\mathfrak{G}}_k$ is G.I.T stable (resp. semi-stable) respectively to the polarization $\mathcal{O}_{\widetilde{\mathfrak{G}}_k}(\varepsilon_0,...,\varepsilon_{m-1})$ and the  action of $SL(V)$ where we have fixed 
\begin{eqnarray*}
\varepsilon_0&=&1-\sum_{i=1}^{m-1} \frac{R_i(k)}{k^n}, \\
\varepsilon_i&=&\frac{R_i(k)}{k^n} \hspace{1cm}(1\leq i\leq m-1).
\end{eqnarray*}
\end{theorem}
\begin{remark}
This G.I.T construction extends the previous work of \cite{D-U-W} 
in which a G.I.T construction is given for one step filtrations
and a similar result is proved.
\end{remark}

\section{G.I.T stability and balanced metrics}

In this part, we apply the Kempf-Ness criterion to the Gieseker spaces
that we have just constructed for holomorphic filtrations.
This means that the condition of G.I.T stability can be transposed
as the existence of a certain sequence of \lq balanced' metrics defined on the
finite dimension vector space $H^0(\mathcal{F}\otimes L^k)$ that are  
critical points of certain functionals of Kempf-Ness type.  \\
 The balance condition for applications was conceptualized by S.K.\ Donaldson \cite{Do2} in the following way. Let's assume that we are given the following objects: a holomorphic map $f:\Xi \rightarrow W$ where
$(\Xi,\omega_0)$ is compact K\"ahler and $W$ is a vector space of finite dimension
embedded by $\pi^W$ as a co-adjoint orbit in 
$Lie(G)^*$ where $G$ is reductive linear. Then, the center of mass of 
 $f$ in $Lie(G)^*$ is given by
\begin{equation*}
\int_{Lie(G)^*}\pi^W_*\left( f_{\ast }\left(\frac{\omega_0^n}{n!}\right)\right).
\end{equation*}
 Under this setting, $f$ is said to be balanced if the orbit of $f$ under the action of   $Lie(G)^*$ contains a center of mass null. This means that we demand that the moment map defined by integration on $C^\infty(\Xi,W)$ respectively
  to the action of $G$ admits a zero in the complex orbit of $f$.

In all the following, our polarization $L$ will be equipped with a smooth hermitian metric $h_{L}$ such that the curvature $c_{1}\left( L,h_{L}\right) $ 
is a K\"ahler metric on $M$ that we denote $\omega,$ which means that
\begin{equation*}
\omega =-\frac{i}{2\pi }\partial \overline{\partial }\log \left(
h_{L}\right).
\end{equation*}
Let $dV=\frac{\omega ^{n}}{n!}$ be the volume form relatively to $\omega$.

\begin{notation}
We define $Met(\Upsilon)$ the space of smooth hermitian metrics for 
the vector space or the vector bundle $\Upsilon$.  Let $F$ be
a hermitian vector bundle on $(M,\omega)$. We associate to a metric $h \in Met(F)$ the $L^2$ hilbertian metric on $H^0(M,F)$ 
\begin{eqnarray*}
Hilb_{\omega}(h)&=&\int_M \langle .,. \rangle_{h} \frac{\omega^n}{n!}\in Met(H^0(M,F)).
\end{eqnarray*}
\end{notation}

We will also need the well-known fact:
\begin{definition}  \label{re}
If $V_1$ and $V_2$ are two vector spaces of finite dimension $N_1,N_2$ equipped with metrics $h_1,h_2$, and $T:V_1 \rightarrow V_2$  is a linear map, then the Hilbert-Schmidt  norm $||T||_{h_1,h_2}$ can be computed using any othonormal basis  $(v^1_i)_{i=1,..,N_1}$ of $V_1$ by $||T||_{h_1,h_2}^2=\sum_{i=1}^{N_1}|T(v^1_i)|^2_{h_2}$.
\end{definition}

\subsection{Gieseker stable holomorphic filtrations and Kempf-Ness type functionals}

We fix a holomorphic filtration $\mathscr{F}$ of length $m$ and for $k$ large, we set $N=h^{0}\left(M,\mathcal{F}\otimes L^{k}\right)$ and $\pi_i$ the 
natural surjection onto
$\mathcal{F}\otimes L^k/ \mathcal{F}_i\otimes L^k$ with the convention that $\pi_0=Id$. Let  $\mathbf{\widetilde{Z}_{ss}}\subset \widetilde{\mathfrak{G}}_k$ be the
open scheme of G.I.T semi-stable points with respect to the 
linearization $\mathcal{O}_{\widetilde{\mathfrak{G}}_k}(\varepsilon_0,..,\varepsilon_{m-1})$
where the constants $\varepsilon_i$ are fixed by Theorem \ref{theo1}. \\

Fix a metric $H$ on the vector space $H^0(\mathcal{F}\otimes
L^k)$. We have a quotient metric on the hermitian bundle  $\mathcal{F}\otimes L^k$ induced by
$V\twoheadrightarrow \mathcal{F}\otimes L^k$ and consequently a
metric $\left\vert \left\vert . \right\vert \right\vert$ on 
$\Lambda^{r-r_i} (\mathcal{F}\otimes L^k)$. By the isomorphism $i:V \tilde{\rightarrow} H^0(M,\mathcal{F}\otimes L^k)$, we get a metric $h_\mathbf{\widetilde{V}_{k,i}}$ on the space
$\mathbf{\widetilde{V}}_{k,i}=Hom(\wedge^{r-r_i} V,H^0(\det(\mathcal{F}\otimes L^k/ \mathcal{F}_i\otimes L^k)))$. We set 
$$\mathbf{\widetilde{V}}_k= \mathbf{\widetilde{V}}_{k,0} \times ... \times \mathbf{\widetilde{V}}_{k,m-1},$$
and get a natural metric  $\left\vert \left\vert \left\vert .,.\right\vert
\right\vert
\right\vert_{\widetilde{\mathbf{V}}_k}=h_{\mathbf{\widetilde{V}}_{k,0}}^{\varepsilon_0}\times ... \times
h^{\varepsilon_{m-1}}_{\mathbf{\widetilde{V}}_{k,m-1}}$. Let $\mathbf{z}$ be a point 
in $\mathbf{\widetilde{Z}_{ss}}$ and $\mathbf{\tilde{z}}\in
\mathcal{O}_{\widetilde{\mathfrak{G}}_k}(-\varepsilon_0,...,-
\varepsilon_{m-1})_{\mathbf{z}}$ be a lifting.
We can evaluate the metric  $\left\vert \left\vert \left\vert
.,.\right\vert \right\vert \right\vert_{\widetilde{\mathbf{V}}_k}$ at that point,
{\small
\begin{eqnarray*}
\left\vert \left\vert \left\vert \mathbf{\tilde{z}}\right\vert
\right\vert \right\vert _{\widetilde{\mathbf{V}}_k}^{2}\hspace{-0.1cm}=  C(i)
\prod_{j=0}^{m-1} \left(\int_{M}\sum_{\tiny{\begin{array}{l}
                                    1\hspace{-0.1cm}\leq \hspace{-0.1cm}i_{1}\hspace{-0.1cm}<\hspace{-0.1cm}... \\
			        ...\hspace{-0.1cm}<\hspace{-0.1cm}i_{r-r_j}\hspace{-0.1cm}\leq\hspace{-0.1cm} N
                                  \end{array}}}
				  \hspace{-0.35cm}\left\Vert \pi_j \circ s_{i_{1}}\left( p\right) \wedge ...\wedge
\pi_j \circ s_{i_{r-r_j}}\left( p\right) \right\Vert ^{2}dV\left( p\right)\right)^{\varepsilon_j}
\end{eqnarray*}}\normalsize where $(s_{i}) _{i=1,..,N}$ is an 
$H$-orthonormal basis of $H^{0}\left(
M,\mathcal{F}\otimes L^{k}\right)$, $C(i)>0$ is a constant depending only on the isomorphism $i$.
\begin{remark}
Our construction does not depend on the choice of the metric on the determinant bundles $\det(\mathcal{F}\otimes L^k/ \mathcal{F}_i\otimes L^k)$.
\end{remark}

\begin{definition}
We define the functional for $g \in SL(V)$, 
\begin{eqnarray*}
\widetilde{F_{k,\mathscr{F}}}(g)\hspace{-0.09cm}=\hspace{-0.09cm} 
{\sum_{j=0}^{m-1} \varepsilon_j \log {\int_M} \frac{ \sum_{
  \tiny{\begin{array}{l}
                                    1\hspace{-0.1cm}\leq \hspace{-0.1cm}i_{1}\hspace{-0.1cm}<\hspace{-0.1cm}... \\
			        ...\hspace{-0.1cm}<\hspace{-0.1cm}i_{r-r_j}\hspace{-0.1cm}\leq\hspace{-0.1cm} N
                                  \end{array}}}\left\Vert \pi_j (g \cdot s_{i_{1}})\wedge ...\wedge
\pi_j (g \cdot s_{i_{r-r_j}}) \right\Vert ^{2}} { \sum_{\tiny{\begin{array}{l}
                                    1\hspace{-0.1cm}\leq \hspace{-0.1cm}i_{1}\hspace{-0.1cm}<\hspace{-0.1cm}... \\
			        ...\hspace{-0.1cm}<\hspace{-0.1cm}i_{r-r_j}\hspace{-0.1cm}\leq\hspace{-0.1cm} N
                                  \end{array}}}\left\Vert \pi_j s_{i_{1}}\wedge ...\wedge
\pi_j s_{i_{r-r_j}} \right\Vert ^{2}}  dV }.
\end{eqnarray*}

\end{definition}
We can sum up our situation by the following lemma:
\begin{lemma}
The following conditions are equivalent:\\
-- The holomorphic filtration $\mathscr{F}$ is $\mathbf{R}$-Gieseker
polystable for a collection $\mathbf{R}=(R_1,...,R_{m-1})$ of rational polynomials of  degree $n-1$ with positive leading coefficient. \\
-- There exists an integer $k_0$ such that for all $k \geq k_0$, the functionals
$\widetilde{F_{k,\mathscr{F}}}:SL(V)\rightarrow \mathbb{R}$
admit a  positive minimum where it has been assumed that
\begin{equation*}
\varepsilon_i=\frac{R_i(k)}{k^n}, \hspace{1cm} \varepsilon_0=1-\sum_{i=1}^{m-1} \frac{R_i(k)}{k^n}  \label{csts}
\end{equation*}
for all $i=1,...,m-1$. 
\end{lemma}
\begin{remark}
If the holomorphic filtration $\mathscr{F}$ is
$\boldsymbol{\tau}$-Mumford stable then there exists an integer $k_0$ such that for all $k \geq k_0$, the functionals 
$\widetilde{F_{k,\mathscr{F}}}$ admit a positive minimum and are proper, under the assumption of $\varepsilon_0=1-\sum_{i=1}^{m-1} \frac{\tau_i}{k} $ and
$\varepsilon_i=\frac{\tau_i}{k}$.
\end{remark}

Note that we can also consider $\widetilde{F_{k,\mathscr{F}}}$
as a functional on the space $Met(V)\times
SL(V)$,  i.e. on a finite dimensional space. We are going to see 
 that we can relate to the functional 
$\widetilde{F_{k,\mathscr{F}}}$ another functional, this time on the
infinite dimensional space $Met(\mathcal{F})\times
SL(V)$. This motivates the following definition and theorem (see \cite[Section 3]{MR} for details).

\begin{definition}
Consider a K\"ahler manifold $(\Xi,\omega)$ and a moment map  $\mu$ associated to the action of a compact linear group $\Gamma$ such that $\Gamma^{\mathbb{C}}$ acts  holomorphically. To the moment map $\mu$  corresponds a functional
\begin{equation*}
I_{\mu}:\Xi \times \Gamma ^{\mathbb{C}}\rightarrow\mathbb{R}\end{equation*}
called the \lq\textit{integral of the moment map} $\mu$\rq \
and which satisfies the properties: \\
-- For all $p\in \Xi,$ the critical points of the restriction $I_{\mu,p}$ of $I_{\mu}$ to $\{p\}\times \Gamma ^{\mathbb{C}}$ coincide with the points in the orbit $\Gamma ^{
\mathbb{C}}p$ for which the moment map vanishes; \\
-- The restriction $I_{\mu,p}$ to the lines \ $\{e^{\lambda u}:u\in
\mathbb{R}\}$ where $\lambda \in Lie\left( \Gamma ^{%
\mathbb{C}}\right) $ is convex.
\end{definition}

\begin{theorem-section} \label{102}
There exists a unique application $I_{\mu}:\Xi \times \Gamma ^{\mathbb{C}}\rightarrow\mathbb{R}
$ which satisfies the two properties: \\
1. $I_{\mu} \left( p,e\right) =0$ for all $p\in \Xi;$ \\
2. $\frac{d}{du}I_{\mu} \left( p,e^{i\lambda u}\right) _{|u=0}=\langle\mu \left(
p\right),\lambda \rangle$ for all $\lambda \in Lie\left( \Gamma \right).$
\end{theorem-section}
 
Now, for $k$ sufficiently large, we have the embeddings of $M$ into the Grassmanians of $r-r_j$ quotients defined by:
\begin{equation}
i_{k,j}:
\begin{array}{ccl}
M & \hookrightarrow & Gr\left(N,r-r_j\right) \\
p & \mapsto & \ker\left( \pi_j \circ ev_{p}: V
\rightarrow \mathcal{F}\otimes L^{k}/\mathcal{F}_j\otimes L^{k} |_{p}\right) ^{\vee }.
\end{array}  \label{i_k3}
\end{equation}
Let $\underline{\mathbf{U}}_{N,r}$ be the universal bundle on the Grassmannian of  $r$-quotients of the Grassmannian
$Gr(N,r)$. We denote $\boldsymbol{\mathrm{\Pi}}=\prod_{i=0}^{m-1} Gr(N,r-r_i)$ and $\pi_{Gr,i}:\boldsymbol{\mathrm{\Pi}} \rightarrow Gr(N,r-r_i)$ for $i=0,...,m-1$ the natural projections. We lift the Fubini-Study metrics on each factor $Gr(N,r-r_i)$ with weight $\varepsilon_i$. This induces a symplectic metric on 
$C^{\infty}\left(M,\boldsymbol{\mathrm{\Pi}}\right)$,
$$\Omega_{(i_{k,*})}(\overrightarrow{x},\overrightarrow{y})= \sum_{i=0}^{m-1} \int_M \varepsilon_i \pi_{Gr,i}^* \omega_{FS}(\overrightarrow{x},\overrightarrow{y}) dV,$$
where $\overrightarrow{x},\overrightarrow{y} \in C^{\infty}\left(M,(i_{k,0},...,i_{k,m-1})^*T\boldsymbol{\mathrm{\Pi}}\right)$. The moment map associated to $\Omega_{(i_{k,*})}$ for the 
action of the special unitary group $SU(N)$ on $C^{\infty}\left(M,\boldsymbol{\mathrm{\Pi}}\right)$ is
\begin{eqnarray*}
\mu_{\mathscr{F},k}(i_{k,*})&=&\int_M \tsum_{j=0}^{m-1} \varepsilon_j  \mathsf{Q}_j\adj{\mathsf{Q}_j} dV-V\frac{\sum_{j=0}^{m-1} (r-r_j)\varepsilon_j}{N}Id,\\
&=&\int_M (\varepsilon_0+\tsum_{j=1}^{m-1} \varepsilon_j)  \mathsf{Q}_0\adj{\mathsf{Q}_0} dV -
\int_M \tsum_{j=1}^{m-1} \varepsilon_j \left( \mathsf{Q}_0\adj{\mathsf{Q}_0}- \mathsf{Q}_j\adj{\mathsf{Q}_j} \right) dV \\
&& -V\frac{r- \sum_{j=0}^{m-1} r_j\varepsilon_j}{N}Id, 
\end{eqnarray*}
\normalsize
where $[\mathsf{Q}_j]$ represents a point of $Gr(N,r-r_j)$ i.e.  $\mathsf{Q}_j:\mathcal{F}\otimes L^k / \mathcal{F}_j \otimes L^k_{|p}\rightarrow V$  is an
isometry respectively to $h$ and $H$, and represents the matrix of the
 endomorphism $\pi_{j} \circ ev_p$ expressed in an orthonormal basis  of $\ker (\pi_j \circ ev_p)^\perp$  and in an orthonormal basis of $V$. \\
Set ${\mathbf{U}}_{r,{N}}=\pi_{Gr,0}^*\underline{\mathbf{U}}_{N,r}$. Since $\mathcal{F}\otimes L^k \simeq j_k^*\mathbf{U}_{r,{N}}$
where $j_k: M\hookrightarrow \boldsymbol{\mathrm{\Pi}}$ is induced by the maps $i_{k,l}$ and $\pi_i$, we get a new smooth hermitian metric on $\mathcal{F} \otimes L^k$
associated to the metric $H$ on $H^0(M,\mathcal{F}\otimes L^k)$.

\begin{definition}
Let  $FS_{k}=FS_{k}(H)\in Met(\mathcal{F}\otimes L^k)$ be the hermitian metric
on $\mathcal{F}\otimes L^k$ induced by
\begin{equation} \label{defFS}
\langle .,. \rangle_{FS_{k}}=\Big\langle \frac{N}{Vr-V\sum_{j=1}^{m-1}\varepsilon_j r_j}\left( Id_{\mathcal{F}} - \displaystyle{\sum}_{j=1}^{m-1}\varepsilon_j
\pi^{\mathscr{F}}_{h,j}\right).,.\Big\rangle_h
\end{equation}
where $h$ is the quotient metric on $\mathcal{F}\otimes L^k$ induced by $H$.
\end{definition}

\begin{remark}
This last expression is well defined since we get for  $k$ sufficiently large, $0<\sum_{j=1}^{m-1}\varepsilon_j <1$.
\end{remark}

\begin{definition}
Define the functional on $SL(V)$,
{
\begin{eqnarray*}
\widetilde{KN_{k,\mathscr{F}}}(g)\hspace{-0.05cm}&=&\sum_{j=0}^{m-1}\frac{\varepsilon_j}{2} \int_M \log
\frac{ \sum_{\tiny{
\begin{array}{l}
                                    1\hspace{-0.1cm}\leq \hspace{-0.1cm}i_{1}\hspace{-0.1cm}<\hspace{-0.1cm}... \\
			        ...\hspace{-0.1cm}<\hspace{-0.1cm}i_{r-r_j}\hspace{-0.1cm}\leq\hspace{-0.1cm} N
                                  \end{array}}}
				  \left\Vert  \pi_j (g \cdot s_{i_{1}})\wedge ...\wedge
 \pi_j (g \cdot s_{i_{r_j}})\right\Vert ^{2}} { \sum_{\tiny{
 \begin{array}{l}
                                    1\hspace{-0.1cm}\leq \hspace{-0.1cm}i_{1}\hspace{-0.1cm}<\hspace{-0.1cm}... \\
			        ...\hspace{-0.1cm}<\hspace{-0.1cm}i_{r-r_j}\hspace{-0.1cm}\leq\hspace{-0.1cm} N
                                  \end{array}}}\left\Vert \pi_j s_{i_{1}}\wedge ...\wedge
\pi_j s_{i_{r_j}} \right\Vert ^{2} } dV
\end{eqnarray*}
} \normalsize
where $(s_{i}) _{i=1,..,N}$ is any  $H$-orthonormal basis of $H^{0}\left(
M,\mathcal{F}\otimes L^{k}\right)$.
\end{definition}

\begin{lemma}
$\widetilde{KN_{k,\mathscr{F}}}$ is the integral of the moment map $\mu_{\mathscr{F},k}$.
\end{lemma}
\begin{proof}
From \cite{Mok}, we know that a potential of the Fubini-Study metric at the point
$[\mathsf{Q}]$ in the Grassmannian $Gr(N,r)$ is 
\begin{equation*}
\log \det \left( \adj{\mathsf{Q}}\mathsf{Q}\right).
\end{equation*}
In order to prove the lemma, we simply need to check that for any trace free  matrix $\mathsf{S}$,
\begin{equation*}
\frac{d}{du}\left(\widetilde{KN_{k,\mathscr{F}}}\left(ge^{\mathsf{S}u}\right)\right)_{|u=0}=\mu_{\mathscr{F},k}(g) \in SL({N}).
\end{equation*}
Let $[\mathsf{Q}_0(p)]$ represent the point in the Grassmannian
$Gr({N},r)$ given by the embedding defined by (\ref{i_k3}),
at $p\in M$. Therefore, we obtain up to a modification of the matrix $\mathsf{Q}_0(p)$ (i.e by considering instead the unitary matrix
$\mathsf{Q}_0\left(^t\overline{\mathsf{Q}_0}\mathsf{Q}_0\right)^{-1/2}$),
\begin{equation*}\int_{M} \log \left(\frac{\sum_{\tiny{
\begin{array}{l}
                                    1\hspace{-0.1cm}\leq \hspace{-0.1cm}i_{1}\hspace{-0.1cm}<\hspace{-0.1cm}... \\
			        ...\hspace{-0.1cm}<\hspace{-0.1cm}i_{r_j}\hspace{-0.1cm}\leq\hspace{-0.1cm} N
                                  \end{array}}}\left\Vert g \cdot s_{i_{1}}\wedge ...\wedge
g \cdot s_{i_{r}} \right\Vert ^{2}}{\sum_{\tiny{\begin{array}{l}
                                    1\hspace{-0.1cm}\leq \hspace{-0.1cm}i_{1}\hspace{-0.1cm}<\hspace{-0.1cm}... \\
			        ...\hspace{-0.1cm}<\hspace{-0.1cm}i_{r_j}\hspace{-0.1cm}\leq\hspace{-0.1cm} N
                                  \end{array}}}\left\Vert s_{i_{1}}\wedge ...\wedge
 s_{i_{r}} \right\Vert ^{2}}\right)dV \hspace{-0.1cm}= \hspace{-0.1cm} \int_M \log  \det  \left( ^{\mathsf{t}}\overline{\mathsf{Q}_0}^{\mathsf{t}}\overline{g}g%
\mathsf{Q}_0\right)dV
\end{equation*}
\normalsize
with $^{\mathsf{t}}\overline{\mathsf{Q}_0(p)}\mathsf{Q}_0(p)=Id_{r\times r}$. \\ 
 From a similar way, let $[\mathsf{Q}_i]$ be the point of $Gr(N,r-r_i)$ induced by
$\pi_i \circ ev_p$.
Thus,
\begin{equation*}
\widetilde{KN_{k,\mathscr{F}}}(g)=\frac{1}{2}\left(\tsum_{j=0}^{m-1} \varepsilon_j \int_{M}\log \frac{\det \left( 
^{\mathsf{t}}\overline{\mathsf{Q%
}_j}\hspace{0.04in}^{\mathsf{t}}\overline{g}g\mathsf{Q}_j\right) }{\det \left( ^{%
\mathsf{t}}\overline{\mathsf{Q}_j}\mathsf{Q}_j\right) }dV \right).
 \end{equation*}
Therefore, we get that for any trace free matrix $\mathsf{S}$ and $g\in SU({N})$,
\small{
\begin{eqnarray*}
\frac{d}{du}\left( \widetilde{KN_{k,\mathscr{F}}}\left( e^{\mathsf{S}u}\right) \right)\hspace{-0.1cm}
&=&  \hspace{-0.1cm}\tsum_{j=0}^{m-1} \frac{\varepsilon_j}{2} \int_{M}\hspace{-0.1cm}\tr\left( \left( ^{%
\mathsf{t}}\overline{\mathsf{Q}_j}\hspace{0.04in}e^{^{\mathsf{t}}\overline{%
\mathsf{S}}u}e^{\mathsf{S}u}\mathsf{Q}_j\right) ^{-1}\hspace{-0.1cm}\left( ^{\mathsf{t}}%
\overline{\mathsf{Q}_j}e^{^{\mathsf{t}}\overline{\mathsf{S}}%
u}\left(\mathsf{S}+^{\mathsf{t}}\hspace{-0.04in}\overline{\mathsf{S}}%
\right) e^{\mathsf{S}u}\mathsf{Q}_j\right) \right) \hspace{-0.03in}dV
\end{eqnarray*}}
\normalsize and therefore,
\begin{eqnarray*}
\frac{d}{du}\left( \widetilde{KN_{k,\mathscr{F}}}\left( e^{\mathsf{S}u}\right) \right)_{u=0}\hspace{-0.3cm} &=&
\tsum_j \varepsilon_j \int_{M}\mathtt{tr}(^{\mathsf{t}}\overline{\mathsf{Q}_j}\mathsf{S}\mathsf{Q}_j),\\
&=&\tsum_j  \varepsilon_j  \int_{M}\mathtt{tr}(^{\mathsf{t}}\overline{\mathsf{Q}_j}\hspace{0.04in}^{\mathsf{t}}\mathsf{S}\mathsf{Q}_j)
- \tsum_j \varepsilon_j \frac{r-r_j}{{N}}%
\int_{M}\mathtt{tr}\left( \mathsf{S}\right), \\
&=& \langle \mu_{\mathscr{F},k}\left(\mathsf{Q}_0,...,\mathsf{Q}_{m-1}  \right),\mathsf{S}\rangle.
\end{eqnarray*}
\normalsize
Finally, $\widetilde{KN_{k,\mathscr{F}}}(Id)=0$, which allows us to conclude.\qed
\end{proof}

\subsection{Balanced metrics for holomorphic filtrations}

We are going to see that the two functionals $\widetilde{KN_{k,\mathscr{F}}}$ and $\widetilde{F_{k,\mathscr{F}}}$ are simultaneously proper. At this point, we will need 
the following classical result of potential theory,
\begin{theorem-section} \label{Br905} Set $Ka(M,\omega^{\prime})=\{\varphi \in C^{\infty}(M,\mathbb{R}):\omega^{\prime} +i\partial \overline{\partial }\varphi>0\}$. There exist some constants $\alpha_{M},C(M,\omega,\omega') >0$ such that for all $\varphi \in Ka(M,\omega'),$ one has
\begin{equation*}
\int_{M}{e^{-\alpha_M (\varphi-\sup_M \varphi)}}\frac{\omega^n}{n!}\leq C.
\end{equation*}
\end{theorem-section}

\begin{lemma} \label{Br928a}
There exist some constants $(\gamma_i)_{i=1..3}$ such that for all $g\in SL\left( N\right)$
\begin{equation*}
\widetilde{F_{k,\mathscr{F}}}(g)-\gamma_1 \geq \frac{1}{\gamma_3}\widetilde{KN_{k,\mathscr{F}}}(g)\geq \widetilde{F_{k,\mathscr{F}}}(g)-\gamma_2.
\end{equation*}
\end{lemma}
\begin{proof}
Let $s_i$ be a  basis of $H^0(\mathcal{F}\otimes L^k)$, at the point $p$,
\begin{eqnarray*}
\varphi_j(p) &=&\log\tsum_{1\leq i_{1}<...<i_{r-r_j}\leq N}\left\Vert 
\pi_j (g \cdot s_{i_{1}} \left(
p\right)) \wedge ...\wedge  \pi_j (g \cdot s_{i_{r-r_j}}(p))\right\Vert ^{2}.
\end{eqnarray*}%
Then  $\varphi_j$  belongs to the K\"ahler cone  $Ka(M,c_1(\det(\mathcal{F}\otimes L^k / \mathcal{F}_j \otimes L^k)))$,  and Theorem \ref{Br905} asserts that there exist two real constants $\alpha_M>0$ and $C>1$ such that
\begin{equation*}
\int_{M}e^{-\alpha_M \left( \varphi_j -\sup_{M}\varphi_j \right) }\frac{\omega
^{n}}{n!}<C,
\end{equation*}%
 which implies that
\begin{equation*}
\log \left( \int_{M}e^{-\alpha_M \left( \varphi_j -\sup_{M}\varphi_j \right) }%
\frac{\omega ^{n}}{n!}\right) <C'.
\end{equation*}%
Now by concavity of $\log ,$ 
\begin{eqnarray*}
\int_{M}\varphi_j dV&\geq &\int_{M}\left(
\sup_{_{M}}\varphi_j \right) dV-\frac{1}{\beta(M,k,\omega) }, \\
&\geq &V\log \left( \sup_{p\in M}\tsum_{i_{1}<...<i_{r-r_j}}\left\Vert
\pi_j (g  s_{i_{1}}\left( p\right)) \wedge ...\wedge \pi_j (g s_{i_{r-r_j}}\left( p\right))
\right\Vert ^{2}\right) \\
&& -\frac{1}{\beta(M,k,\omega) }.
\end{eqnarray*}
Indeed, by concavity of $\log$, we have also
\begin{eqnarray*}
\log \int_M \tsum_{1\leq i_{1}<...<i_{r-r_j}\leq N}\left\Vert
\pi_j (g \cdot s_{i_{1}}) \wedge ...\wedge \pi_j(g \cdot s_{i_{r-r_j}})
\right\Vert ^{2} dV \geq \int_M \varphi_j dV.
\end{eqnarray*}%
Now, summing previous inequalities for all $j$, we obtain the lemma with $\gamma_i$ depending on the data $\{k,\mathscr{F},L,dV\}$. \qed
\end{proof}

\begin{definition}[Balanced metrics for holomorphic filtrations]  \label{equili} \\
-- Let $\mathbf{p}=(i_{k,0},...,i_{k,m-1}) \in \boldsymbol{\mathrm{\Pi}}$ be the point induced by the metric $H\in Met(H^0(M,\mathcal{F}\otimes L^k))$. If $\mu_{\mathscr{F},k}(\mathbf{p})=0$ then the holomorphic filtration $\mathscr{F}$ and $H$ are said to be $k$-balanced. \\
-- If $H\in Met(H^0(M,\mathcal{F}\otimes L^k))$ is $k$-balanced, the
metric $h\in Met(\mathcal{F}\otimes L^k)$  given by
$h=FS_{k}(H)$ is said to be $k$-balanced. 
\item  We will say that the filtration
$\mathscr{F}$ is balanced if there exists an integer $k_0$ such that 
for all $k \geq k_0$, $\mathscr{F}$ is
$k$-balanced.
\end{definition}

Since the sections $s_i\in H^0(M,\mathcal{F}\otimes L^k)$ are also
coordinate sections of the universal bundle,  we see that the
 metric $H\in Met(H^0(M,\mathcal{F}\otimes L^k))$  is
balanced if and only if it is a fixed point of the map
$Hilb_{\omega} \circ FS_{k}$. From a similar way, $h\in Met(\mathcal{F}\otimes L^k)$
is $k$-balanced if and only if it is a fixed point 
of the map $FS_{k} \circ Hilb_{\omega}$. \\

The balanced condition for a holomorphic filtration $\mathscr{F}$
can be translated in terms of the  Bergman kernel of the bundle $\mathcal{F}$. We shall now make explicit what we mean by \lq Bergman kernel' in the following definition.
\begin{definition}
The Bergman kernel of a globally generated bundle $(F,h_F)$ is an endomorphism of the 
bundle associated to the  $L^{2}$ orthonormal projection
 from the space of sections $L^{2}\left( M,F\right)$ onto the space of holomorphic sections  $H^{0}\left(M,F \right)$,
\begin{equation}
\widehat{\mathsf{B}}_{F,{h_F}}=\tsum_{i=1}^{h^0(M,F)}s_{i}\langle
.,s_{i}\rangle_{h_F}\in C^{\infty}(M,End(F))
\end{equation}
where $s_i$ is any basis of $H^0(M,F)$, orthonormal for $Hilb_{\omega}(h_F)$.
\end{definition}

The following result will allow us to consider the balanced condition on the space of infinite dimension $Met(\mathcal{F})$ as the vanishing of a certain moment map related to the action of the Gauge group of the bundle $\mathcal{F}$. We will denote $\mathscr{F}\otimes L^k$ the associated holomorphic filtration obtained by tensorizing each subbundle of the filtration $\mathscr{F}$ by $L^k$.

\begin{lemma}
\label{Br60aa}The holomorphic filtration $\mathscr{F} $ of length $m$ is balanced if and only if there exists an integer $k_{0}>0$ such that for all $k\geq k_{0},$ there
exists a hermitian balanced metric $h_{k}\in Met(\mathcal{F}\otimes L^k)$  such that we have
\begin{equation}
\widehat{\mathsf{B}}_{\mathcal{F} \otimes L^k,h_{k}}+\epsilon_k\sum_{j=1}^{m-1} \varepsilon_j \pi^{\mathscr{F}\otimes L^k}_{j,h_k}= \frac{N+\epsilon_k \sum_{j=1}^{m-1} \varepsilon_j r_j}{rV} Id_{\mathcal{F}\otimes L^{k}}
\label{BrH2c}
\end{equation}
where 
$\epsilon_k=\frac{\chi(\mathcal{F}\otimes L^k)}{Vr-V\sum_{j=1}^{m-1} \varepsilon_j r_j}$.
\end{lemma}

\begin{proof}
Assume that $H$ is a balanced metric on 
$H^0(M,\mathcal{F}\otimes L^k)$ and $s_i$ is an orthonormal basis of
the space $H^0(M,\mathcal{F}\otimes L^k)$ for $H$. Let $h$  be the quotient
metric on $\mathcal{F}\otimes L^k$ induced by
$V\twoheadrightarrow \mathcal{F}\otimes L^k$ and
$FS_{k}(H)$ the metric on $\mathcal{F}\otimes L^k$ constructed
as before. Remember at that point that the Bergman kernel is independent of the choice
of an orthonormal basis. We choose now a $H$-orthonormal basis of sections 
of $H^0(\mathcal{F}\otimes L^k)$ by the following procedure: let $s_{1},...,s_{r_1}$ be
orthogonal to the kernel of $(\pi_{1,h}^{\mathscr{F}\otimes L^k} \circ ev_p)$ with $H$-norm 1, 
and $s_{r_1+1},...,s_{r_2} \in \ker(\pi_{1,h}^{\mathscr{F}\otimes L^k} \circ ev_p)$ be orthogonal  to the kernel of $(\pi_{2,h}^{\mathscr{F}\otimes L^k} \circ ev_p)$, and so on until the sections $s_{r_{m-1}+1},...,s_{r}\in \ker(\pi_{m-1,h}^{\mathscr{F}\otimes L^k} \circ ev_p)$ are orthogonal  to the kernel of $(\pi_{m,h}^{\mathscr{F}\otimes L^k} \circ ev_p)$. Eventually, we can do the assumption that $s_i(p)=0$ for $i>r$. Introduce the sets $\mathscr{I}_j=\{r_{j-1}+1,...,r_{j}\}$ and the map $f:i\mapsto j$ where $j$ is such that $i\in \mathscr{I}_j$. 
Hence, by (\ref{defFS}), we compute
\begin{eqnarray*}
|s_i(p)|_h&=& \frac{N}{Vr-V\sum_{j>0} \varepsilon_j r_j}, \\
|s_i(p)|_{FS_{k}}&=&  \frac{N}{Vr-V\sum_{j>0} \varepsilon_j r_j} -
 \frac{N}{Vr-V\sum_{j>0} \varepsilon_j r_j} \varepsilon_{f(i)}.
\end{eqnarray*}
Note that the term $\frac{s_i \otimes s_i^{*_{FS_{k}}}}{|s_i|^2_{FS_{k}}}\in End(\mathcal{F}\otimes L^k)$
is simply the orthogonal projection onto the image of $s_i$ respectively
to the metric $FS_{k}$ (and also $h$). At $p\in M$,
\begin{eqnarray}
\tsum_{i=1}^{N}s_{i}\left\langle .,s_{i}\right\rangle _{FS_{k}}
&=& \tsum_{i=1}^{N}s_{i}\left\langle .,s_{i}\right\rangle _{h} + \tsum_i \left(s_{i}\left\langle .,s_{i}\right\rangle _{FS_{k}} - s_{i}\left\langle .,s_{i}\right\rangle_{h}\right), \nonumber \\
&=& Cst \times Id - \tsum_i \left(\frac{s_i \otimes s_i^{*_{FS_{k}}}}{|s_i|^2_{FS_{k}}}
\frac{N}{Vr-V\tsum_j \varepsilon_j r_j} \varepsilon_{f(i)}\right), \nonumber \\
&=& Cst \times Id - \frac{N}{Vr-V\tsum_j \varepsilon_j r_j} \tsum_j \varepsilon_j \pi^{\mathscr{F}\otimes L^k}_{j,h_k}. \label{berg}
\end{eqnarray}
Here we have used the fact that for the quotient metric, the Bergman kernel
is constant, since it can be considered as the identity isomorphism
of the universal vector bundle over the Grassmannian.
Clearly, this implies the existence of a metric $h_k=FS_{k}$ on the bundle $\mathcal{F}\otimes L^k$ satisfying (\ref{BrH2c}).

Conversely, if (\ref{BrH2c}) is satisfied, as the $s_i$
are coordinate sections of the universal bundle, we obtain that they are also 
 $Hilb_{\omega}(FS_{k})$-orthonormal, i.e. they are orthonormal for
$$\int_M \Big\langle \frac{N}{Vr-V\sum_j \varepsilon_j r_j} \left( Id -
\displaystyle{\sum}_j \varepsilon_j \pi^{\mathscr{F}\otimes L^k}_{j,h}\right)  .,.\Big\rangle_h dV,$$
and therefore $Hilb_{\omega}(FS_{k})$ is a metric on $H^0(M,\mathcal{F}\otimes L^k)$ which is a zero of the moment map $\mu_{\mathscr{F},k}$. \qed
\end{proof}

\begin{theorem} \label{Br403a}
Let $\mathscr{F}$ be a holomorphic filtration of length $m$  over a projective manifold. Then
$\mathscr{F}$ is $\mathbf{R}$-Gieseker stable if and only if $Aut(\mathscr{F})=\mathbb{C}$ and for $k$ large, there exists a metric $h_k\in Met(\mathcal{F}\otimes L^k)$ such that
\begin{equation*}
\widehat{\mathsf{B}}_{\mathcal{F} \otimes L^k,h_{k}}+\epsilon_k\sum_{j=1}^{m-1} \varepsilon_j \pi^{\mathscr{F}\otimes L^k}_{j,h_k}= \frac{N+\epsilon_k \sum_{j=1}^{m-1} \varepsilon_j r_j}{rV} Id_{\mathcal{F}\otimes L^{k}}
\end{equation*}
where
$$\epsilon_k=\frac{\chi(\mathcal{F}\otimes L^k)}{Vr-V\sum_{j>0} \varepsilon_j r_j}.$$
\end{theorem}

\begin{proof}
One already knows by Theorem \ref{102} that the zeros of the moment map $\mu_{\mathscr{F},k}$ are the critical points of 
$\widetilde{KN_{k,\mathscr{F}}}$.  To apply Kempf-Ness stability criterion \cite{K-N}, we simply need to remark that the functionals $\widetilde{KN_{k,\mathscr{F}}}$ and $\widetilde{F_{k,\mathscr{F}}}$ are simultaneously proper by Lemma \ref{Br928a}. We are done with Lemma \ref{Br60aa}.\qed
\end{proof}

\section{$\boldsymbol{\tau}$-Hermite-Einstein metrics and holomorphic filtrations \label{section3}}

We now consider the Hermite-Einstein equation for a holomorphic filtration $\mathscr{F}$ of length $m$,
\begin{equation} \label{Bradlo}
 \sqrt{-1}\Lambda F_{h} +
\sum_{i=1}^{m-1} {\tau}_i \pi_{h,i}^{\mathscr{F}}= \mu_{\boldsymbol{\tau}}(\mathscr{F}) Id
\end{equation}
on a smooth projective manifold $M$ of complex dimension $n$. The goal of this part is to give an approximation of the metric solution of equation (\ref{Bradlo}) using the balanced metrics that we have just defined. We will need the following expansion proved in \cite{Ca,W2,Ze} of the Bergman kernel of $\mathcal{F}\otimes L^k$ when $k\rightarrow \infty$,

\begin{theorem-section} \label{ca}
Let $(M,\omega)$ be a K\"ahler manifold and $(L,h_L)$ an ample line bundle on  $M$
such that $\omega$ represents the curvature of $L$. Let $(\mathcal{F},h_\mathcal{F})$ be  a holomorphic hermitian vector bundle.
For any integer $\alpha\geq 0$, we have the following asymptotic expansion when $k\rightarrow +\infty$ of the Bergman kernel $\widehat{\mathsf{B}}_{h_{\mathcal{F}}\otimes h_{L^k}}$
\begin{equation}
\left\Vert \widehat{\mathsf{B}}_{h_{\mathcal{F}}\otimes h_{L^k}}-k^{n}Id_{r\times r}-\left( \frac{1}{2}%
Scal(g)Id_{r\times r}+\sqrt{-1}\Lambda F_{h_\mathcal{F}} \right)
k^{n-1}\right\Vert _{C^{\alpha}}\hspace{-0.4cm} \hspace{-0.05cm}\leq \hspace{-0.05cm} C_{\alpha}k^{n-2}  \label{As64}
\end{equation}%
where $Scal\left( g\right) $ denotes the scalar curvature of the metric $g$ associated
to the K\"ahler form
\begin{equation*}
\omega=\frac{i}{2}\sum g_{i\overline{j}}dz^{i}\wedge dz^{\overline{j}}.
\end{equation*}
This estimate is uniform on $M$ when $h_{\mathcal{F}}$ and $h_{L}$ belong to a 
 compact for the $C^{\alpha}$ topology.
\end{theorem-section}

Moreover, we denote in all the following for any integer $k$, $$\epsilon_k=\frac{\chi(\mathcal{F}\otimes L^k)}{Vr-V\sum_{j=1}^{m-1} \varepsilon_j r_j}.$$

\subsection{Action of the Gauge group \label{sectionpositif}}

\qquad Let $\mathcal{A}(\mathcal{F})$ be the space of  $C^{\infty }$ connections on the vector bundle $\mathcal{F}$. By the Newlander-Nirenberg's theorem, it is equivalent 
to considering a holomorphic structure on the $\mathbb{C}^{r}$ vector bundle $\mathcal{F}$  or an operator $\overline{\partial }$%
\begin{equation*}
\overline{\partial }:\Omega ^{0}(\mathcal{F})\rightarrow \Omega ^{0,1}(\mathcal{F})
\end{equation*}%
with $\overline{\partial }(f.s)=\overline{\partial }f.s+f.\overline{\partial
}s$ and $\overline{\partial }^{2}=0.$ We will denote    
$\mathcal{A}(\mathcal{F},h_{\mathcal{F}})$ the space of smooth connections which are compatible with the hermitian metric $h_{\mathcal{F}}$ (i.e. unitary). Any connection on the holomorphic vector bundle ${\mathcal{F}}$ which is integrable, i.e. belongs to the subset
\begin{equation*}
\mathcal{A}^{1,1}({\mathcal{F}},h_{\mathcal{F}})=\{A\in \mathcal{A}({\mathcal{F}},h_{{\mathcal{F}}}):F_{A}^{0,2}=F_{A}^{2,0}=0\}
\end{equation*}%
where $F_{A}\in \Omega ^{2}(M,End({\mathcal{F}}))$ denotes the curvature of the connection, defines a holomorphic structure on ${\mathcal{F}}.$ It is well known
that a holomorphic vector bundle equipped with a hermitian metric admits a
unique connection which is compatible with the holomorphic structure (i.e. integrable) and the hermitian structure $h_{{\mathcal{F}}}$, called Chern connection
 (see \cite{L-T} for the details). This means that there exists an isomorphism (see \cite[Section 8]{A-B}) $$A^{1,1}(\mathcal{F},h_\mathcal{F})\tilde{\rightarrow}\mathcal{C}(\mathcal{F}),$$  where we have set 
\begin{eqnarray*}
\mathcal{C}(\mathcal{F})=\{\nabla^{0,1}:\Omega^0(\mathcal{F})\rightarrow \Omega^{0,1}(\mathcal{F})
\text{ s.t. }&& \left(\nabla^{0,1}\right)^2=0, \\
&&\nabla^{0,1}(f.s)=\nabla^{0,1}f.s+f.\nabla^{0,1}s\}.
\end{eqnarray*}
$\mathcal{A}^{1,1}({\mathcal{F}},h_{\mathcal{F}})$ is a subvariety
 (possibly with singularities) of infinite dimension of the symplectic space  $\mathcal{A}({\mathcal{F}})$ which can be equipped of the symplectic structure
(cf. \cite[p.587]{A-B} or \cite{Do0}):
\begin{equation*}
\Omega(A,B)=\int_M Tr(A \wedge B)\frac{\omega^{n-1}}{(n-1)!}
\end{equation*}
 for any $A,B\in \Omega^1(End({\mathcal{F}}))$. One denotes by $\mathcal{G}$ the Gauge group of $\mathcal{F}$, i.e. the group of unitary automorphisms of ${\mathcal{F}}$:%
\begin{equation*}
\mathcal{G}=\{U\in C^{\infty }(GL(\mathcal{F})):\,U^*U=I\}.
\end{equation*}
The complexified Gauge group, i.e. the space of smooth sections of automorphisms $\mathcal{G}^{\mathbb{C}}=C^{\infty }(GL({\mathcal{F}}))$ acts (on the right) on $\mathcal{A}^{1,1}({\mathcal{F}},h_{\mathcal{F}})$ by
\begin{eqnarray}
\overline{\partial }_{g(A)} &=&g^{-1}\circ \overline{\partial }_{A}\circ g, \label{action1}\\
\partial _{g(A)} &=& ^{\mathsf{t}}\overline{g} \circ
\partial _{A}\circ \left(\hspace{0.1cm}^{\mathsf{t}}\overline{g}\right)^{-1}, \nonumber
\end{eqnarray}%
where $A\in \mathcal{A}^{1,1}(\mathcal{F},h_{\mathcal{F}})$, $g\in \mathcal{G}^{\mathbb{C}
}$ and $\partial _{A}$ is the $\left( 1,0\right) $ part of the covariant derivate induced naturally from $A$. In particular, the action of $\mathcal{G}^\mathbb{C}$ is holomorphic on $\mathcal{C}(\mathcal{F})$ which admits a complex structure. Therefore $\mathcal{A}^{1,1}({\mathcal{F}},h_\mathcal{F})$ 
inherits a complex structure for which $\Omega$ is a K\"ahler metric. 

\subsection{Limit of a sequence of balanced metrics \label{Sec2}}

We may naturally ask what can be the limit of a sequence of balanced metrics.
In this section, we show that if this limit does exist, it is necessarily a conformally  $\boldsymbol{\tau}$-Hermite-Einstein metric. 

\begin{definition}
Let $\left(\Xi,\omega \right) $ be a complex manifold and let $\mathscr{F}$ be a holomorphic filtration of length $m$ over $\Xi$. A hermitian metric $h$ on $\mathcal{F}$ is said to be conformally $\boldsymbol{\tau}$-Hermite-Einstein if the curvature $F_{h}$ of the Chern connection associated to  $h\in Met(\mathcal{F})$ satisfies the equation
\begin{equation} \label{faibspe}
\sqrt{-1}\Lambda _{\omega }F_{h}+ \sum_{i=1}^{m-1} {\tau}_i \pi_{h,i}^{\mathscr{F}}
=\lambda _{h}Id_{\mathcal{F}},
\end{equation}%
where $\lambda _{h}$ is a real valued function.
\end{definition}

\begin{proposition}
\label{40}Assume that $\left(\Xi,\omega \right) $ is a compact K\"ahler manifold. If $%
h$ is conformally $\boldsymbol{\tau}$-Hermite-Einstein then there exists $f\in C^{\infty
}\left( M,\mathbb{R}\right) $ unique up to a constant, such that the new metric  $e^{f}\cdot h$ is $\boldsymbol{\tau}$-Hermite-Einstein with parameter $\lambda =
\frac{1}{V}\tint_{M}\lambda _{h}dV$ in equation (\ref{faibspe}).
\end{proposition}

\begin{proof}
Indeed, with this conformal change, we get
$$ \sqrt{-1}\Lambda_{\omega} F_{e^f \cdot h}=\sqrt{-1}\Lambda_{\omega} F_h+\sqrt{-1}\Lambda_{\omega}\partial\overline{\partial}(f)Id $$
and $$\pi_{e^f \cdot h,i}^{\mathscr{F}}= \pi_{h,i}^{\mathscr{F}}.$$
Now, it can be found a function $f\in C^{\infty}(M,\mathbb{R})$ by classical theory of elliptic operators over a compact manifold $M$ (see \cite[Corollary 7.2.9]{L-T}) such that
$$\sqrt{-1}\Lambda_{\omega}\partial\overline{\partial}(f) = \lambda_h - \frac{1}{V}\int_M \lambda_{h} dV.$$\qed
\end{proof}

\begin{notation}
To a sequence of balanced metrics $h_k\in Met(\mathcal{F}\otimes L^k)$ for a holomorphic filtration $\mathscr{F}$, we shall associate the sequence of normalized metrics $\boldsymbol{h}_k=h_k \otimes h_L^{-k}\in Met(\mathcal{F})$ that will still be qualified as balanced metrics.
\end{notation} 

We use use the asymptotic expansion provided by Theorem \ref{ca} and the previous proposition to get 

\begin{theorem}
\label{Br38}Let $\mathscr{F}$ be a balanced holomorphic filtration of length $m$ over $M$. If the sequence of balanced metrics $\boldsymbol{h}_{k}\in Met(\mathcal{F})$ admits a limit $h_{\infty }$ in $C^2$ topology when  $k\rightarrow \infty ,$
then the metric $h_\infty$ is conformally $\boldsymbol{\tau}$-Hermite-Einstein
satisfying
\begin{equation}
\sqrt{-1}\Lambda F_{h_{\infty}} +\sum_{i=1}^{m-1} {\tau}_i \pi_{h,i}^{\mathscr{F}} =\left( \mu_{\boldsymbol{\tau}} \left( \mathcal{F}\right) +\frac{1}{2}\left(
\int_{M}c_{1}\left( M\right) \omega ^{n-1}-Scal(g)\right) \right)Id_{\mathcal{F}},
\label{Br105}
\end{equation}
and up to a conformal change, this metric is $\boldsymbol{\tau}$-Hermite-Einstein.
\end{theorem}

\subsection{Stable holomorphic filtrations and natural moment maps \label{Br44}}

We adapt the method of \cite{Do3} to make apparent the balanced condition
as the vanishing of two moment maps, one induced by the unitary group 
and the other one by the Gauge group of the vector bundle $\mathcal{F}$.
We know that the space $C^{\infty }(M,\mathcal{F})$ of smooth sections
 of $\mathcal{F}$ has a natural symplectic form $\Omega_{[0]} $ associated to the 
 hermitian metric $h_{\mathcal{F}}$ on $\mathcal{F}$:%
\begin{equation*}
\Omega_{[0]} (s_{1},s_{2})=2\Im\left( \int_{M}\langle s_{1},s_{2}\rangle
_{h_{\mathcal{F}}}dV\right) .
\end{equation*}%
It is not difficult to check that  
\begin{equation*}
\mu _{C^{\infty }(M,\mathcal{F})}(s)= \sqrt{-1} s\langle .,s\rangle _{h_{\mathcal{F}}}\frac{%
\omega ^{n}}{n!}\in \Omega ^{2n}(M,End(\mathcal{F}))\simeq Lie(\mathcal{G}^{\mathbb{C}}\mathcal{)}^{\ast }
\end{equation*}%
is a moment map associated to the action of the group $\mathcal{G}$ on $C^{\infty }(M,\mathcal{F}).$

For a holomorphic filtration $\mathscr{F}$, we can consider a family 
  $\theta_i:\mathcal{F}_i\rightarrow \mathcal{F}$ of smooth sections of the bundle in Grasmanians that we denote   $Gr(r_i,\mathcal{F})$ and whose fibers in $p\in M$ are the $r_i$ planes of $\mathcal{F}_{|p}$.
  Such a section $\theta_i$ gives naturally a projection $h_{\mathcal{F}}$-orthogonal onto the orthogonal of its kernel with respect to $h_{\mathcal{F}}$, i.e. the projection $\pi^{\mathscr{F}}_{h_{\mathcal{F}},i}.$ Moreover, the metric $h_{\mathcal{F}}$ on the fiber  $\mathcal{F}_{|p}$ induces a K\"ahler form $\omega^{Gr}_{h_{\mathcal{F}}}$ on $Gr(r_i,\mathcal{F}_{|p})$. We obtain a symplectic form using the evaluation map $ev_i:C^{\infty}(M,Gr(r_i,\mathcal{F}))\times M \rightarrow Gr(r_i,r)$ and the projection on the first component $p_1:C^{\infty}(M,Gr(r_i,\mathcal{F}))\times M \rightarrow C^{\infty}(M,Gr(r_i,\mathcal{F})),$
  $$\Omega_{(i)}=(p_1)_*\left(ev_i^*(\omega^{Gr}_{h_{\mathcal{F}}}) \wedge dV\right).$$
  The action of the Gauge group on $C^{\infty}(M,Gr(r_i,\mathcal{F}))$ respectively to $\Omega_{(i)}$ is 
  then given by  $$\mu_{C^{\infty}(M,Gr(r_i,\mathcal{F}))}(\theta_i)=\sqrt{-1}\pi^{\mathscr{F}}_{h_{\mathcal{F}},i}\frac{\omega^n}{n!}\in Lie(\mathcal{G}^{\mathbb{C}})^*.$$
  
  We know that there exists for $k$ large enough an embedding $i_{k}$ of $M$ into the
Grassmannian $Gr(N,r)$ using the holomorphic sections of $\mathcal{F}\otimes L^{k}.$
Nevertheless  the action  of $\mathcal{G}$ on $C^{\infty }(M,\mathcal{F}\otimes L^{k})$ does not preserve the
set of holomorphic sections for a connection  $A \in
\mathcal{A}^{1,1}(\mathcal{F}\otimes L^{k},h_{\mathcal{F}}\otimes h^{k}_L)$ defined \textit{a priori}. Note that in general, the dimension of the space of holomorphic sections of of $\mathcal{F}\otimes L^k$ depends on the choice of the connection $A$. \\
  In order to consider global holomorphic sections and their variations with respect to the Gauge group, 
  we are constrained to modify simultaneously the considered connection. But for any connection $A$ and for  $k$ large enough, there exists an open set of the complex orbit of $A$ in $\mathcal{G}^{\mathbb{C}}$ such that for any connection belonging to this set,  $\dim(H^i(M,\mathcal{F}\otimes L^k))=0$ (by semi-continuity \cite[Section 9.3]{V}) and $\dim(H^0(M,\mathcal{F}\otimes L^k))$ is constant. Finally for such a $k$, we need to introduce the following   manifold of infinite dimension,

\begin{definition}
Let $\mathscr{F}$ be a holomorphic filtration of length $m$, and $\mathcal{Q}_{0}$ the subset of
$$C^{\infty }(M,\mathcal{F}\otimes L^{k})^{N}\times \mathcal{A}^{1,1}(\mathcal{F},h_{\mathcal{F}})\times \prod_{i=1}^{m-1} C^{\infty}(M,Gr(r_i,\mathcal{F}))$$  formed by $(N+m)$-tuples of the form
\begin{equation*}
\Big\{s_{1},...,s_{N},A,\theta_1,...,\theta_{m-1} \Big\}
\end{equation*}
 such that the sections $\left( s_{i}\right) _{i=1..N}$ are linearly independent,  and
\begin{eqnarray}
\overline{\partial }_{A}s_{i} &=&0\quad \text{ }\forall i=1,..,n  \label{Br58} \\
\overline{{\partial}}_{[A]}\theta_j&=&0 \quad \text{ }\forall j=1,..,m-1, \label{Br58f}
\end{eqnarray}%
where $\overline{\partial }_{A}$ represents the $\left(0,1\right) $ part of the covariant derivative induced naturally using the unitary connection $A$ and the Chern connection on $L$
over the space $C^{\infty}(M,\mathcal{F}\otimes L^k)$ and $\overline{{\partial}}_{[A]}$ represents the $\left(0,1\right) $ part of the covariant derivative induced naturally using $A$ on $C^\infty(M,Gr(r_i,\mathcal{F}))$.
\end{definition}

Under these conditions, the diagonal action of $\mathcal{G}$ preserves $
\mathcal{Q}_{0}$. Also notice that the unitary group $U(N)$
acts naturally on  $\mathcal{Q}_{0}$ by
\begin{equation*}
\left( \mathsf{u}_{ij}\right)\ast \{s_{1},...,s_{N},A,\theta_1,...,\theta_{m-1}\}=\left\{
\tsum_{j=1}^{N}\mathsf{u}_{1j}s_{j},...,\tsum_{j=1}^{N}\mathsf{u}_{Nj}s_{j},A,\theta_1,...,\theta_{m-1}\right\},
\end{equation*}%
and if we denote $h=h_{\mathcal{F}}\otimes h_{L^k}$ then the moment map for this action is,
\begin{equation*}
\mu _{U(N)}(s_{1},...,s_{N},A,\theta_1,...,\theta_{m-1})=\int_M \langle s_{i},s_{j}\rangle_{h}dV.
\end{equation*}%
Let $\pi
_{0}:\mathcal{Q}_{0}\rightarrow C^{\infty }(M,\mathcal{F}\otimes L^{k})^{N}\times \prod_{i=1}^{m-1} C^{\infty}(M,Gr(r_i,\mathcal{F}))$
be the natural  projection. Actually, $\pi _{0}$ is an immersion.
 Indeed, one has by (\ref{Br58}) that $0=\overline{%
\partial }_{A}\left( \varepsilon s_{i}\right) +\varepsilon \overline{%
\partial }_{A}\left( s_{i}\right) =\varepsilon \overline{\partial }%
_{A}\left( s_{i}\right) $ (and similary  $\varepsilon \overline{\partial }%
_{A}\left( \theta_i\right)=0$) for a small variation of $\left(s_1,..,s_N, \overline{
\partial}_{A},\theta_1,..,\theta_{m-1}\right) $. Since $i_{k}$ is an embedding,
 we get that $\varepsilon \overline{\partial }_{A}=0$ and hence
$d\pi _{0}$ is injective, which gives the immersion property.

We consider the symplectic standard form 
$\Omega_{[k]}$ on $C^{\infty }(M,\mathcal{F}\otimes L^k)$. We take the sum
of $\Omega_{[k]}$ over $N$ copies of the spaces $C^{\infty }(M,\mathcal{F}\otimes L^k)$ and consider the symplectic form $\Omega_{(i)}$ on each $C^{\infty}(M,Gr(r_i,\mathcal{F}))$ with weight $\epsilon_k\varepsilon_i$. We get a symplectic form that we can lift using
the injective immersion $\pi _{0}$. We denote $\Omega
_{\mathcal{Q}_{0}}$ the symplectic form obtained on ${\mathcal{Q}_{0}}$.  The action of $\mathcal{G}$ over $\mathcal{Q}_{0}$
admits a moment map associated to $\Omega _{\mathcal{Q}_{0}},$%
\begin{eqnarray*}
\mu _{\mathcal{G}}(s_{1},...,s_{N},A,\theta_1,...,\theta_{m-1})&=&\sum_{i=1}^{N}s_{i}\left\langle .,s_{i}\right\rangle_{h} + \epsilon_k
\sum_{i=1}^{m-1}\varepsilon_i\pi^{\mathscr{F}\otimes L^k}_{h,i}.
\end{eqnarray*}%
Indeed, $\mathcal{Q}_{0}$ admits a K\"ahler structure by \cite[Theorem 3]{Hi} and with the fact $\mathcal{A}^{1,1}(\mathcal{F},h_\mathcal{F})$ admits a complex structure. Moreover, the actions of $\mathcal{G}$ and $U(N)$ commute and these two groups have center of dimension
$1$, given respectively by the constant functions and the multiples of identity.
This allows us to restrict to the action $%
SU(N)$ by considering another natural moment map with values in the Lie algebra $\sqrt{-1}\mathfrak{su}(N)$,
\begin{equation*}
\mu _{SU(N)}(s_{1},..,s_{N},A,\theta_1,..,\theta_{m-1})= \int_M \langle s_{i},s_{j}\rangle_h dV-\frac{1}{N}\left( \tsum_{i} \int_M \vert s_{i}\vert_h^2 dV
\right) \delta _{ij}.
\end{equation*}

\subsection{Complex orbits and double symplectic quotient \label{Sec1}}

The moment map for the action of the product $\mathcal{G}\times
SU\left( N\right) $ on $\mathcal{Q}_{0}$ will be given by the direct sum $\mu _{\mathcal{G%
}}\oplus \mu _{SU(N)}$. For a real number $\lambda,$  we will consider the symplectic quotient 
\begin{equation*}
\mathcal{Q}_{0}\boldsymbol{//}\left( \mathcal{G}\times SU\left( N\right) \right) =\frac{%
\mu _{\mathcal{G}}^{-1}\left( \lambda Id\right) \cap \mu _{SU(N)}^{-1}\left(
0\right) }{\mathcal{G}\times SU\left( N\right) }.
\end{equation*}
This must be understood as the symplectic quotient of $\mathcal{Q}_{0}$ by $\mathcal{G}$ in a first step, via
\begin{equation*}
\mathcal{Q}_{0}\boldsymbol{//}\mathcal{G}= \mu _{\mathcal{G}}^{-1}\left( \lambda
Id\right) /\mathcal{G}
\end{equation*}%
which admits as we have seen, a natural symplectic structure. Secondly, we do the symplectic quotient of $%
\mathcal{Q}_{0}\boldsymbol{//}\mathcal{G}$ by $SU\left( N\right) $, since
$SU\left( N\right) $ acts naturally on
$\mathcal{Q}_{0}\boldsymbol{//}\mathcal{G}$. Every $\mathcal{G}^{\mathbb{C}}$-orbit in $\mathcal{Q}_{0}$ contains a point in $\mu
_{\mathcal{G}}^{-1}\left( \lambda Id\right) $, (resp. $\mu
_{SU(N)}^{-1}(0))$ unique up to the action of $\mathcal{G}$
(resp.\ $SU\left( N\right) )$. Our situation is summed up by the following proposition.

\begin{proposition}
There exists a metric $k$-balanced for $\mathscr{F}$ if and only if the complex orbit in $\mathcal{Q}_{0}$ given by the action of $\mathcal{G}^{\mathbb{C}}\times SL\left( N\right)$, contains a point in   $\mu _{\mathcal{G}}^{-1}\left( \lambda Id\right) \cap
\mu _{SU(N)}^{-1}(0)$ for all $\lambda >0.$ This is equivalent to saying that the complex orbit is represented by a point in the double symplectic quotient $\mathcal{Q}_{0}\boldsymbol{//}\left( \mathcal{G}\times SU\left( N\right) \right)$.
\end{proposition}

\begin{proof}
A point $z_{0}=\{s_{1},...,s_{N},A,\theta_1,...,\theta_{m-1}\}\in \mathcal{Q}_{0}$
belongs to $\mu _{\mathcal{G}}^{-1}\left( \lambda Id\right) $ if and only if 
$ \sum_{i}s_{i}\langle
.,s_{i}\rangle_h +\epsilon_k \sum_i \varepsilon_i \pi^{\mathscr{F}}_{h,i}
=\lambda Id $ and $z_{0}$ belongs to $\mu _{SU(N)}^{-1}(0)$ if and only if there exits a 
constant $c$ such that the sections 
$s_{i}$ form a $L^{2}$-orthonormal basis. Therefore we get $\lambda=\frac{1}{rV}\left(\chi(\mathcal{F}\otimes L^k)+\epsilon_k \sum_{j} \varepsilon_j r_j\right)$ by taking the trace and hence we obtain the balance condition of Lemma \ref{Br60aa}.
\end{proof}

Consider the action $\cdot $ of $SU(N)$ on the symplectic quotient
\begin{equation*}
\mathcal{Z}=\mathcal{Q}_{0}\boldsymbol{//}\mathcal{G}.
\end{equation*}%
Since $\mathcal{Q}_0$ is K\"ahler, $\mathcal{Z}$ admits a K\"ahler structure (in fact it is an orbifold
if  the stabilizers group is finite everywhere) since the extension of the action 
of $\mathcal{G}$ to $\mathcal{G}^{\mathbb{C}}$ preserves the complex structure. For all point $%
z\in \mathcal{Z},$ the infinitesimal action of $\mathfrak{su}%
(N),$ provides us with an application
\begin{equation*}
\nu _{z}^{\mathcal{Z},SU\left( N\right) }:\mathfrak{su}(N)\rightarrow T%
\mathcal{Z}_{z}.
\end{equation*}%
Let us consider the following operator on $\mathfrak{su}(N)$,
\begin{equation*}
\mathfrak{q}_{z}^{SU\left( N\right) }=\left( \nu _{z}^{\mathcal{Z},SU\left( N\right) }\right) ^{\ast }\nu _{z}^{\mathcal{Z},SU\left( N\right)}, 
\end{equation*}
where $\left( \nu _{z}^{\mathcal{Z},SU\left( N\right) }\right) ^{\ast }$ is the adjoint of $\nu _{z}^{\mathcal{Z},SU\left( N\right) }$ formed using the invariant metric
 on  $\mathfrak{su}(N)$ and the metric on $T\mathcal{Z}_{z}.$
Assume that the stabilizers of a point in $\mathcal{Z}$ under the action
of $SU(N)$ are discrete; then $\nu _{z}^{\mathcal{Z},SU\left( N\right) }$
is injective and $\mathfrak{q}_{z}^{SU\left( N\right) }$ is invertible.

\begin{notation} 
Let $\mathsf{Q}$ be a hermitian matrix. The Hilbert-Schmidt norm and the operator norm for $\mathsf{Q}$ are given by
\begin{eqnarray*}
||\mathsf{Q}||^{2} &=&\sum_{i ,j }|\mathsf{Q}_{i j }|^{2},
\\
\left\vert \left\vert \left\vert \mathsf{Q}\right\vert \right\vert
\right\vert &=&\underset{||v||\leq 1}{\sup }\frac{|\mathsf{Q}v|}{|v|}.
\end{eqnarray*}%
\end{notation}

\begin{notation}
Let $\mathit{\Lambda }_{z}$ (resp. $\boldsymbol{\Lambda }_{z}$) be the Hilbert-Schmidt norm (resp. operator norm) of $\left(
\mathfrak{q}_{z}^{SU\left( N\right) }\right) ^{-1}:\mathfrak{su}(N)\rightarrow \mathfrak{su}(N)$ with respect to the invariant euclidian metric on $\mathfrak{su}(N)$.
\end{notation}
In particular, the inequality $\boldsymbol{\Lambda }_{z}\leq \lambda$ is induced from the fact that for all $\mathsf{A}\in \sqrt{-1}\mathfrak{su}(N)$ one has
$|\mathsf{A}|^2\leq \lambda \left\vert \nu_{z}^{\mathcal{Z},SU(N)}(\sqrt{-1}\mathsf{A}) \right\vert^2_{T\mathcal{Z}}$.\newline

We will need the following key result inspired directly from the formalism for moment maps from \cite{Do3}.
\begin{proposition}\label{Br16}
Let $z_{0}\in \mathcal{Z}$ and $\lambda
,\delta >0$ be two real numbers such that \\
1. $\lambda ||\mu_{SU(N)}(z_{0})||<\delta$; \\
2. $\boldsymbol{\Lambda}_{z}\leq \lambda$ pour tout $z=e^{i\mathsf{S}}\cdot z_{0}$ and $||\mathsf{S}||\leq \delta $.\\
Then, there exists a zero $z_{1}=e^{i\mathsf{S}^{\prime }}\cdot z_{0}$ of $\mu_{SU(N)}$,
\begin{equation*}
\mu_{SU(N)}(z_{1})=0
\end{equation*}%
with $||\mathsf{S}^{\prime }||\leq \lambda ||\mu(z_{0})||.$ Here $||.||$ is the norm induced
by the standard $SU(N)$-invariant inner product on $\mathfrak{su}(N)$. 
\end{proposition}

\begin{proof}
See \cite[Proposition 17]{Do3}. Indeed, let $z_0$ be a point in $\mathcal{Z}$ and  denote $(z_t)_{t>0}$ the trajectory of $z_0$ along the flow $-\overrightarrow{Grad}(||\mu_{SU(N)}(z)||^2)$ and set $\mathcal{Z}^{min}=\{ z_0\in \mathcal{Z}: \lim_{t\rightarrow +\infty}z_t \in\mu^{-1}(0)\}$. Then it is well-known (see \cite[Theorem 7.4]{Ki2}) that the complex orbits of the zeros of $\mu_{SU(N)}$ can be identified with the set $\mathcal{Z}^{min}$. Finally, this proposition is an effective version of this result.\qed
\end{proof}

The aim of the next three sections is the approximation Theorem \ref{Br39} of a metric solution  of the $\boldsymbol{\tau}$-Hermite-Einstein equation using the tools defined in sections \ref{sectionpositif}, \ref{Sec2}, \ref{Br44}, \ref{Sec1}. We are now ready to give the main ideas to prove this result. 

An irreducible $\boldsymbol{\tau}$-Hermite-Einstein holomorphic filtration is in particular Gieseker $\mathbf{R}$-stable when one has fixed $R_i=\tau_i k^{n-1}$. Now, we are seeking to construct for such a  filtration a sequence of balanced metrics which converge towards the (conformally) $\tau$-Hermite-Einstein metric that satisfies equation (\ref{Br105}). To a $k$-balanced filtration $\mathscr{F}$, we know that there corresponds
a point in a certain space of parameters, inside a complex orbit (under the action of $\mathcal{G}^{\mathbb{C}}\times SL\left( N\right) )$ and this point is a zero of the moment map $\mu _{\mathcal{G}}\oplus \mu _{SU(N)}$ as previously defined. Finally this leads us to 
look for this point.
\begin{itemize}
\item From one side, via Propositions \ref{Br75} and  \ref{obtenirpe}, we get a point in the symplectic quotient by $\mathcal{G}$, i.e. a zero of $\mu_{\mathcal{G}}$: the existence \textit{a priori} of a conformally  $\boldsymbol{\tau}$-Hermite-Einstein metric will permit us to construct an \lq almost' balanced metric
 $h_{k,q}$ that will provide us with a point in $\mathcal{Z}$. To such a point will correspond a metric denoted $\widetilde{h_q}$.
\item From another side, we will switch in a finite dimensional problem by looking for a
zero (unique up to action of $SU(N)$) of the moment map
$\mu_{SU(N)}$ in a $SL(N)$-orbit. Thus, we will study the gradient flow of 
 $||\mu _{SU(N)}||^{2}$ on $\mathcal{Z}$ and give an estimate of 
$\boldsymbol{\Lambda }_{z}$ in order to apply Proposition
\ref{Br16}. 
\end{itemize}
Therefore, we finally obtain a $k$-balanced metric and also by construction,
 the convergence of this sequence of metrics (when $k\rightarrow \infty$) towards the conformally $\boldsymbol{\tau}$-Hermite-Einstein metric solution of (\ref{Br105}).

\subsection{Construction of almost balanced metrics}

 From now on and until the end of this section, we will consider that there exists a conformally $\boldsymbol{\tau}$-Hermite-Einstein metric  $h_{\infty }$ for the irreducible filtration $\mathcal{F}$. 
 By singular perturbation of the conformally $\boldsymbol{\tau}$-Hermite-Einstein metric, we can get a metric almost balanced, as will be made explicit in the following proposition.

\begin{proposition} \label{Br75}
Let $\mathscr{F}$ be an irreducible holomorphic filtration over $M$ of length $m$ such there exists a  conformally $\boldsymbol{\tau}$-Hermite-Einstein metric $h_\infty$ on $\mathcal{F}$ satisfying equation (\ref{Br105}). 
Then, there exists a family of smooth hermitian endomorphisms
 $\left(\mathsf{\boldsymbol{\eta }}_{i}\right) _{i\in\mathbb{N}}\in C^{\infty }(End(\mathcal{F}))$ such that the metrics defined on $\mathcal{F}$ for all  $q\geq 1$ by
\begin{equation*}
\left\langle .,.\right\rangle _{h_{k,q}}=\left\langle \left(
Id+\sum_{i=1}^{q}\frac{1}{k^{i}}\boldsymbol{\eta }_{i}\right) .,.\right\rangle
_{h_{\infty }}\otimes {h_{L^k}},
\end{equation*}%
are hermitian and smooth for $k$ large enough and there exists a constant $C_{q,\alpha}$ such that
\begin{equation}
\widehat{\mathsf{B}}_{\mathcal{F}\otimes L^k,h_{k,q}}+
\epsilon_k\sum_{i=1}^{m-1} \varepsilon_i\pi^{\mathscr{F}\otimes L^k}_{h_{k,q},i} =\frac{\chi (\mathcal{F}\otimes L^{k})+\epsilon_k
\sum_i \varepsilon_i r_i}{rV}Id+\mathsf{\boldsymbol{\sigma }}_{q}(k),
\label{B79}
\end{equation}%
where $\Vert\boldsymbol{\sigma }_{q}(k)\Vert_{C^{\alpha+2}}\leq
C_{q,\alpha}k^{n-q-1}$. 
The metrics $h_{k,q}$ will be called \lq almost balanced'. Here $C_{q,\alpha}$ is a constant 
that depends only on $q,\alpha,h_{\infty }$ and $\omega .$
\end{proposition}

\begin{proof}
First of all, we know that under these assumptions, $\mathscr{F}$ is simple by Lemma \ref{lemme1}.
The  Theorem \ref{ca} asserts that we have an asymptotic expansion in the variable $k$
of the Bergman kernel
\begin{equation*}
\widehat{\mathsf{B}}_{\mathcal{F}\otimes L^k,h_{\infty}\otimes h_{L^k}}
=k^{n}Id+%
\mathsf{a}_{1}(h_{\infty })k^{n-1}+...+\mathsf{a}_{q}(h_{\infty })k^{n-q}+%
\mathbf{O}(k^{n-q-1}),
\end{equation*}%
and the $\mathsf{a}_{i}$ are polynomial expressions of the curvature tensors of $%
h_{\infty }$ and $h_{L}$ and their covariant derivatives. The approximation term 
is uniformly bounded in $C^{\alpha}$ norm when ${h_{\infty }\otimes h_{L}^{k}}$ belongs to a  bounded family in $C^{\alpha^{\prime }}$
norm  (where $\alpha^{\prime }$ depends on $\alpha).$ 
We notice that we have also $\mathsf{a}_1(h_\infty)=\sqrt{-1}\Lambda_{\omega}F_{h_{\infty}}$. \\
Consequently, $%
\mathsf{a}_{i}\left( h_{\infty }\left( 1+\mathsf{\boldsymbol{\eta }}\right)
\right) =\mathsf{a}_{i}(h_{\infty })+\sum_{l=1}^{q}\mathsf{a}_{i,l}(\mathsf{%
\boldsymbol{\eta }})+O(\left\Vert \mathsf{\boldsymbol{\eta }}\right\Vert
_{C^{s}}^{q+1})$ with $s$ sufficiently large depending on $\alpha$ and $q.$
For all $\left( \mathsf{\boldsymbol{\eta }}_{i}\right) _{i\in\mathbb{N}}\in C^{\infty }(End(\mathcal{F})),$ we can write 
\begin{equation*}
\mathsf{a}_{i}\left( h_{\infty }\left( 1+{\tsum}_{j=1}^{q}\mathsf{\boldsymbol{\eta }%
}_{j}k^{-j}\right) \right) =\mathsf{a}_{i}(h_{\infty })+\tsum_{l=1}^{q}%
\mathsf{b}_{i,l}k^{-l}+\mathbf{O}(k^{-q-1}),
\end{equation*}%
where the $\mathsf{b}_{i,l}$ are multilinear expressions in $%
\mathsf{\boldsymbol{\eta }}_{j}$ and their covariant derivatives, beginning by
\begin{equation*}
\mathsf{b}_{i,1}=\mathsf{a}_{i,1}\left( \mathsf{\boldsymbol{\eta }}\right).
\end{equation*}%
If we now set $\mathsf{a}_{i}=\mathsf{a}_{i}(h_{\infty }),$
then we get
\begin{eqnarray}
\widehat{\mathsf{B}}_{\mathcal{F}\otimes L^k,h_{k,q}} &=&\tsum_{p=0}^{q}k^{n-p}%
\mathsf{a}_{p}+\tsum_{i,l=1}^{r}\mathsf{b}_{i,l}k^{n-i-l}+\mathbf{O}%
(k^{n-q-1}),  \notag \\
&=&k^{n}+\mathsf{a}_{1}k^{n-1}+\left( \mathsf{a}_{2}+\mathsf{b}_{1,1}\right)
k^{n-2}  \notag \\
&&+k^{n-3}\left( \mathsf{a}_{3}+\mathsf{b}_{1,2}+\mathsf{b}_{2,1}\right)
+...+\mathbf{O}(k^{n-q-1}),\qquad \qquad  \label{Br98}
\end{eqnarray}%
and we choose inductively the $\mathsf{\boldsymbol{\eta }}_{j}$ in such a way that
the coefficients that appear with $k^{n-j}$ $(j<q)$
are constants, which means that the RHS of (\ref{Br98}) is exactly
 (up to order $k^{n-q-1}$)
\begin{equation*}
\frac{\chi (\mathcal{F}\otimes L^{k})+\epsilon_k
\sum_i \varepsilon_i r_i}{rV}Id - \epsilon_k\sum_i \varepsilon_i\pi^{\mathscr{F}\otimes L^k}_{h_{k,q},i}.
\end{equation*}
We set the asymptotic expansions (in the variable $k)$,
\begin{eqnarray*}
\frac{\chi (\mathcal{F}\otimes L^{k})+\epsilon_k
\sum_i \varepsilon_i r_i}{rV}Id &=&\mathsf{c}_0k^n+\mathsf{c}_1k^{n-1}+... \\
\epsilon_k \frac{1}{k} &=&d_1k^{n-1}+d_2k^{n-2}+...\\
\end{eqnarray*} 
By Lemma \ref{annex1}, we get also the expansion
\begin{eqnarray*}
\tsum_j \tau_j \pi^{\mathcal{F}\otimes L^k}_{h_{k,q},j}&=&\tsum_j \tau_j \pi^{\mathcal{F}}_{h_{\infty},j}+ \frac{1}{k}
 \Pi_{h_{\infty}}^{\mathscr{F},\boldsymbol{\tau}}(\boldsymbol{\eta}_1)+...\\
&=& \mathsf{e}_1+k^{-1}\mathsf{e}_2+...
\end{eqnarray*}
where we have done the substitution $\epsilon_k \sum_j \varepsilon_j \pi^{\mathcal{F}\otimes L^k}_{h_{k,q},j}=
 \frac{\epsilon_k}{k} \sum_j \tau_j \pi^{\mathcal{F}\otimes L^k}_{h_{k,q},j}$. In particular, we have $d_1=1$ by our choice of $\epsilon_k$. \\
 From another point of view, we know that if $F_{H}$ denotes the curvature of the $H,$
then $F_{H\left( 1+\varepsilon \right) }=F_{H}+\overline{\partial }\partial
\varepsilon +O\left( ||\varepsilon ||^{2}\right),$ and therefore 
\begin{equation*}
\mathsf{b}_{1,1}=\sqrt{-1}\Lambda \left( \overline{\partial }\partial \mathsf{\boldsymbol{%
\eta }}_{1}\right).
\end{equation*}%
Moreover, when $k$ is sufficiently large
\begin{eqnarray*}
e_1 &=& \tsum_j \tau_j \pi^{\mathcal{F}}_{h_{\infty},j}, \\
e_2 &=& \tsum_j \tau_j \pi^{\mathcal{F}}_{h_{\infty},j}\left(\boldsymbol{\eta}_1\right)(Id-\pi^{\mathcal{F}}_{h_{\infty},j}).
\end{eqnarray*}
In order to get $\boldsymbol{\eta}_1$, we aim to solve
$$\mathsf{b}_{1,1}+d_1\mathsf{e}_2 =\mathsf{c}_{2}-\mathsf{a}_{2}-d_2\mathsf{e}_1.$$
But the operator $Q:u \mapsto \sqrt{-1}\Lambda \overline{\partial }\partial u +d_1  \Pi^{\mathscr{F},\boldsymbol{\tau}}_h(u)$
is elliptic of order 2. We can apply Lemma \ref{Br36e} noticing that $\int_M \tr( \mathsf{c}_{2}-\mathsf{a}_{2}-d_2\mathsf{e}_1)=0$
and the endomorphism $\left(\mathsf{c}_{2}-\mathsf{a}_{2}-d_2\mathsf{e}_1\right)$ preserves the filtration. Then, we get a solution $\boldsymbol{\eta}_1$
which is self-adjoint since the term $\mathsf{c}_2-\mathsf{a}_{2}-d_2\mathsf{e}_1 $ is also self-adjoint. 
Now, if we are looking for an almost balanced metric up to order $3$, it is sufficient to solve
$$ \mathsf{b}_{2,1} + d_1\mathsf{e}_3= \mathsf{c}_3-\mathsf{a}_3-\mathsf{b}_{1,2}-d_3 \mathsf{e}_1 -d_2\mathsf{e}_2.$$
We find all the $\mathsf{\boldsymbol{\eta }}_{j}$ by solving at each step a 
differential equation of the form
\begin{equation*}
\sqrt{-1}\Lambda \left( \overline{\partial }\partial \mathsf{\boldsymbol{\eta }}_j%
\right)+ \tsum_j \tau_j \pi^{\mathcal{F}}_{h_{\infty},j}\boldsymbol{\eta}_j(Id-\pi^{\mathcal{F}}_{h_{\infty},j}) =\mathsf{c}_{j}- \mathsf{a}_{j} - P_j(\boldsymbol{\eta}_1,\boldsymbol{\eta}_2,...,\boldsymbol{\eta}_{j-1}),
\end{equation*}%
where $P_j$ is self-adjoint, $\int_M \mathtt{tr}(P_j+\mathsf{a}_j-\mathsf{c}_j)=0$, $(P_j+\mathsf{a}_j-\mathsf{c}_j)$ preserves  the filtration, and $P_j$ is totally
determined by the $\boldsymbol{\eta}_l$ computed previously for $l<j$.
The fact that $h_{k,q}$ is hermitian is clear since  the endomorphisms $\mathsf{a}_{i}$, 
the generalized Bergman kernel and the operator $P_j$ are hermitian.\qed
\end{proof}

In order to get from the conformally $\boldsymbol{\tau}$-Hermite-Einstein metric a point
in the symplectic quotient $\mathcal{Z}$, we will need the following lemma.

\begin{lemma} \label{avantpe}
Let
$$\mathfrak{q}=(s_1,...,s_N,A,\theta_1,...,\theta_{m-1})\in \mathcal{Q}_0$$ 
be a point corresponding to the holomorphic filtration $\mathscr{F}$. Fix $\widetilde{h}$ a hermitian metric on $\mathcal{F}\otimes L^k$ and consider the map in the Sobolev space $End(\mathcal{F})^{l,\alpha}$ of hermitian endomorphisms of $C^{l,\alpha}$ class of  $\mathcal{\mathcal{F}}$ over the compact manifold $M$,
$$ \mathcal{B}_{\mathfrak{q},\widetilde{h},\rho}:\eta \mapsto  \left(\tsum_{i=1}^{N} s_{i }\langle .,s_{i }\rangle_{\widetilde{h}}\right)\eta +\rho\tsum_{j=1}^{m-1} \frac{\tau_j}{k}  \pi^{\mathscr{F}\otimes L^k}_{\widetilde{h}(\eta.,),j}$$
where $\rho \in \mathbb{R}$. If we assume that \\
--  $\tr \left( \tsum_{i }s_{i }\langle .,s_{i }\rangle_{\widetilde{h}}\right)=ck^n +O(k^{n-1})$, \\
--  $0\leq \rho_0 \leq \epsilon_k=\frac{\chi(\mathcal{F}\otimes L^k)}{Vr-V\sum_{j>0} \varepsilon_{j} r_j}$, \\
--  $k\geq k_0$ where $k_0$ only depends on the choice of the data  $\{(\tau_i,r_i)_{i=1,..,m},c\}$, \\
--  $\Gamma$ is an $\widetilde{h}$-hermitian smooth endomorphism such that $\tr(\Gamma)=O(k^n)$, \\
then there exists for all $0 \leq \rho \leq \rho_0$ a smooth solution $\eta_{\rho}$ of
\begin{equation}
\mathcal{B}_{\widetilde{h},\rho}(\eta_{\rho})=\Gamma. \label{leqn}
\end{equation}
\end{lemma}
\begin{proof}
We will use a continuity method on the Banach space $End(\mathcal{F})^{l,\alpha}$ with 
respect to the parameter $\rho$. First of all, for $\rho=0$, one sees that it is possible to solve (\ref{leqn}) by choosing 
$$\eta_0=\left(\tsum_{i }s_{i }\langle .,s_{i }\rangle_{\widetilde{h}}\right)^{-1}\Gamma. $$
Therefore, if one denotes $I\subset \mathbb{R}_{+}$ the interval such that $\rho \in I$ if $\eta_{\rho}$
is solution of (\ref{leqn}), we have just proved that $I \neq \emptyset$. Apply the Implicit function Theorem to see that  $I$ is open. By Lemma \ref{annex1}, we know that the differential in $\eta$ of $\mathcal{B}_{\mathfrak{q},\widetilde{h},\rho}$ 
is given by,
$$D_{\eta}\mathcal{B}_{\mathfrak{q},\widetilde{h},\rho}(\varsigma)=  \left(\tsum_{i }s_{i }\langle .,s_{i }\rangle_{\widetilde{h}}\right)\varsigma
+ \frac{\rho}{k}{\Pi}^{\mathscr{F}\otimes L^k,\boldsymbol{\tau}}_{\widetilde{h}\cdot \eta}(\varsigma ). $$
But, from another side, we know that
$$\normop{\Pi^{\mathscr{F}\otimes L^k,\boldsymbol{\tau}}_{\widetilde{h}\cdot \eta}} \leq
\tsum_i \tau_i r_i(r-r_i).$$ Now, with our choice of $\rho_0$ we get that  $\mathcal{B}_{\mathfrak{q},\widetilde{h},\rho}$ is invertible since we have $\tr( \sum_{i }s_{i }\langle .,s_{i }\rangle_{\widetilde{h}})=O(k^n)$. \newline

Finally, we prove that $I$ is closed. If one has a solution $\eta$ of  
\begin{equation}
\mathcal{B}_{\mathfrak{q},\widetilde{h},\rho}(\eta)=\Gamma, \label{eqndiff}
 \end{equation}
 then for all $U\in \mathcal{F}\otimes L^k_{|p}$,
\begin{equation}
\tsum_{i }\langle s_{i },U \rangle_{\widetilde{h}\cdot \eta} \langle U,s_{i }\rangle_{\widetilde{h}\cdot \eta} +\rho \tsum_j  \frac{\tau_j}{k} \langle U, \pi^{\mathscr{F}\otimes L^k}_{\widetilde{h}\cdot\eta,j}U \rangle_{\widetilde{h}\cdot \eta}= \langle U,\Gamma U\rangle_{\widetilde{h}\cdot \eta}. \label{2ineq}
\end{equation}
Since $\rho$ and the $\tau_i$ are non negative, we immediately get that
$$\tsum_{i }\langle s_{i },U \rangle_{\widetilde{h}\cdot \eta} \langle U,s_{i }\rangle_{\widetilde{h}\cdot \eta} \leq  \langle U,\Gamma U\rangle_{\widetilde{h}\cdot \eta}.$$
Moreover,  since we have $\tr(\Gamma)=O(k^n)$ and $\frac{\rho}{k}=O(k^{n-1})$, there exists a constant $c'(k_0)$ such that
$$\tsum_{i }\langle s_{i },U \rangle_{\widetilde{h}\cdot \eta} \langle U,s_{i }\rangle_{\widetilde{h}\cdot \eta} \geq  \frac{1}{1+c'}\langle U,\Gamma U\rangle_{\widetilde{h}\cdot \eta}.$$
Now, if one considers $\lambda_{max}\geq 0$ the maximal eigenvalue of $\eta$ at $p$ and $v$ an associated eigenvector, we obtain that
\begin{equation} \label{3ineq}
\frac{1}{1+c'} \langle v,\Gamma v\rangle_{\widetilde{h}} \leq  \lambda_{max} \tsum_i \vert \langle v, s_i \rangle_{\widetilde{h}} \vert^2 \leq \langle v, \Gamma v \rangle_{\widetilde{h}}.
\end{equation}
Since the $(s_i)_{i=1,..,N}$ form a free family, one gets that  $\lambda_{max}$ belongs to a compact set in $C^0$ norm and also the solution $\eta$ of the considered equation, which is a definite positive hermitian endomorphism.
With Lemma \ref{annex1} we can see that differentiating (\ref{eqndiff}), we get
\begin{eqnarray}
\tsum_i s_i \langle .,s_i \rangle_{\widetilde{h}}\partial \eta
+ \frac{\rho}{k} \Pi^{\mathscr{F}\otimes L^k,\boldsymbol{\tau}}_{\widetilde{h}\cdot \eta,i}\partial \eta = \partial \Gamma 
- \partial  \left(\tsum_i s_i \langle .,s_i \rangle_{\widetilde{h}}\right)\eta.  \label{diffl}
\end{eqnarray}
Since $\Gamma$ is smooth and $M$ compact, one gets from the $C^0$ estimate of $\eta$,
a $C^0$ bound on $\partial \eta$, and by a similar way a $C^0$ bound on 
  $\overline{\partial}\eta$. By differentiating  (\ref{diffl}), one sees again that
  $\overline{\partial}\partial\eta$ is bounded in $C^0$ norm and consequently,
our solution $\eta$ is bounded in $C^{1,1}$ norm. Finally, with Arzela-Ascoli's Theorem 
we get $I$ closed and therefore $I=[0,\epsilon_k]$.\qed
\end{proof}

\begin{proposition} \label{obtenirpe}
Let us fix an integer $q\geq 1$. To each almost balanced metric $h_{k,q}$ corresponds a zero $\widetilde{h_q}$ of the moment map $\mu_{\mathcal{G}}$ on $\mathcal{Q}_0$ such that $$\Vert h_{k,q}- \widetilde{h_q} \Vert_{C^\alpha}=O\left(\frac{1}{k^{q-1-\alpha}}\right).$$
Moreover, we have the decomposition $$\int_M \langle s_i, s_j \rangle_{\widetilde{h_q}}\frac{\omega^n}{n!}=\delta_{ij}+ \boldsymbol{\eta}_{\widetilde{h_q}}$$
where $\boldsymbol{\eta}_{\widetilde{h_q}}$ is a matrix $N\times N$ such that
$$\vert \vert \vert \boldsymbol{\eta}_{\widetilde{h_q}} \vert \vert \vert = O\left(\Vert \boldsymbol{\sigma}_q(k)\Vert_{C^0} \right), $$
where $\boldsymbol{\sigma}_q(k)$ is given by Proposition \ref{Br75}.
\end{proposition}
\begin{proof}
Indeed, by Proposition \ref{Br75}, there exists a metric $h_{k,q}\in Met(\mathcal{F}\otimes L^k)$ such that
\begin{equation}\label{eqnpe2}
\tsum_{i }s_{i }\langle .,s_{i }\rangle _{{h}_{k,q}}=%
\frac{N+\epsilon_k V \sum_j \varepsilon_j r_j}{rV}Id 
+\boldsymbol{\sigma}_{q}(k)-\epsilon_k \tsum_j \varepsilon_j \pi^{\mathscr{F}\otimes L^k}_{{h}_{k,q}},
\end{equation}
with $\left\Vert \mathsf{\boldsymbol{\sigma}}_{q}(k)\right\Vert
_{C^{\alpha+2}}\leq C_{q,\alpha}k^{n-q-1}$, $(s_i)_{i=1,..,N}$
a $Hilb_{\omega}(h_{k,q})$-orthonormal basis and still $\varepsilon_j=\frac{\tau_j}{k}$.
Consider the metric $\widetilde{h_q}= h_{k,q}(\eta .,.)$ 
and a point $(s_1,...,s_N,A,\theta_1,...,\theta_{m-1})\in \mathcal{Q}_0$. By perturbation of our almost balanced metric, we can get a zero of the moment map $\mu_{\mathcal{G}}$ on $\mathcal{Q}_0$. Indeed, by Lemma \ref{avantpe} with the data $$\Gamma=\frac{N+\epsilon_k V \sum_j \varepsilon_j r_j}{rV}Id=\mathbf{O}(k^n) \hspace{1.5cm} \widetilde{h}=h_{k,q} \hspace{1.5cm}\rho=\epsilon_k,$$ we see that for $k$ fixed such that $k\geq k_0$ we can find a smooth solution $\eta$ of (\ref{leqn}). Then, the metric $\widetilde{h_q}=h_{k,q}\cdot \eta$ is  a  zero of the moment map $\mu_{\mathcal{G}}$.

 Now the fact that $\widetilde{h_q}=h_{k,q}\cdot \eta$ is close to  $h_{k,q}$ is a consequence of the relation (\ref{2ineq}) which gives us for  $\rho=\epsilon_k$ and $v$ an eigenvector associated to the eigenvalue $\lambda \geq 0$ of $\eta$,
 \begin{eqnarray}
 (\lambda-1) \displaystyle{\tsum}_i |\langle v,s_i\rangle_{h_{k,q}}|^2 &=& \langle v,\Gamma v \rangle_{h_{k,q}} - \epsilon_k \displaystyle{\tsum}_j \varepsilon_j \langle \nonumber v,\pi^{\mathscr{F}\otimes L^k}_{\widetilde{h_q},j}v\rangle_{h_{k,q}} \\
 &&   - \displaystyle{\tsum}_i |\langle v, s_i\rangle_{h_{k,q}}|^2,  \nonumber \\
 &=& \langle v,\Gamma v \rangle_{h_{k,q}} - \epsilon_k \displaystyle{\tsum}_j \varepsilon_j \langle v,\pi^{\mathscr{F}\otimes L^k}_{h_{k,q}\cdot (Id + (\eta-Id)),j}v\rangle_{h_{k,q}} \nonumber \\
 && - \displaystyle{\tsum}_i |\langle v, s_i\rangle_{h_{k,q}}|^2,  \nonumber \\
 &=& \langle v,\Gamma v \rangle_{h_{k,q}} - \epsilon_k \displaystyle{\tsum}_j \varepsilon_j \langle v,\pi^{\mathscr{F}\otimes L^k}_{h_{k,q},j}v\rangle_{h_{k,q}} \nonumber \\
  &&- \langle v, \frac{1}{k}\Pi^{\mathscr{F},\mathbf{\tau}}_{h_{k,q}}(\eta-Id)v \rangle_{h_{k,q}} \nonumber \\
   &&
   +\langle v,\frac{1}{k}\vartheta(Id-\eta)v\rangle_{h_{k,q}}
  - \displaystyle{\tsum}_i |\langle v, s_i\rangle_{h_{k,q}}|^2, \label{leqn5}
\end{eqnarray}
using Lemma \ref{annex1}. Here $\vartheta(Id-\eta)$ is an endomorphism of $\mathcal{F}\otimes L^k$ such that $\Vert\vartheta(Id-\eta)\Vert_{C^0}=O(\Vert Id-\eta \Vert^2_{C^0})$.
  From another side, we get that 
 \begin{eqnarray}
 \epsilon_k \left\Vert \frac{1}{k}\Pi^{\mathscr{F},\mathbf{\tau}}_{h_{k,q}}(Id-\eta) \right\Vert_{C^0} &\leq& \frac{\epsilon_k}{k} \left( {\tsum}_j \tau_j r_j(r-r_j) \right) \Vert \eta-Id \Vert_{C^0}.
  \label{leqn4}
\end{eqnarray}
Since $h_{k,q}$ satisfies equation (\ref{eqnpe2}) by definition,  we get from (\ref{leqn4}) and (\ref{leqn5}) that for some constants $c_0,c_0',c_0''$, independent of $k$,
 \begin{eqnarray*}
 \Vert \eta - Id \Vert_{C^0} \left(1-\frac{c_0\epsilon_k}{k^{n+1}}\right) - 
 \frac{c'_0\epsilon_{k}}{k^{n+1}}\Vert \eta - Id \Vert^2_{C^0} &\leq& \frac{c_0''}{k^n}\left(\Vert \Gamma - \Gamma \Vert_{C^0} +
  \Vert \boldsymbol{\sigma}_{q}(k) \Vert_{C^0} \right), \\
  &\leq & \frac{c_0''}{k^n} \Vert \boldsymbol{\sigma}_{q}(k) \Vert_{C^0}.
 \end{eqnarray*}
 We obtain the expected estimate for $k$ sufficiently large. 
Finally, since the $s_i$ are orthonormal respectively to $Hilb_{\omega}(h_{k,q})$, we have
\begin{eqnarray*}
\int_M \langle s_i,s_j \rangle_{\widetilde{h_q}}\frac{\omega^n}{n!}&=&
\int_M \langle s_i,s_j \rangle_{h_{k,q}}\frac{\omega^n}{n!}+\int_M \langle (\eta- Id)s_i,s_j \rangle_{h_{k,q}}\frac{\omega^n}{n!}, \\
&=&\delta_{ij}+\int_M \langle (\eta- Id)s_i,s_j \rangle_{h_{k,q}}\frac{\omega^n}{n!}.
\end{eqnarray*}
We use now the first part of the proof and Cauchy-Schwartz inequality.\qed
\end{proof}

\subsection{Explicit formulas and analytic estimates\label{Br52}}
Fix $q$ a  positive integer. From the almost balanced metrics 
 $h_{k,q}$, we have just got a point $z$ in the symplectic quotient $\mathcal{Z}=\mathcal{Q}_{0}\boldsymbol{//}\mathcal{G}$,
  with a metric $\langle .,. \rangle=\widetilde{h_q}\in Met(\mathcal{F}\otimes L^k)$ which satisfies
\begin{equation}  \label{prec}
\sum_{i=1}^N s_{i}\left\langle .,s_{i}\right\rangle _{\widetilde{h_q}}+\epsilon_k
\sum_{j=1}^{m-1} \varepsilon_j \pi^{\mathscr{F}\otimes L^k}_{\widetilde{h_q},j}=Cst_k Id,
\end{equation}
where $\mathfrak{q}=(s_{1},...,s_{N},A,{\theta_1},...,{\theta_{m-1}})\in
\mathcal{Q}_{0}$ is a lifting of $z$.  We now study the $SL\left( N\right) $-orbit of the point $z$. We are looking for  an estimate of the quantity $\boldsymbol{\Lambda }_{z}$ in order to apply 
 Proposition \ref{Br16}. We will apply the following lemma inspired from \cite[Lemma 18]{Do3} in order  to ease the computation of $\boldsymbol{\Lambda }_{z}$.

\begin{lemma} \label{double}
Let $\widehat{\nu }_{\widehat{z}}^{\mathcal{Q}_{0},SU(N)}:\mathfrak{su}(N)\rightarrow T_{\widehat{%
z}}\mathcal{Q}_{0}$ and $\widehat{\nu }_{\widehat{z}}^{\mathcal{Q}_{0},
\mathcal{G}}:Lie(\mathcal{G}^{\mathbb{C}})\rightarrow T_{\widehat{z}}\mathcal{Q}_{0}$
be the infinitesimal actions induced by $SU(N)$ and $\mathcal{G}^{\mathbb{C}}$ on $\mathcal{Q}_{0}$ and let $z\in
\mathcal{Z}$ be represented by $\widehat{z}\in \mathcal{Q}_{0}.$ Then for all
 $\xi \in \mathfrak{su}(N),$%
\begin{equation*}
\left\langle \mathfrak{q}_{z}^{SU\left( N\right) }(\xi ),\xi \right\rangle
=\left\vert \underline{\pi }\left( \widehat{\nu }_{\widehat{z}}^{\mathcal{Q}%
_{0},SU(N)}(\xi )\right) \right\vert ^{2},
\end{equation*}%
with $\underline{\pi }:T_{\widehat{z}}\mathcal{Q}_{0}\rightarrow T_{\widehat{%
z}}\mathcal{Q}_{0}$ orthogonal projection onto $\Im\left( \widehat{\nu
}_{\widehat{z}}^{\mathcal{Q}_{0},\mathcal{G}}\right) ^{\bot }$. In
particular,
\begin{equation*}
\boldsymbol{\Lambda }_{z}=\left( \underset{\xi \in \mathfrak{su}(N)}{\min }\frac{%
\left\vert \underline{\pi }\left( \widehat{\nu }_{\widehat{z}}^{\mathcal{Q}%
_{0},SU(N)}(\xi )\right) \right\vert }{\left\vert \xi \right\vert }\right)
^{-2}.
\end{equation*}
\end{lemma}

\vspace{1cm}
\noindent Consider the matrix
$$\mathsf{A}=(a_{i j })_{i j }\in \sqrt{-1}\mathfrak{su}(N)$$
and let  
\begin{equation*}
\sigma _{i }=\tsum_{j=1}^N a_{i j }s_{j }, \hspace{0.6cm}i=1,...,N
\end{equation*}%
be the induced basis by the infinitesimal action of $\mathsf{A}$. We look for the 
projection of $$\underline{\sigma }=(\sigma _{1},...,\sigma
_{N},0,...,0)$$ from the space  $C^{\infty}(M,\mathcal{F}\otimes L^{k})^{N}\times \prod_i C^{\infty}(M,\theta_i^*TGr(r_i,\mathcal{F}))$
onto the orthogonal complement of the subspace
\begin{eqnarray*}
\mathcal{P}=\left\{ \hspace{-0.1cm}
\begin{array}{r}
\left(\mathfrak{g}  s_{1},...,\mathfrak{g} s_{N},(Id-\pi_{\widetilde{h_q},1}^{\mathscr{F}\otimes L^k})\mathfrak{g} \pi_{\widetilde{h_q},1}^{\mathscr{F}\otimes L^k},...,(Id-\pi_{\widetilde{h_q},m-1}^{\mathscr{F}\otimes L^k})\mathfrak{g} \pi_{\widetilde{h_q},m-1}^{\mathscr{F}\otimes L^k}\right), \\
\text{s.t. }\mathfrak{g}\in Lie(\mathcal{G}^{
\mathbb{C}})
\end{array}
\right\}
\end{eqnarray*}
which is the image of the infinitesimal action
of $\mathcal{G}^{\mathbb{C}}$ of the point $\mathfrak{q}$,
 since $\theta^*_i TGr(r_i,\mathcal{F})\simeq Hom(\Im(\theta_i),\Im(\theta_i)^{\perp_{\widetilde{h_q}}})$.

\begin{notation}
We set as self-adjoint operator on $End(\mathcal{F}\otimes L^k)$,
$$\mathscr{B}:X \mapsto \frac{\epsilon_k}{Cst_k} \sum_{j=1}^{m-1} \varepsilon_j \pi_{\widetilde{h_q},j}^{\mathscr{F}\otimes L^k} X \pi_{\widetilde{h_q},j}^{\mathscr{F}\otimes L^k}.$$
If the condition 
\begin{equation}
\frac{\epsilon_k}{Cst_k}\tsum_{j=1}^{m-1} \varepsilon_j <1 \label{condst}
\end{equation}
holds, we notice that it makes sense to consider the operator $(Id-\mathscr{B})^{-1}$.
\end{notation}
\begin{lemma}
With the previous notations, assume that one has $Cst_k=O(k^n)$
in equation (\ref{prec}). Let
\begin{eqnarray*}
B_{\mathsf{A}}&=&(Id-\mathscr{B})^{-1}\left(\frac{1}{Cst_k}\displaystyle{\tsum}_{i ,j }a_{j i }s_{i }\langle
.,s_{j }\rangle\right) \\
 &=&  \hspace{-0.07cm}\frac{1}{Cst_k}\displaystyle{\tsum}_{i ,j }a_{j i }s_{i }\langle
.,s_{j }\rangle \hspace{-0.05cm} + \hspace{-0.05cm}\frac{\epsilon_k}{(Cst_k)^2} \displaystyle{\tsum}_l \varepsilon_l \pi_{\widetilde{h_q},l}^{\mathscr{F}\otimes L^k} \left(\displaystyle{\tsum}_{i ,j }a_{j i }s_{i }\langle
.,s_{j }\rangle\right) \pi_{\widetilde{h_q},l}^{\mathscr{F}\otimes L^k}\hspace{-0.2cm} +..
\end{eqnarray*}
be the hermitian endomorphism of $End(\mathcal{F}\otimes L^k)$ induced by the matrix $%
\mathsf{A}$. Then the orthogonal projection of $\underline{\sigma }$ onto $\mathcal{P}$
is
\begin{equation*}
\underline{p}\hspace{-0.05cm}=\hspace{-0.1cm}\left(\hspace{-0.1cm}B_{\mathsf{A}}s_{1},.., B_{\mathsf{A}}s_{N},
(Id- \pi^{\mathscr{F}\otimes L^k}_{\widetilde{h_q},1})B_{\mathsf{A}}\pi_{\widetilde{h_q},1}^{\mathscr{F}\otimes L^k},..,(Id- \pi^{\mathscr{F}\otimes L^k}_{\widetilde{h_q},m-1})B_{\mathsf{A}}\pi_{\widetilde{h_q},m-1}^{\mathscr{F}\otimes L^k}\right)\hspace{-0.05cm}.
\end{equation*}
\end{lemma}

\begin{proof}
We need to prove that for all $\mathfrak{g}\in Lie(\mathcal{G}^{
\mathbb{C}})$,
\begin{equation*}
\tsum_i \hspace{-0.1cm}\big\langle B_{\mathsf{A}}s_{i}-\sigma _{i},\mathfrak{g}s_{i}\big\rangle + \epsilon_k \tsum_i \big\langle  \varepsilon_i (Id- \pi^{\mathscr{F}\otimes L^k}_{\widetilde{h_q},i})B_{\mathsf{A}}\pi^{\mathscr{F}\otimes L^k}_{\widetilde{h_q},i}\hspace{-0.1cm},(Id- \pi^{\mathscr{F}\otimes L^k}_{\widetilde{h_q},i})\mathfrak{g}\pi^{\mathscr{F}\otimes L^k}_{\widetilde{h_q},i}\big\rangle\hspace{-0.08cm} =\hspace{-0.08cm}0
\end{equation*}
which is equivalent to saying that
\begin{equation*}
\tsum_i B_{\mathsf{A}}s_{i}\otimes s_{i}^*-\sigma _{i}\otimes s_i^* + \epsilon_k \tsum_i \varepsilon_i (Id-\pi^{\mathscr{F}\otimes L^k}_{\widetilde{h_q},i}) B_{\mathsf{A}} \pi^{\mathscr{F}\otimes L^k}_{\widetilde{h_q},i} =0,
\end{equation*}
i.e. 
$$(Id-\mathscr{B}) B_{\mathsf{A}} =\frac{1}{Cst_k}\tsum_i \sigma _{i}\otimes s_i^*,$$
since one considers a point of $\mu _{%
\mathcal{G}}^{-1}\left( Cst\times Id\right)$ and therefore we have a determined metric   via (\ref{prec}). The fact that $$\varepsilon_i \epsilon_k =\frac{\tau_i k^{n-1}}{k^n}\frac{\chi(\mathcal{F} \otimes L^k)}{Vr-V\sum_j \frac{\tau_j k^{n-1}} {k^n} } = O(k^{n-1})$$ allows us to reach our conclusion since condition (\ref{condst}) is satisfied for $k$ large.  \qed
\end{proof}

\noindent We set
\begin{equation*}
\begin{array}{lllll}
\psi_{i }&=&\sigma _{i }-B_{\mathsf{A}}s_{i} &\hspace{1cm} &  1\leq i  \leq N,\\
\psi _{N+i }&=&(\pi^{\mathscr{F}\otimes L^k}_{\widetilde{h_q},i}-Id)B_{\mathsf{A}}\pi_{\widetilde{h_q},i}^{\mathscr{F}\otimes L^k}& &  0\leq i \leq m-1, \\
\underline{\psi }&=&(\psi_{1},...,\psi_{N+m-1}).
\end{array}
\end{equation*}

\begin{notation}
In all the following, we will denote $$\left\Vert \underline{\psi }\right\Vert^2_{L^{2}(\omega
)}=\tsum_{i=1}^{N+m-1}||\psi _{i }||_{Hilb_{\omega}(\widetilde{h_q})}^{2}=\tsum_{i }\int_M |\psi _{i }|^2_{\widetilde{h_q}}\frac{\omega^n}{n!},$$
and $\left\vert \left\vert \left\vert \mathsf{Q} \right\vert \right\vert \right\vert
_{L^{2}(\omega)}^{2}$ will stand for the $L^{2}$ operator norm
 induced by $\omega$ and the metric $\widetilde{h_q}$ 
of an endomorphism $\mathsf{Q}\in C^{\infty}(M,End(\mathcal{F}\otimes L^k))$.
In particular, we can write the quantities $\mathit{\Lambda }_{z},\boldsymbol{\Lambda }_{z}$ as
\begin{eqnarray}
\mathit{\Lambda }_{z}^{-1}=\underset{i\mathsf{A}\in \mathfrak{su}(N),||%
\mathsf{A}||=1}{\min } \left\Vert \underline{\psi }\right\Vert^2_{L^{2}(\omega
)},  \hspace{0.4cm}
\boldsymbol{\Lambda }_{z}^{-1}= \underset{i\mathsf{A}\in \mathfrak{su}%
(N),\left\vert \left\vert \left\vert \mathsf{A} \right\vert \right\vert \right\vert =1}%
{\min }\left\Vert \underline{\psi }\right\Vert^2_{L^{2}(\omega
)}. \hspace{0.4cm} \label{Br18} 
\end{eqnarray}
\end{notation}
We also need the following definitions.
\begin{definition}
Let $h_\mathcal{F}$ be a hermitian metric that we fix as a reference metric on  $\mathcal{F}$ and fix an integer $\alpha>2$. To the integer $k,$ we associate the metric $$\widetilde{h_{\mathcal{F}}}=h_{\mathcal{F}}\otimes h_{L}^{k}$$ on $\mathcal{F}\otimes L^{k}.$ We will say that for $R>0,$ another hermitian metric
 $\widetilde{h_1}=h_1\otimes h_L^k$ on $\mathcal{F}\otimes L^{k}$ constructed in a similar way
has $(R,\alpha)$-bounded geometry if the two following conditions are satisfied:
\begin{equation*}
\widetilde{h_1}>\frac{1}{R}\widetilde{h_{\mathcal{F}}}, \hspace{0.8cm} \left\Vert \widetilde{h_1}-\widetilde{h_{\mathcal{F}}}\right\Vert _{C^{\alpha}}<R,
\end{equation*}
where $\left\Vert .\right\Vert _{C^{\alpha}}$ designs the standard norm  $%
C^{\alpha}$ determined by the reference metric $%
\widetilde{h_{\mathcal{F}}}.$ 
These conditions are equivalent to 
\begin{equation*}
h_1>\frac{1}{R}h_{\mathcal{F}},\quad \left\Vert h_1-h_{\mathcal{F}}\right\Vert _{C^{\alpha}\left(
h_1\right) }<k^{\alpha/2}R.
\end{equation*}%
Clearly, up to a modification of $R,$ this definition is independent 
of the choice of the metric $h_{\mathcal{F}}.$
\end{definition}

\begin{definition}
Consider a point $(s_0,...,s_N,A,\theta_1,...,\theta_{m-1})\in \mathcal{Q}_0$ and $R$ a positive real. We will say that the basis $\left( s_{i}\right) _{i=1,..,N}$ has 
$(R,\alpha)$-bounded geometry if there exists a smooth hermitian metric  $
\widetilde{h_q}$ which satisfies the condition (\ref{prec}) and has $(R,\alpha)$-bounded geometry.
\end{definition}

We are ready to introduce the main result of this subsection.
\begin{theorem}\label{Br14}
Let $\mathscr{F}$ be a simple holomorphic filtration. For all
 $R>0,$ there exist some constants $C=C(R,h_{\mathcal{F}},h_{L})$ and $\varepsilon
(R,h_{\mathcal{F}},h_{L})<\frac{1}{2}$ such that if, for all $k,$ the basis $\left(
s_{i}\right) _{i=1..N}\in H^{0}(M,\mathcal{F}\otimes L^{k})$ satisfying (\ref{Br83})
has $(R,\alpha)$-bounded geometry with $\left\vert \left\vert
\left\vert \mathsf{\boldsymbol{\eta}}\right\vert \right\vert \right\vert
<\varepsilon $ and $Cst_k=O(k^n)$, then for all matrix $\mathsf{A}=(a_{i
j })_{i j }\in \sqrt{-1}\mathfrak{su}(N)$, we have
\begin{equation*}
\left\Vert \mathsf{A}\right\Vert \leq Ck\left\Vert \underline{\psi }%
\right\Vert _{L^{2}(\omega )},
\end{equation*}%
where $\underline{\psi }\in \mathcal{P}^{\bot }$ is the orthogonal projection
to $\mathcal{P}$ of $\underline{\sigma }$. For the corresponding point $z\in \mathcal{Z}$, we have
\begin{eqnarray*}
\mathit{\Lambda }_{z} \leq C^{2}k^{2}, \hspace{0.8cm}\boldsymbol{\Lambda}_{z}\leq C^{2}k^{2}.
\end{eqnarray*}
\end{theorem}

We shall postpone the proof of the theorem and introduce now some useful identities and lemmas.

First of all, we set the following decomposition with respect to the metric satisfying equation (\ref{prec}),
\begin{equation}
\left\langle s_{i },s_{j }\right\rangle _{L^{2}(\omega )}=\int_M \left\langle s_{i },s_{j }\right\rangle \frac{\omega^n}{n!}= \delta
_{i j }+\mathsf{\boldsymbol{\eta }}_{i j },  \label{Br83}
\end{equation}%
where $\boldsymbol{\eta }=(\boldsymbol{\eta }_{ij })$
is a trace free hermitian matrix $N\times N$ . Thus,
$\boldsymbol{\eta }\equiv 0$ if and only if $\mathscr{F}$ is $k$-balanced.

We will also use the following well-known facts (for all hermitian matrices $\mathsf{S,R}$ of size $N\times N$).
\begin{eqnarray}
\left\vert \left\vert \left\vert \mathsf{R} \right\vert \right\vert \right\vert & \leq & 
\left\vert \left\vert \mathsf{R} \right\vert \right\vert  \leq \sqrt{N} \left\vert \left\vert \left\vert \mathsf{R} \right\vert \right\vert \right\vert, \label{Br102aaa}  \\
\left\vert \mathtt{tr}(\mathsf{SRS})\right\vert &\leq &\left\Vert \mathsf{S}%
\right\Vert ^{2}\left\vert \left\vert \left\vert \mathsf{R}\right\vert
\right\vert \right\vert ,  \label{Br102aa} \\
\left\vert \mathtt{tr}(\mathsf{RS})\right\vert &\leq &\sqrt{N}\left\Vert \mathsf{S}%
\right\Vert \left\vert \left\vert \left\vert \mathsf{R}\right\vert
\right\vert \right\vert.  \label{Br102}
\end{eqnarray}%

\begin{lemma}
Under previous assumptions, if the basis of holomorphic sections $\left( s_{i}\right) _{i=1..N}\in H^{0}(M,\mathcal{F}\otimes L^{k})$ has  $(R,\alpha)$-bounded geometry, and
 $\left\vert \left\vert \left\vert \mathsf{\boldsymbol{\eta }}\right\vert
\right\vert \right\vert <\frac{1}{2},$ then
\begin{equation*}
||\mathsf{A}||^{2}\leq 2\left( \left\vert \left\vert \left\vert
B_{\mathsf{A}}\right\vert \right\vert \right\vert _{L^{2}(\omega
)}^{2}+\left\Vert \underline{\psi }\right\Vert _{L^{2}(\omega )}^{2}\right).
\end{equation*}
\end{lemma}
\begin{proof}
 To prove the second inequality, we use the fact that we have an orthogonal decomposition
of $C^{\infty}(M,\mathcal{F}\otimes L^{k})^N \times \prod_{i=1}^{m-1} C^{\infty}(M,\theta_i^*TGr(r_i,\mathcal{F}))$:%
\begin{equation*}
\underline{\sigma }=\underline{\psi }+\underline{p},
\end{equation*}%
with $\underline{p}\in \mathcal{P},$ $\underline{\psi }\in \mathcal{P}^{\bot
}.$ But,
\begin{eqnarray*}
\left\Vert \underline{\sigma }\right\Vert _{L^{2}(\omega )}^{2}&=&%
\tsum_{i ,j }|a_{i j }|^{2}+\tsum_{i ,j
,l }a_{i j }\boldsymbol{\eta } _{j l
}a_{l i }=\Vert \mathsf{A} \Vert^2+\tr(\mathsf{A}\boldsymbol{\eta}\mathsf{A}),
\end{eqnarray*}
and since $\left\vert \left\vert \left\vert \mathsf{%
\boldsymbol{\eta}}\right\vert \right\vert \right\vert <\frac{1}{2}$ we have
by (\ref{Br102aa}) that
\begin{equation*}
||\mathsf{A}||^{2}<2\left\vert \left\vert \underline{\sigma }%
\right\vert \right\vert _{L^{2}(\omega )}^{2},
\end{equation*}%
and thus
\begin{equation*}
\frac{1}{2}||\mathsf{A}||^{2}\leq \left\Vert \underline{\psi }\right\Vert
_{L^{2}(\omega )}^{2}+\left\Vert \underline{p}\right\Vert _{L^{2}(\omega
)}^{2}\leq \left\Vert \underline{\psi }\right\Vert _{L^{2}(\omega
)}^{2}+\left\vert \left\vert \left\vert B_{\mathsf{A}}\right\vert
\right\vert \right\vert _{L^{2}(\omega )}^{2}.
\end{equation*} \qed
\end{proof}

We use Poincaré inequality to evaluate the term $\left\vert \left\vert \left\vert B_{\mathsf{A}}\right\vert \right\vert \right\vert _{L^{2}(\omega )}^{2}$.

\begin{lemma}
Assume that the holomorphic filtration $\mathscr{F}$ is simple and that $\mathsf{A} \in \sqrt{-1}\mathfrak{su}(N)$. If the basis $\left( s_{i}\right) _{i=1..N}$ of $H^0(M,\mathcal{F}\otimes L^k)$ has $(R,\alpha)$-bounded geometry, then there exist two constants $C_{1},C_2$ that depend only
on $R$ and of the reference metric $h_{\mathcal{F}}$ on $\mathcal{F}$ such that for $k$ large enough, 
\begin{equation*}
\left\vert \left\vert \left\vert B_{\mathsf{A}}\right\vert \right\vert
\right\vert _{L^{2}(\omega )}^{2}\leq C_{1}\left\Vert \overline{%
\partial }B_{\mathsf{A}}\right\Vert _{L^{2}(\omega )}^{2}+C_2\left(\left\vert
\left\vert \left\vert \mathsf{\boldsymbol{\eta}}\right\vert \right\vert
\right\vert+\frac{1}{k}\right) ^{2}\left\Vert \mathsf{A}\right\Vert ^{2}.
\end{equation*}
\end{lemma}

\begin{proof}
The fact that $\mathscr{F}$ is simple implies that for any $\gamma >1$, there exists
 $c(h_{\mathcal{F}},\gamma)$ such that if $h\in Met(\mathcal{F})$ is a metric satisfying
$$ \gamma h_{\mathcal{F}}>h>\frac{1}{\gamma} h_{\mathcal{F}},$$
then 
$$ \left\vert \left\vert \varpi \right\vert \right\vert _{L^{2}(\omega )}^{2}\leq c\left\Vert \overline{%
\partial }\varpi \right\Vert _{L^{2}(\omega )}^{2}+\frac{1}{rV}%
\mid \int_{M}\mathtt{tr}(\varpi)dV \mid ^{2}$$
for all $\varpi \in End(\mathcal{F})$ such that                 $\varpi(\mathcal{F}_i)\subset\mathcal{F}_i$. This
 is simply the Poincaré inequality with respect to the metric $\omega $ (see also \cite[Lemma 25]{Do3}) whose volume is $V$. Now, since $B_{\mathsf{A}}$ is hermitian, we can decompose $B_{\mathsf{A}}$ under the form $$B_{\mathsf{A}}= \mathsf{T}_{B_{\mathsf{A}}}+ \mathsf{D}_{B_{\mathsf{A}}}+ \mathsf{T}^*_{B_{\mathsf{A}}},$$
where $\mathsf{T}_{B_{\mathsf{A}}}$ is upper triangular and $\mathsf{D}_{B_{\mathsf{A}}}$  diagonal. Let
$\Pi(B_{\mathsf{A}})= \mathsf{T}_{B_{\mathsf{A}}}+ \frac{1}{2}\mathsf{D}_{B_{\mathsf{A}}}$.  Then $\Pi(B_{\mathsf{A}})$ is an endomorphism of $\mathcal{F}$ such that
 $\Pi(B_{\mathsf{A}})(\mathcal{F}_i) \subset \mathcal{F}_i$. Therefore, 
\begin{equation*}
 \left\vert \left\vert \Pi(B_{\mathsf{A}})\right\vert 
\right\vert _{L^{2}(\omega )}^{2}\leq C_2\left\Vert \overline{
\partial }\Pi(B_{\mathsf{A}})\right\Vert _{L^{2}(\omega )}^{2}+\frac{1}{rV}%
\left( \int_{M}\frac{1}{2}\mathtt{tr}(B_{\mathsf{A}})\frac{\omega^n}{n!}\right) ^{2}.
\end{equation*}
Nevertheless, the fact that $B_{\mathsf{A}}$ is hermitian also gives us that
$$\left\vert \left\vert \Pi(B_{\mathsf{A}})\right\vert 
\right\vert _{L^{2}(\omega )}^{2}=\frac{1}{2}\left\vert \left\vert B_{\mathsf{A}}\right\vert 
\right\vert _{L^{2}(\omega )}^{2}$$ and $\left\Vert \overline{%
\partial }\Pi(B_{\mathsf{A}})\right\Vert _{L^{2}(\omega )}^{2}=\frac{1}{2}\left\Vert\overline{%
\partial }B_{\mathsf{A}}\right\Vert _{L^{2}(\omega )}^{2}$. Thus, we have
\begin{equation*}
 \left\vert \left\vert \left\vert B_{\mathsf{A}}\right\vert  \right\vert 
\right\vert _{L^{2}(\omega )}^{2}\leq \left\vert \left\vert B_{\mathsf{A}}\right\vert 
\right\vert _{L^{2}(\omega )}^{2}\leq C_{1}\left\Vert \overline{%
\partial }B_{\mathsf{A}}\right\Vert _{L^{2}(\omega )}^{2}+\frac{1}{rV}%
\left( \int_{M}\mathtt{tr}(B_{\mathsf{A}})\frac{\omega^n}{n!}\right) ^{2}.
\end{equation*}
Now, $\mathsf{A}$ is trace free and $\frac{1}{Cst_k}\int_M \tr(\sum_i s_i\langle .,s_i\rangle)dV=O(1)$. Therefore we notice with (\ref{Br102}) that for $k$ large enough,  there exists $c$ such that
\begin{eqnarray*}
\int_{M}\mathtt{tr}(B_{\mathsf{A}})\frac{\omega^n}{n!} 
&\leq& \frac{\sum_{i ,j }a_{i j } \boldsymbol{\eta}_{i j }}{Cst_k} +\sum_{p=1}^{\infty} \frac{1}{k^p} \left(\frac{\epsilon_k}{Cst_k}\right)^p \sum_{l=1}^{m-1} {c\tau_l^p r_l} \Vert \mathsf{A} \Vert^p,  \\
&\leq&  {C_2}\left(\left\vert \left\vert
\left\vert \mathsf{\boldsymbol{\eta}}\right\vert \right\vert \right\vert +\frac{1}{k}\right)\left\Vert \mathsf{A}\right\Vert 
\end{eqnarray*}
for $C_2 \geq 1$ large enough, and we get the expected estimate.\qed
\end{proof}

To get rid of the term $\left\vert \left\vert \left\vert \overline{\partial
}B_{\mathsf{A}}\right\vert \right\vert \right\vert _{L^{2}(\omega )}^{2}$, we need
the following lemmas.

\begin{lemma}
\label{Br37}Under previous assumptions, there exist two independent constants $C^{(1)}(j_0)$, $ C^{(2)}(j_0)$ such that for all $%
j\leq j_0,$
\begin{eqnarray*}
 \sum_{i }\left\vert \nabla ^{j}s_{i 
}(z)\right\vert ^{2} &\leq
&C^{(1)}k^{j+n}\text{ for all }z \in M, \\
\left\Vert \nabla ^{j}B_{\mathsf{A}}\right\Vert^2_{L^{2}(\omega )} &\leq
&C^{(2)}k^{j}||\mathsf{A}||^2.
\end{eqnarray*}
\end{lemma}

\begin{proof}
We have pointwise the following Poincaré inequality for any holomorphic section $f$,
\begin{equation*}
\left\vert \nabla ^{j}f(x)\right\vert _{h}^{2}\leq c\int_{B(x)}\left\vert f\right\vert _{h}^{2}\frac{\omega ^{n}}{n!},
\end{equation*}%
and for $\widetilde{h_q}$,
\begin{equation*}
\left\vert \nabla ^{j}f(x)\right\vert
_{\widetilde{h_q}}^{2}\leq ck^{j}\int_{B(x)}\left\vert f\right\vert _{\widetilde{h_q}}^{2}%
\frac{\omega ^{n}}{n!}\leq ck^{j}\int_{M}\left\vert f\right\vert _{\widetilde{h_q}}^{2}%
\frac{\omega ^{n}}{n!},
\end{equation*}%
where $B(x)$ is a geodesic ball centered at $x\in M$ and  $c$ depends only on $R$ and $\omega .$ Now, we sum up for all $i$ and use the fact that
\begin{eqnarray*}\tsum_{i }\left\vert s_{i
}(x)\right\vert ^{2} \leq \tr \left( {\tsum}_{i }s_{i }(x)\langle
.,s_{i }(x)\rangle \right)+ \tr \left(\epsilon_k{\tsum}_j \varepsilon_j \pi^{\mathscr{F}\otimes L^k}_{\widetilde{h_q},j}\right)\leq rCst_k.
\end{eqnarray*}
 Hence, we get the inequality since $Cst_k=O(k^n)$. \newline
Moreover, since $B_{\mathsf{A}}$ is not holomorphic, we look at $%
M^{\prime }=M\times \overline{M}$ (cf. \cite[p. 507]{Do3}) where $%
\overline{M}$ is equipped with the opposite complex structure and $%
p_{1},p_{2}$ are the projections on the first and second factor.\ To the 
connection and the metric on $L\rightarrow M$ correspond a connection and a metric on $\overline{L}\rightarrow \overline{M}$ equipped with the opposite complex structure. Let $%
\mathcal{F}^{\prime }\rightarrow M'$ be defined by $p_{1}^{\ast }\left( \overline{\mathcal{F}}%
\otimes \overline{L}^{k}\right) ^{\vee }\otimes p_{2}^{\ast }\left(\mathcal{F}\otimes
L^{k}\right) $. Let $s^{\vee }\in H^0 \left(\overline{M},\overline{\mathcal{F}}\otimes \overline{L}^{k}\right) ^{\vee }$ be the holomorphic section associated to $s\in H^0(M,\mathcal{F}\otimes L^k)$ via the $C^{\infty }$ isomorphism of bundle defined by the metric. Then, we set 
\begin{equation*}
\widetilde{B_{\mathsf{A}}}=(Id-\mathscr{B}^{-1})\left(\frac{1}{Cst_k}\tsum_{i ,j }a_{j i
}s_{i }\otimes s_{j }^{\vee }\right),
\end{equation*}%
which is a holomorphic  section of $\mathcal{F}'$. Thus, if we denote $\Sigma_{\mathsf{A}}=\sum_{i ,j }a_{j i }s_{i }\langle
.,s_{j }\rangle$ we notice by Cauchy-Schwarz inequality that 
$$\langle \Sigma_{\mathsf{A}}, \mathscr{B}^p(\Sigma_{\mathsf{A}}) \rangle \leq \Vert \Sigma_{\mathsf{A}} \Vert
\Vert \mathscr{B}^p(\Sigma_{\mathsf{A}}) \Vert \leq \left(\frac{\epsilon_k\sum_i \varepsilon_i}{Cst_k}\right)^p\Vert \Sigma_{\mathsf{A}} \Vert^2.$$ Then, we get
\begin{eqnarray*}
\left\Vert \widetilde{B_{\mathsf{A}}}\right\Vert _{L^{2}\left(
\omega \right) }^{2} \hspace{-0.23cm}&=& \langle \Sigma_{\mathsf{A}}\hspace{-0.05cm}+ \hspace{-0.05cm}\mathscr{B}(\Sigma_{\mathsf{A}})\hspace{-0.05cm}+\hspace{-0.05cm}\mathscr{B}^2(\Sigma_{\mathsf{A}})\hspace{-0.05cm}\hspace{-0.05cm}+\hspace{-0.05cm}..., \Sigma_{\mathsf{A}}\hspace{-0.05cm}+\hspace{-0.05cm}\mathscr{B}(\Sigma_{\mathsf{A}})\hspace{-0.05cm}+\hspace{-0.05cm}\mathscr{B}^2(\Sigma_{\mathsf{A}})\hspace{-0.05cm}+\hspace{-0.05cm}... \rangle_{L^2(\omega)}, \\
&=& \langle \Sigma_{\mathsf{A}}, \Sigma_{\mathsf{A}}\rangle_{L^2(\omega)}(1+ O(1/k)),\\
&=& \frac{(1 +O(1/k))}{(Cst_k)^2}\int_M \dsum\limits_{i,j,k,l}a_{ij}a_{kl}\langle
s_{i},s_{j}\rangle
\langle s_{k},s_{l}\rangle \frac{\omega^n}{n!},\\
&=&\frac{(1 +O(1/k))}{(Cst_k)^2}\mathtt{tr}\left( \mathsf{A}\left(Id 
+\mathsf{\boldsymbol{\eta }}\right) \left(Id+\,^{t}\overline{%
\mathsf{\boldsymbol{\eta }}}\right) \,^{t}\overline{\mathsf{A}}\right).
\end{eqnarray*}%
But, we know that $Cst_k=O(k^n)$. Since we have chosen a basis with $(R,\alpha)$-bounded geometry, we obtain by inequality (\ref{Br102aa}),
\begin{equation*}
\left\Vert \widetilde{B_{\mathsf{A}}}\right\Vert _{L^{2}\left(
\omega \right) }\leq Ck^{-n}||\mathsf{A}||.
\end{equation*}%
Now, we notice that $B_{\mathsf{A}}$ is simply the restriction 
of $\widetilde{B_{\mathsf{A}}}$ over the diagonal of $M^{\prime }$.
 We apply the inequality for a holomorphic section as in the first
step of the proof, and get
\begin{equation*}
\left\Vert \nabla ^{j}\widetilde{B_{\mathsf{A}}}\right\Vert
_{L^{2}(\omega )}^{2}\leq C^{(2)}k^{j}||\mathsf{A}||^{2},
\end{equation*}%
which allows us to conclude.\qed
\end{proof}

\begin{lemma} 
There exists a constant $C_{3}$ which depends only on $R$ and of the reference metric $h_{\mathcal{F}}$ such that for $k$ large enough, we have
\begin{equation}
\label{Br665}
\tsum_{i=1}^{N+m-1}||\overline{\partial }\psi _{i }||_{L^{2}(\omega
)}^{2}=\left\vert \left\vert \left\vert \overline{\partial }B_{\mathsf{A}%
}\right\vert \right\vert \right\vert _{L^{2}(\omega )}^{2}\leq
kC_{3}\left\Vert \underline{\psi }\right\Vert _{L^{2}(\omega )}||\mathsf{A}%
||.
\end{equation}
\end{lemma}

\begin{proof}
We shall begin to prove the LHS equality. Pointwise, we have
\begin{eqnarray}
 \tsum_{i} |\overline{\partial }\psi _{i}|^2&=&\tsum_i
\left\vert \overline{\partial }(B_{\mathsf{A}}s_i)\right\vert^2 +
 \tsum_i \epsilon_k \varepsilon_i \left\vert \overline{\partial }
 \left((Id-\pi_{\widetilde{h_q},i}^{\mathscr{F}\otimes L^k})B_{\mathsf{A}}\pi_{\widetilde{h_q},i}^{\mathscr{F}\otimes L^k}\right)\right\vert^2, \nonumber \\
 &=&\tsum_i
\left\vert \overline{\partial }(B_{\mathsf{A}})s_i\right\vert^2 +
 \tsum_i \epsilon_k \varepsilon_i \left\vert (Id-\pi_{\widetilde{h_q},i}^{\mathscr{F}\otimes L^k})\overline{\partial }
(B_{\mathsf{A}})\pi_{\widetilde{h_q},i}^{\mathscr{F}\otimes L^k}\right\vert^2, \label{Br505} 
 \end{eqnarray}
since the $s_i$ and $\theta_i$ are holomorphic. The RHS of (\ref{Br505}) is, by definition of the metric we have fixed, the operator norm $\left\vert \left\vert \left\vert  B_{\mathsf{A}} \right\vert  \right\vert  \right\vert^2$ at the considered point $p$ of $M$. By integration over the manifold, we have the first part of the lemma. \\
Now, we have with previous equality,
\begin{eqnarray*}
\left\vert \left\vert \left\vert \overline{\partial }B_{\mathsf{A}%
}\right\vert \right\vert \right\vert _{L^{2}(\omega )}^{2}
&\leq &\sqrt{\tsum_{i=1}^N||\Delta B_{\mathsf{A}}s_{i }||_{L^{2}(\omega
)}^{2}\tsum_{i=1}^N||B_{\mathsf{A}} s_{i }||_{L^{2}(\omega )}^{2}} \nonumber \\
&+&\epsilon_k\sqrt{\tsum_{j=1}^{m-1} \varepsilon_{j}^2 \left\Vert \Delta \left( B_{\mathsf{A}} {\pi_{h,j}^{\mathscr{F}\otimes L^k} }\right) \right\Vert _{L^{2}\left(\omega \right) }^{2} \tsum_{j=1}^{m-1} \left\Vert B_{\mathsf{A}}{\pi_{\widetilde{h_q},j}^{\mathscr{F}\otimes L^k}} \right\Vert_{L^{2}(\omega )}^{2}} 
\end{eqnarray*}%
and from another side,
\begin{eqnarray*}
\left\Vert \Delta \left( B_{\mathsf{A}}s_{i }\right) \right\Vert
_{L^{2}(\omega )}&\leq& \left\Vert \nabla ^{2}(B_{\mathsf{A}})\right\Vert
_{L^{2}(\omega )}+\left\Vert 2\nabla B_{\mathsf{A}}\cdot \nabla s_{i
}\right\Vert _{L^{2}(\omega )}, \\
\left\Vert \Delta ( B_{\mathsf{A}}{\pi_{\widetilde{h_q},j}^{\mathscr{F}\otimes L^k}})\right\Vert _{L^{2}(\omega
)} &\leq &\left\Vert \nabla ^{2}(B_{\mathsf{A}})\right\Vert
_{L^{2}(\omega )}+\left\Vert 2\nabla B_{\mathsf{A}}\cdot \nabla
{\pi_{\widetilde{h_q},j}^{\mathscr{F}\otimes L^k}}\right\Vert _{L^{2}(\omega )}.
\end{eqnarray*}%
We conclude using the last lemma (2$^{nd}$ inequality) and the first part of the proof.
\qed
\end{proof}

\textit{Proof of Theorem \ref{Br14}.} \ With the last lemmas, we have straightforward,
\begin{eqnarray*}
\left\Vert \mathsf{A}\right\Vert ^{2} &\leq &2\left( \left\vert
\left\vert \left\vert B_{\mathsf{A}}\right\vert \right\vert \right\vert
_{L^{2}(\omega )}^{2}+\left\Vert \underline{\psi }\right\Vert
_{L^{2}(\omega )}^{2}\right) , \\
&\leq &2\left( C_{1}\left\Vert \overline{\partial }B_{%
\mathsf{A}}\right\Vert _{L^{2}(\omega )}^{2}+C_2\left(\left\vert \left\vert
\left\vert \mathsf{\boldsymbol{\eta}}\right\vert \right\vert \right\vert
+\frac{1}{k}\right)^{2}\left\Vert \mathsf{A}\right\Vert ^{2}+\left\Vert \underline{\psi }%
\right\Vert _{L^{2}(\omega )}^{2}\right) , \\
&\leq &C_4\left( k\left\Vert \underline{\psi }\right\Vert
_{L^{2}(\omega )}||\mathsf{A}||+\left(\left\vert \left\vert \left\vert \mathsf{%
\boldsymbol{\eta}}\right\vert \right\vert \right\vert +\frac{1}{k}\right)^{2}\left\Vert \mathsf{A}%
\right\Vert ^{2}+\left\Vert \underline{\psi }\right\Vert _{L^{2}(\omega
)}^{2}\right) .
\end{eqnarray*}%
With the assumption $\left\vert \left\vert \left\vert \mathsf{\boldsymbol{\eta}}%
\right\vert \right\vert \right\vert <\varepsilon $ and $k$ sufficiently large, we choose $\varepsilon $ 
such that we have $C_4\left(\epsilon+\frac{1}{k}\right) ^{2}<\frac{1}{2}.$ Then, 
there exists a constant $C$, independent of $k$, such that
\begin{equation*}
\left\Vert \mathsf{A}\right\Vert ^{2}\leq C\left( k\left\Vert \underline{%
\psi }\right\Vert _{L^{2}(\omega )}||\mathsf{A}||+\left\Vert \underline{%
\psi }\right\Vert _{L^{2}(\omega )}^{2}\right) .
\end{equation*}%
If $k\left\Vert \underline{\psi }\right\Vert _{L^{2}(\omega )}||\mathsf{%
A}||\leq \left\Vert \underline{\psi }\right\Vert _{L^{2}(\omega )}^{2}$
then the result is clear since $k\geq 1$.\ Otherwise, after 
simplification, we get exactly the expected inequality.\newline
The second part of the theorem is a consequence of the equality (\ref{Br18}) and of the fact that $\left\vert \left\vert  \left\vert \mathsf{A} \right\vert \right\vert \right\vert \leq  \left\vert  \left\vert \mathsf{A}  \right\vert \right\vert$.\qed

\subsection{Approximation Theorem for $\boldsymbol{\tau}$-Hermite-Einstein metrics}

\quad In this section, we achieve the proof for Theorem \ref{Br39} 
 using the analytical estimate of the previous part and the construction of 
 an almost balanced metric. We fix now as a reference metric the metric  $h_{\mathcal{F}}=h_{\infty }$, which is a conformally  $\boldsymbol{\tau}$-Hermite-Einstein metric solution of equation (\ref{Br105}).

 For any trace free matrix $\mathsf{S}\in \sqrt{-1}\mathfrak{su}(N),$ we can consider the action of $SL(N)$  on $z \in \mathcal{Z}$ to obtain another point $e^{\mathsf{S}}\ast z$ of the  symplectic quotient. This gives us a new hermitian metric 
$\widetilde{h_{\mathsf{S}}} \in Met(\mathcal{F}\otimes L^{k})$ that still satisfies
 equation (\ref{prec}) and depends on $q$. Let $\boldsymbol{\eta}({\mathsf{S}})$ be the matrix satisfying the decomposition (\ref{Br83}) for this new hermitian metric  $\widetilde{h_{\mathsf{S}}}.$ Under these conditions and with the notations of Proposition \ref{Br75}, we have the following estimates.

\begin{lemma}
\label{Br17}Fix a real number $R>0$ and $\mathsf{S} \in \sqrt{-1}\mathfrak{su}(N)$ with
$\left\vert \left\vert \left\vert \mathsf{S}\right\vert
\right\vert \right\vert \leq \frac{1}{2}$. 
\newline
1. If $q>\alpha+1$, then there exists a constant $C_5$ (independent of  $k$ and $R$) such that if one has
\begin{equation*}
\left\vert \left\vert \left\vert \mathsf{S}\right\vert \right\vert
\right\vert+\frac{1}{k}%
\leq C_5 R,
\end{equation*}%
then the metric $\widetilde{h_{\mathsf{S}}}$ has $(R,\alpha)$-bounded geometry.\newline
2. There exists a constant $C_6$ (independent  of $k$) such that
\begin{equation*}
\left\vert \left\vert \left\vert {\boldsymbol{\eta}}({\mathsf{S}%
})\right\vert \right\vert \right\vert \leq {C_6}\left( \left\vert
\left\vert \left\vert \mathsf{S}\right\vert \right\vert \right\vert +||%
\mathsf{\boldsymbol{\sigma}}_{q}(k)||_{C^{0}}\right) .
\end{equation*}
\end{lemma}

\begin{proof}
First of all, the construction is invariant under $SU(N)$ action. So, we can assume
that $\mathsf{S}=diag(\lambda _{i})$ is a diagonal matrix
with $\sum_{i }\lambda _{i }=0$. We consider the two points ${z}=(s_1,...,s_N,A,\theta_1,...,\theta_{m-1})\in \mathcal{Z}$
and ${z}'=(e^{\lambda_1} s_1,...,e^{\lambda_N}s_N,A,\theta_1,...,\theta_{m-1})\in \mathcal{Z}$.
 From Lemma \ref{avantpe}, there exists a smooth endomorphism $\eta_\mathsf{S}\in End(\mathcal{F}\otimes L^k)$ such that
$$\mathcal{B}_{\widetilde{h_q},\epsilon_k}(\eta_\mathsf{S})=C_k Id.$$
Fix $\widetilde{h_{\mathsf{S}}}=\widetilde{h_q}\cdot \eta_\mathsf{S}$. Now, by definition,
$$\tsum_i s_i \langle .,s_i \rangle_{\widetilde{h_{\mathsf{S}}}} + \epsilon_k \tsum_j \varepsilon_j  \pi^{\mathscr{F}\otimes L^k}_{\widetilde{h_{\mathsf{S}}},j} = Cst_k Id + \tsum_i (1-e^{2\lambda_i})s_i \langle .,s_i \rangle_{\widetilde{h_{\mathsf{S}}}}.$$
We apply the same reasoning that in the proof of Proposition \ref{obtenirpe} using the estimate of  Lemma \ref{Br37} (1$^{\text{rst}}$ inequality). Thus, we have that $$\left\Vert \widetilde{h_q}-\widetilde{h_{\mathsf{S}}}\right\Vert_{C^{\alpha}} \leq
c_{1}\left\vert \left\vert \left\vert \mathsf{S}\right\vert \right\vert
\right\vert .$$
Moreover, by Proposition \ref{obtenirpe}, the metric $\widetilde{h_q}$ differs from the reference metric  $\widetilde{h_{\infty}}$ by a term of the form $O\left( 1/k\right) $ in $C^{\alpha}$ norm. Hence  $\left\Vert \widetilde{h_{\infty}}-\widetilde{h_{\mathsf{S}}}\right\Vert_{C^{\alpha}} < R$ by choosing suitably $C_5$. Finally, for $C_5$ chosen small enough, we can also ask that $\widetilde{h_{\mathsf{S}}} > \frac{1}{R} \widetilde{h_{\infty}}$ since the quantity $\vert \vert \vert \mathsf{S} \vert \vert \vert$ is bounded by assumption. \newline
For the second assertion, we notice that
\begin{eqnarray}
\hspace{-0.5cm}{\boldsymbol{\eta}\left( \mathsf{S}\right)}_{ij} \hspace{-0.1cm}
&=&\hspace{-0.1cm}\int_M \hspace{-0.1cm}\left\langle \eta_\mathsf{S} e^{\lambda_i}s_{i}, e^{\lambda_j}s_{j}\right\rangle 
_{\widetilde{h_q}}dV -\delta_{ij}, \nonumber \\
&=&\hspace{-0.1cm} \int_M \hspace{-0.1cm}\left\langle (\eta_\mathsf{S}-Id) e^{\lambda_i}s_{i}, e^{\lambda_j}s_{j}\right\rangle 
_{\widetilde{h_q}}dV\hspace{-0.1cm} + \hspace{-0.1cm}\left(
\int_M \hspace{-0.1cm}\left\langle e^{\lambda_i}s_{i}, e^{\lambda_j}s_{j}\right\rangle 
_{\widetilde{h_q}}dV \hspace{-0.05cm}-\hspace{-0.05cm}\delta_{ij} \hspace{-0.05cm}\right)\hspace{-0.05cm}. \label{leqn1}
\end{eqnarray}%
 From the first part of the proof and Proposition \ref{obtenirpe}, the first term of the RHS of (\ref{leqn1}) is  bounded in $C^0$ norm by a multiple of $\left\vert \left\vert \left\vert \mathsf{S}\right\vert\right\vert \right\vert^2$. By Proposition \ref{obtenirpe}, the second term of the RHS is bounded by a multiple of $\left\vert \left\vert \left\vert \mathsf{S}\right\vert\right\vert \right\vert+ ||\mathsf{\boldsymbol{\sigma}}_{q}(k)||_{C^{0}}$.\qed
\end{proof}

With the next proposition we check that all the assumptions of Proposition \ref{Br16}
are satisfied.
\begin{proposition} \label{lastp}
Let $R$ and $q$ be positive real numbers. \\ If we assume $\vert \vert \vert \mathsf{S} \vert \vert \vert \leq \min\{\delta,\delta k^{n-q+1}\}$ with $\delta(R,M,\mathscr{F})$
small enough, then the metric $\widetilde{h_{\mathsf{S}}}$ has $(R,\alpha)$-bounded geometry, the basis $(s_i)_{i=1,..,N}$ has $(R,\alpha)$-bounded geometry and $z=(s_1,...,s_N,A,\theta_1,...,\theta_{m-1})$
is a zero of the moment map $\mu_{\mathcal{G}}$. Moreover, $\mu_{SU(N)}=\boldsymbol{\eta}({\mathsf{S}})$ where the matrix ${\boldsymbol{\eta}}\left( \mathsf{S}%
\right)_{ij}$ satisfies (\ref{Br83}) with
$$\Vert \boldsymbol{\eta}({\mathsf{S}}) \Vert =O(k^{3n/2-q-1}).$$
\end{proposition}
\begin{proof}
 With the choice of $\left\vert \left\vert \left\vert \mathsf{S}\right\vert \right\vert
\right\vert +\frac{1%
}{k}\leq C_5 R,$ Lemma \ref{Br17} gives us that if  $\left\vert
\left\vert \left\vert \mathsf{S}\right\vert \right\vert \right\vert \leq \delta$
with $$\delta=
\min (C_5 R/2,1/2),$$
 then the metric $h_{\mathsf{S}}$ has $(R,\alpha)$-bounded geometry.
  From the almost balanced metrics $h_{k,q}$ which satisfy
\begin{equation}
\tsum_{i}s_{i}\langle .,s_{i}\rangle _{h_{k,q}}=C_k 
Id-\epsilon_k \tsum_j \varepsilon_j \pi^{\mathscr{F}\otimes L^k}_{h_{k,q},j}+{\boldsymbol{\sigma }}_{q}(k),
\end{equation}%
and $\left\Vert {\boldsymbol{\sigma }}_{q}(k)\right\Vert
_{C^{\alpha+2}}=O(k^{n-q-1})$ we can apply Proposition \ref{obtenirpe}. Under these conditions, we obtain from inequality (\ref{Br102aaa}) and Lemma \ref{Br17} that
\begin{equation*}
\Vert \mathsf{\boldsymbol{\eta }}({\mathsf{S}}) \Vert \leq \sqrt{N} \vert \vert \vert 
\mathsf{\boldsymbol{\eta }}({\mathsf{S}}) \vert \vert \vert \leq 
C_6 k^{n/2}(\delta +C_{7})k^{n-q-1}.
\end{equation*}
\qed
\end{proof}

Finally, we get the main result of this paper. 

\begin{theorem} \label{final}
\label{Br39}Let $\mathscr{F}$ be an irreducible holomorphic filtration equipped with a  $\boldsymbol{\tau}$-Hermite-Einstein metric $h_{HE}$ on a smooth projective manifold. Then $\mathscr{F}$ is balanced and there exists a sequence of balanced metrics $\boldsymbol{h}_{k}$
which converges in $C^{\infty}$ sense towards a metric $h_\infty$ conformally $\boldsymbol{\tau}$-Hermite-Einstein, i.e. towards $h_{HE}$ up to a conformal change.
\end{theorem}

\begin{proof}
We begin to prove that we can construct a sequence of balanced metrics which converges  in $C^{\alpha}$ topology towards the conformally   $\boldsymbol{\tau}$-Hermite-Einstein metric  $\widetilde{h_{\infty}}$, solution of equation (\ref{Br105}).

Let $\varepsilon$ be fixed by Theorem \ref{Br14} and $\delta$ by Proposition \ref{lastp}. Apply  Proposition \ref{lastp} with $R>0$, $q>\frac{3n}{2}+2+\alpha$
and $\Vert \mathsf{S} \Vert \leq \min\{\delta k^{n-q+1},\varepsilon\} \leq \delta$.
 We obtain a point $z\in \mathcal{Z}$,
represented by $\{s_{1},...,s_{N},A,\theta_1,...,\theta_{m-1}\}\in \mathcal{Q}_{0}.$
 By Theorem \ref{Br14}, we know that at that point, $\boldsymbol{%
\Lambda }_{z}\leq C^{2}k^{2}$. Again, by Proposition \ref{lastp}, we get
\begin{equation}
 \lambda \Vert\mathsf{\boldsymbol{\eta }}_{\mathsf{S}}\Vert \leq \lambda C_8k^{3n/2-q-1}\leq C_9k^{3n/2+1-q}. \label{thmf}
\end{equation}
Since $\left\vert
\left\vert \left\vert \mathsf{S}\right\vert \right\vert \right\vert \leq ||%
\mathsf{S}||,$ we want to apply Proposition \ref{Br16} with the data $$\mu _{SU\left(
N\right) }\left( z_{0}\right) =\boldsymbol{\eta }_{\mathsf{S}}, \hspace{1cm} \lambda =C^{2}k^{2},$$ and  $\delta$ given by Proposition \ref{lastp}. 
But inequality (\ref{thmf}) asserts that the quantity $\lambda \Vert \mu _{SU\left(
N\right) }(z_{0})\Vert $ can be chosen smaller than $\delta$ for $k$ large enough since $q>\frac{3n}{2}+2$. 
Then, by Proposition \ref{Br16}, we obtain that
\begin{equation*}
\left\vert \left\vert \left\vert \mathsf{S}\right\vert \right\vert
\right\vert \leq ||\mathsf{S}||\leq C_9k^{3n/2+1-q},
\end{equation*}%
and also the existence of a metric $\widetilde{h_{\mathsf{S}}}$ close in $C^0$ topology of $\widetilde{h_{\infty }}$ up to an error in $O\left( k^{3n/2-q+1}\right) $ and $k$-balanced. In fact, we also have
\begin{equation*}
\left\Vert {h_{\mathsf{S}}}-h_{\infty }\right\Vert _{C^{\alpha}}=O\left(
k^{3n/2-q+1+\alpha}\right) ,
\end{equation*}
and consequently, the convergence in $C^{\alpha}$ topology towards $h_{\infty}$
 for all $\alpha>2$. The theorem is proved using now Proposition \ref{40}.\qed
\end{proof}

\begin{corollary}
The Hermite-Einstein connection on the bundle $\mathcal{F}$ of the stable filtration $\mathscr{F}$
is unique up to an holomorphic automorphism of $\mathscr{F}$.
\end{corollary}

\section{Applications to the case of Vortex type equations}

In this part, we give some applications of Theorem \ref{Br39}, and in particular we study the case of coupled Vortex equations for which we have found it impossible to develop a direct method. Instead of working with the projective manifold $M$, we consider the Vortex equations as $\boldsymbol{\tau}$-Hermite-Einstein equations on a higher dimensional manifold and apply a dimensional reduction procedure.

\subsection{Equivariant holomorphic filtrations}
Consider $G$ a connected simply connected semisimple complex Lie group and $P\subset G$ a
parabolic subgroup of $G$. In particular, $G/P$ is a Flag manifold. We denote $K$
the maximal compact subgroup of $G$ and consider $X=M\times G/P$ with trivial action on $M$.
The K\"ahler structure on $M$ and $G/P$ define a K\"ahler structure on $X$.
\begin{remark}
One could also have construct $X$ as a projectively flat $G/P$-bundle. See \cite{B-G-K} for details.
\end{remark}
\begin{definition} \label{defequiv}
We will say that a coherent sheaf $\mathcal{F}$ is $G$-equivariant if the  action of $G$ on $X$ can be lifted holomorphically to $\mathcal{F}$. A $G$-equivariant filtration $\mathscr{F}$ on $X$ is a  filtration of $G$-invariant coherent subsheaves of a $G$-equivariant sheaf $\mathcal{F}$ on $X$, 
$$ \mathscr{F}:0= \mathcal{F}_0 \hookrightarrow ... \hookrightarrow \mathcal{F}_m=\mathcal{F}. $$
If the sheaves $\mathcal{F}_i$ are locally free, we say that the filtration is holomorphic.
\end{definition}

\begin{definition} 
A filtration $\mathscr{F}$ is $G$-equivariantly $\boldsymbol{\tau}$-stable (resp. semi stable) if $\mathscr{F}$ is $G$-equivariant and for any  $G$-invariant proper subfiltration $\mathscr{F}'\hookrightarrow \mathscr{F}$, we have
$$\mu_{\boldsymbol{\tau}}(\mathscr{F}') < \mu_{\boldsymbol{\tau}}(\mathscr{F}) \hspace{1cm}(\text{resp. }\leq).$$
Moreover, a filtration is $G$-equivariantly polystable if it is a direct sum of $\boldsymbol{\tau}$-stable $G$-equivariant filtrations with same slope $\mu_{\boldsymbol{\tau}}$.
\end{definition}

\begin{definition}
A filtration $\mathscr{F}$ is said to be $G$-equivariantly Gieseker $\mathbf{R}$-stable (resp. semi-stable) if $\mathscr{F}$ is $G$-equivariant and for $k$ large, one has for all proper  $G$-invariant subfiltration $\mathscr{F}'$ of $\mathscr{F}$,
\begin{eqnarray*}
 \frac{\mathscr{P}_{\mathbf{R},\mathscr{F}'}(k)}{r(\mathcal{F}')} &<& \frac{\mathscr{P}_{\mathbf{R},\mathscr{F}}(k)}{r(\mathcal{F})}  \hspace{1cm} (\text{resp.} \leq).
\end{eqnarray*}
\end{definition}
 
\begin{proposition}
Let $\mathscr{F}$ a $G$-equivariant holomorphic filtration on $X$. Then, $\mathscr{F}$
is $G$-equivariantly $\boldsymbol{\tau}$-stable (resp. Gieseker $\mathbf{R}$-stable) if and only if it is $G$-equivariantly indecomposable, and considered as a holomorphic filtration, has a direct sum decomposition into $\boldsymbol{\tau}$-stable (resp. Gieseker $\mathbf{R}$-stable) holomorphic filtrations $\mathscr{F}_j$, $1 \leq j \leq j_0$
such that these filtrations are images one from another by an element of $G$.
\end{proposition}
\begin{proof}
This is similar to \cite[Theorem 6]{GP0} and \cite[Theorem 2.2]{AC-GP} and thus will be omitted.
\end{proof}

One obtains a HKDUY correspondence for the holomorphic $G$-equivariant filtrations.
\begin{theorem-section}  Let $\boldsymbol{\tau}\in \mathbb{R}^{m-1}_{+}$ and $\mathscr{F}$ be a holomorphic filtration of length $m$.  A holomorphic filtration $\mathscr{F}$ is $G$-equivariantly $\boldsymbol{\tau}$-polystable
if and only if there exists a smooth $K$-invariant hermitian metric $h$ solution of the $\boldsymbol{\tau}$-Hermite-Einstein equation (\ref{Br904}).
\end{theorem-section}
\begin{proof}
See \cite[Theorem 4.7]{AC-GP3}.
\end{proof}

\begin{proposition}
Let $\mathscr{F}$ a holomorphic filtration over $X$ of length $m$. Then
$\mathscr{F}$ is $G$-equivariantly $\mathbf{R}$-Gieseker stable if and only if $Aut(\mathscr{F})=\mathbb{C}$ and for $k$ large, there exists a $K$-invariant metric $h_k\in Met(\mathcal{F}\otimes L^k)$ such that
\begin{equation*}
\widehat{\mathsf{B}}_{\mathcal{F} \otimes L^k,h_{k}}+\epsilon_k\sum_{j=1}^{m-1} \varepsilon_j \pi^{\mathscr{F}\otimes L^k}_{j,h_k}= \frac{N+\epsilon_k \sum_{j=1}^{m-1} \varepsilon_j r_j}{rV} Id_{\mathcal{F}\otimes L^{k}}
\end{equation*}
where
$$\epsilon_k=\frac{\chi(\mathcal{F}\otimes L^k)}{Vr-V\sum_{j>0} \varepsilon_j r_j}.$$
\end{proposition}
\begin{proof}
First of all, since $\mathscr{F}$ is $G$-equivariantly $\mathbf{R}$-Gieseker stable
it is in particular $\mathbf{R}$-Gieseker semistable and so one can consider for $k$ sufficiently large a Gieseker space $\widetilde{\mathfrak{G}}_k$ as it is done in Section 2.2.
The fact that the action lifts holomorphically implies that we have an action of $G$ on $\widetilde{\mathfrak{G}}_k$. By definition of the embeddings $i_{k,j}$ (cf. p.\pageref{i_k3}), the same holds for the space $\boldsymbol{\mathrm{\Pi}}$ and therefore on the zero set of the $\mu_{\mathscr{F},k}$. Now by uniqueness, this implies that the balanced point in $\boldsymbol{\mathrm{\Pi}}$ induced by the $i_{k,j}$ is $K$-invariant (the existence of this point is clear from Theorem \ref{Br403a}). But this point only depends on the choice of the metric $H_k\in Met(H^0(\mathcal{F}\otimes L^k))$, and therefore the balanced metric is $K$-invariant. This gives the result using the Fubini-Study map $FS_k$.
\end{proof}

With the previous results in hand, it is now clear that we have an equivariant version of Theorem \ref{Br39}, i.e that a $K$-equivariant $\boldsymbol{\tau}$-Hermite-Einstein metric
can be approximated by $K$-invariant balanced metrics. Indeed, the balanced metrics that we construct in the proof of Theorem \ref{Br39} are by uniqueness necessarily $K$-invariant.
\begin{theorem} Let $\mathscr{F}$ be an irreducible holomorphic filtration over $X$ equipped with a $K$-equivariant $\boldsymbol{\tau}$-Hermite-Einstein metric $h_{HE}$. Then $\mathscr{F}$ is balanced and there exists a sequence of $K$-equivariant balanced metrics $\boldsymbol{h}_{k}$
which converges in $C^{\infty}$ topology towards $h_{HE}$ up to a conformal change. \label{Approxequiv}
\end{theorem}

 \subsection{Filtrations, quivers and dimensional reduction}
We now suppose that $X=M \times G/P$, and the action on $M$ is trivial. 
We denote  $p:X\rightarrow  M$ and $q:X \rightarrow G/P$ the natural projections. 
By restriction, any $G$-equivariant vector bundle on $X$ defines a $P$-equivariant holomorphic bundle on $M\times P / P\simeq M$. Conversely to any holomorphic $P$-equivariant bundle $E$ on $M$, one can associate a $G$-equivariant bundle by considering the quotient of $G\times E$ by the action of $u \in P$ given by $u \cdot (g,e)=(g\cdot u^{-1},u\cdot e)$ for which one has an action of $g' \in G$ by $g'\cdot (g,e)=(g'g,e)$. This principle of induction and restriction can also be applied to coherent sheaves. In fact, this equivalence of categories between $G$-equivariant holomorphic vector bundle 
on $X$ and $P$-equivariant bundle on $M$ can be extended in the following framework developped in \cite{AC-GP3,AC-GP}.

\begin{proposition}
Any coherent $G$-equivariant sheaf $\mathcal{F}$ on $X$ admits a $G$-equivariant sheaf filtration
\begin{eqnarray*}
\mathscr{F}:0 \hookrightarrow \mathcal{F}_0 \hookrightarrow ... \hookrightarrow \mathcal{F}_m= \mathcal{F} & &\\
\mathcal{F}_i/\mathcal{F}_{i-1}\cong p^*(\mathcal{E}_i) \otimes q^*(\mathcal{O}(\lambda_i)) & & \hspace{0.75cm} 0\leq i \leq m,
\end{eqnarray*}
where $\lambda_i$ are increasing numbers and $\mathcal{E}_i$ are non zero coherent sheaves
on $M$ with trivial $G$-action. If $\mathcal{F}$ is a holomorphic vector bundle, then the $\mathcal{E}_i$ are also holomorphic vector bundles.
\end{proposition}

This motivates the following definitions (see \cite{AC-GP3,AC-GP} for details).
\begin{definition}
A quiver is a pair of sets $\mathcal{Q}=\lbrace\mathcal{Q}_v,\mathcal{Q}_a\rbrace$
together with two maps $h,t:\mathcal{Q}_a\rightarrow \mathcal{Q}_v$. The elements
of $\mathcal{Q}_v$ are called the vertices, and the elements of $\mathcal{Q}_a$
are called the arrows. For each $\overrightarrow{a} \in \mathcal{Q}_a$, the vertex $h\overrightarrow{a}$
is called the head of the arrow $\overrightarrow{a}$, and $t\overrightarrow{a}$ its tail. Moreover, $\mathcal{Q}$
will be locally finite, i.e we require $h^{-1}(v)$ and $t^{-1}(v)$ to be finite
for all $v \in \mathcal{Q}_v$. A trivial path at $v \in \mathcal{Q}_v$ consists of the vertex $v$ with no arrows. A non trivial path in $\mathcal{Q}$ is a sequence of arrows $p=\overrightarrow{a_0} \circ... \circ \overrightarrow{a_l}$ that can be composed, i.e $t\overrightarrow{a_{i-1}}=h\overrightarrow{a_i}$. A relation of a quiver $\mathcal{Q}$ is a formal finite sum $r=c_1p_1+...+c_np_n$ of paths $p_i$ with coefficients $c_i\in \mathbb{C}$. A quiver with relations $(\mathcal{Q},\mathcal{R})$,
is a pair consisting of a quiver $\mathcal{Q}$ and a set of relations $\mathcal{R}$ for 
$\mathcal{Q}$.
\end{definition}

\begin{definition} Let $\mathcal{Q}$ be a quiver. A $\mathcal{Q}$-sheaf $(\mathscr{E},\boldsymbol{\phi})$ is a collection of $\mathscr{E}$ of coherent sheaves $\mathcal{E}_v$ for each $v\in \mathcal{Q}_v$ together with a collection of morphisms
$\phi_{\overrightarrow{a}}:\mathcal{E}_{t\overrightarrow{a}} \rightarrow \mathcal{E}_{h\overrightarrow{a}}$ for each arrow
$\overrightarrow{a}\in \mathcal{Q}_a$ such that $\mathcal{E}_v$ is zero for all but finitely many $v\in \mathcal{Q}_v$. A holomorphic $\mathcal{Q}$-bundle is a $\mathcal{Q}$-sheaf such that all the sheaves are holomorphic vector bundles. For a $\mathcal{Q}$-sheaf, any path in $\mathcal{Q}$ induces a morphism of sheaves and the trivial path induces the identity morphism $id:\mathcal{E}_v \rightarrow \mathcal{E}_v$. The $\mathcal{Q}$-sheaf $(\mathscr{E},\boldsymbol{\phi})$ satifies a relation 
$r=c_1p_1+...+c_np_n$ if $\sum_i c_i\phi_{\overrightarrow{a_{i,0}}}\circ ... \circ \phi_{\overrightarrow{a_{i,l_i}}} =0$ 
where $l_i$ is the length of the path $p_i=\overrightarrow{a_{i,0}} \circ... \circ \overrightarrow{a_{i,l_i}}$. Let $\mathcal{R}$ be a set of relations of $\mathcal{Q}$. A $(\mathcal{Q},\mathcal{R})$-sheaf (resp. $(\mathcal{Q},\mathcal{R})$-bundle) is a $\mathcal{Q}$-sheaf (resp. $\mathcal{Q}$-bundle) satisfying the relations $\mathcal{R}$. 
\end{definition}

\begin{remark}
The notion of quiver is a natural generalisation of the notion of chain that appeared
in \cite{AC-GP}. A sheaf chain is a pair $\mathscr{C}=(\mathscr{E},\boldsymbol{\phi})$
where $\mathscr{E}=(\mathcal{E}_0,...,\mathcal{E}_{m})$ is a $(m+1)$-tuple of 
coherent sheaves and a $m$-tuple $\boldsymbol{\phi}=(\phi_1,...,\phi_m)$
of homomorphisms $\phi_i\in Hom(\mathcal{E}_i,\mathcal{E}_{i-1}).$ We will later need the notion of a stability for a chain. If we denote the $\boldsymbol{\alpha}$-slope  $$\mu_{\boldsymbol{\alpha}}(\mathscr{C})=\frac{\sum_{i=1}^m \deg(\mathcal{E}_i)-\sum_{i=1}^m \alpha_i rk(\mathcal{E}_i)}{\sum_{i=1}^m rk(\mathcal{E}_i)}$$ where $\boldsymbol{\alpha}=(\alpha_1,...,\alpha_m)$ is a collection of real numbers, then the chain $\mathscr{C}$ is called $\boldsymbol{\alpha}$-stable if $\mu_{\boldsymbol{\alpha}}(\mathscr{C}')<\mu_{\boldsymbol{\alpha}}(\mathscr{C})$ for all proper
subchains $\mathscr{C}'$ of $\mathscr{C}$.
\end{remark}

We have the following theorem from \cite[Theorem 2.5]{AC-GP3}.

\begin{theorem-section} There exists a one-to-one correspondence between the categories of 
$G$-equivariant holomorphic vector bundles on $X=M\times G/P$ and of holomorphic $(\mathcal{Q},\mathcal{R})$-bundles on $M$.
\end{theorem-section}

Let's give a simple example. We have a correspondence at the level of holomorphic objects between the extensions on $X$ of the form
$$ 0 \rightarrow p^*\mathcal{E}_0 \rightarrow E \rightarrow p^*\mathcal{E}_1 \otimes q^*\mathcal{O}(2)
\rightarrow 0$$
and the triple $(\mathcal{E}_0,\mathcal{E}_1,\phi_1)$ where $\phi_1 \in Hom(\mathcal{E}_1,\mathcal{E}_0).$ Indeed, 
 Kunneth formula gives us $H^1(X,p^*(\mathcal{E}_0\otimes \mathcal{E}_1^*)\otimes q^*\mathcal{O}(-2))\simeq H^0(M,\mathcal{E}_0\otimes \mathcal{E}_1^*)\otimes H^1(\mathbb{P}^1,\mathcal{O}(-2))$ and if we now fix an element in $H^1(\mathbb{P}^1,\mathcal{O}(-2))\simeq \mathbb{C}$, the homomorphism $\phi_1$ can be identified with the extension class defining $E$.

\subsection{Dimensional reduction and applications} \label{lastpar}

Let $\omega'=p^*\omega + q^*\omega_{\epsilon}$ (where $\omega_{\epsilon}$ is the $K$-invariant smooth K\"ahler form constructed in \cite[Lemma 4.8]{AC-GP3}) be a K\"ahler form on $X=M\times  G/P$. We have the following central theorem \cite[Theorem 4.13]{AC-GP3} of dimensional reduction between quivers and $\boldsymbol{\tau}$-Hermite-Einstein filtrations.

\begin{theorem-section} \label{centralAC-GP}
Let $\mathcal{F}$ be a $G$-equivariant vector bundle on $X$, and 
let $\mathscr{F}$ be a $G$-equivariant filtration associated to $\mathcal{F}$ of length $(m+1)$. Let  $(\mathscr{E},\boldsymbol{\phi})$ the corresponding holomorphic $(\mathcal{Q},\mathcal{R})$-bundle on $M$. Then $\mathscr{F}$ admits a $K$-invariant $\boldsymbol{\tau}$-Hermite-Einstein metric with respect to $\omega'$ if and only if for each $v\in\mathcal{Q}_v$ such that the holomorphic vector bundle $\mathcal{E}_v$ of $(\mathscr{E},\boldsymbol{\phi})$ is non trivial, there exists smooth hermitian metrics $h_v\in Met(\mathcal{E}_v)$ satisfying  
\begin{equation}
\sqrt{-1}n_v\Lambda F_{h_v}+\sum_{\overrightarrow{a}\in h^{-1}(v)}\phi_{\overrightarrow{a}}\circ \phi_{\overrightarrow{a}}^{*_{h_v}}- \sum_{\overrightarrow{a}\in t^{-1}(v)}\phi_{\overrightarrow{a}}^{*_{h_v}} \circ \phi_{\overrightarrow{a}}= \tau_v' Id_{\mathcal{E}_v}  \label{QuivVortex}  
\end{equation}
where $n_v=\dim(\rm{M_v})$ is the multiplicity of the irreducible representation $\rm{M_v}$ of $P$ attached to $v$.
\end{theorem-section}
The relation between the real positive numbers $(\tau_i)$ and $(\tau_v')$
is made clear in \cite[Section 4.2.2]{AC-GP3}. We will refer to the system (\ref{QuivVortex}) as 
a quiver Vortex equation. Many equations from litterature can be obtained as a particular case of a quiver Vortex equation. We shall now give the main examples of such equations.
\begin{itemize}
\item The case of coupled Vortex equations. Indeed, consider $X=M\times \mathbb{P}^1$ and the group action $SL(2)=SL(2,\mathbb{C})$ given by the trivial action over $M$ and the standard action on $\mathbb{P}^1$ via the natural identification $$\mathbb{P}^1=SL(2)/ \mathfrak{P},$$ where $\mathfrak{P}$ stands for the parabolic subgroup of lower triangular matrices of $SL(2)$. In that case studied in \cite{AC-GP}, the previous theorem can be rephrased in the following way. A holomorphic filtration $\mathscr{F}$ on $X$ admits a $SU(2)$-invariant $\boldsymbol{\tau}$-Hermite-Einstein metric 
respectively to $p^*\omega + q^*\omega_{FS}$ (where $\omega_{FS}$ denotes the Fubini-Study 
metric on $\mathbb{P}^1$) if and only if the chain $\mathscr{C}=(\mathscr{E},\boldsymbol{\phi})$ admits a $(m+1)$-tuple of hermitian metrics $\mathbf{h}=(h_0,...,h_{m})$ satisfying the following chain of Vortex equations (also called coupled Vortex equations),
\begin{eqnarray}
\sqrt{-1}\Lambda F_{h_0}+\frac{1}{2}\phi_1\circ \phi_1^{*_{h_0}}&=&\tau_0 Id_{\mathcal{E}_0} \nonumber, \\
\sqrt{-1}\Lambda F_{h_i}+\frac{1}{2}\left(\phi_{i+1}\circ\phi_{i+1}^{*_{h_i}}- \phi_i^{*_{h_i}}\circ \phi_i\right) &=& (\tau_i-2i) Id_{\mathcal{E}_i} \hspace{0.3cm} (1\leq i\leq  \hspace{-0.06cm}m\hspace{-0.06cm}-\hspace{-0.06cm}1), \nonumber\\
\sqrt{-1}\Lambda F_{h_m}-\frac{1}{2}\phi_m^{*_{h_m}}\circ \phi_m&=&(\tau_m-2m) Id_{\mathcal{E}_m}.\nonumber \\
\label{chaine3}
\end{eqnarray}
Note that these equations have been related to the computation of some twisted Gromov-Witten invariants \cite{OT0}. It covers the interesting case of holomorphic triples $(E_1,E_2,\phi)$ where $\phi \in H^0(X,Hom(E_1,E_2))$ which is simply a holomorphic chain of length 1. The coupled  Vortex equations relative to these triples have been extensively studied (see for instance \cite{GP,B-GP2,B-GP4}).
\item The case of Hitchin's self-duality equation \cite{Hi2} over a complex curve $C$,
\begin{equation}
 F_h^{\perp} + [\Phi, \Phi^{*_h}]= 0, \label{Higgs}
\end{equation}
where $\Phi\in H^0(M,End(E)\otimes K_C)$ and $F_h^{\perp}$ is the trace free part of the curvature $F_h$.  The notion of stability considered here is simply the usual one for the holomorphic vector bundles restricted to the $\Phi$-invariant subbundles. On the moduli space $\mathcal{M}$  of stable Higgs bundles of degree $0$ over $C$, we have the so-called  $S^1$ action of Hitchin,
$$g \cdot (A,\Phi)=(A, g\Phi)\in \mathcal{A}(E)\times  H^0(M,End(E)\otimes K_C)$$
which preserves the natural K\"ahler form on  $\mathcal{M}$. A stable $(E,A,\Phi)$ bundle represents a fixed point of this action if and only if there exists a Gauge transformation  $\vartheta$ such that $D_{A}\vartheta=0$ and $[\vartheta, \Psi]=\sqrt{-1}\Psi$ (cf. \cite{Hi3,Sim3}). Then, the Higgs bundle $E$ is called \textit{critical} and can be decomposed holomorphically $$E = \bigoplus_{i=0}^{d} E_i,$$ and $\vartheta$ acts with increasing weights $\lambda_{i} \in \mathbb{R}$ on each factor $E_i$, i.e one has a variation of Hodge structure. We also have some non trivial morphims $\Phi_i:E_i \rightarrow E_{i+1}\otimes K_{C}$ with $\Phi=\bigoplus_i \Phi_i$. If one denotes now  $\mathcal{E}_{i}=E_{d-i}\otimes K_{C}^{d-i}$, then we can define a holomorphic chain $(\mathscr{E},\boldsymbol{\phi})$ by considering the morphisms $$\phi_i=\Phi_{d-i} \otimes Id: \mathcal{E}_{i} \rightarrow \mathcal{E}_{i-1}.$$ Now, to find a solution $h\in Met(E)$ of (\ref{Higgs}) for a critical Higgs bundle over $C$ is equivalent to find for $0\leq i \leq d$ smooth metrics $h_i \in Met(\mathcal{E}_{i})$ 
solutions of the couplex Vortex equations associated to the holomorphic chain $(\mathscr{E},\boldsymbol{\phi})$. Indeed, this chain is $\boldsymbol{\alpha}$-stable with respect to the weight $\boldsymbol{\alpha}=((m-i)\deg(K))_{i=0,..,d}$ and up to a conformal change $e^\rho$ given by the potential of the fixed metric on $K_C$, we have $h=\bigoplus_{i=0}^{d} e^{-(d-i)\rho}h_i$.
\item The case of Witten triples (Vafa-Witten equations studied in \cite{Wi}). Let $L$ be a line bundle on a complex surface $S$ and $({\mathcal{L}}, \phi, \theta)$ a triple formed by a holomorphic structure  ${\mathcal{L}}$ on $L$, a holomorphic section $\phi \in H^0(M,\mathcal{L})$ and a morphism $\theta: \mathcal{L} \rightarrow K_S$. The triple $({\mathcal{L}}, \phi, \theta)$ is called $\beta$-stable if $\deg(\mathcal{L})<\beta$ and $\phi \neq 0$ or $\beta< \deg(\mathcal{L})$ and $\theta \neq 0$. A triple is $\beta$-stable if and only if $(\phi,\theta)\neq (0,0)$ and there exists a metric $h$ on $L$ satisfying the equation 
\begin{equation}
\sqrt{-1}\Lambda F_h + \frac{1}{2}(|\phi|^2_h-|\theta|^2_h)=\beta.  \label{VW}
\end{equation}
Suppose now that $\deg(\mathcal{L})<\beta$ and $\phi$ does not vanish (we will say that the triple is \textit{special}), then we can consider the chain  $$((K_S,\mathcal{L},\mathcal{O}),(\theta,\phi))$$ which is $\boldsymbol{\alpha}$-stable with $\boldsymbol{\alpha}=(\alpha_0,\beta,\alpha_2)$ if $\deg(\mathcal{L})<\beta+\deg(K_S)-\alpha_0$. If we choose $\alpha_0-\deg(K_S)> 0$ sufficiently small, then we know the existence of solutions for the associated couplex Vortex equations and in particular there exist two metrics $h_0\in Met(K_S)$ and $h_1\in Met(L)$ such that \begin{eqnarray}
\sqrt{-1}\Lambda F_{h_0}+\frac{1}{2}|\theta|^2&=&\alpha_0\label{1eme}, \\
\sqrt{-1}\Lambda F_{h_1}+\frac{1}{2}\left(|\phi|^2-|\theta|^2\right)&=&\beta. \label{2eme}
\end{eqnarray}
Note that here the norms  $|\theta|^2$ (resp. $|\phi|^2$) are computed with respect to both metrics $h_0$ and $h_1$  (resp. $h_1$ and a trivialisation frame on the structure sheaf). For $\alpha_0-\deg(K_S)$ sufficiently small, the term $A=|\phi|^2-|\theta|^2$ is going to be positive by (\ref{1eme}) and we can recover (\ref{VW}) from (\ref{2eme}) if we do a conformal change. Indeed, this leads us to find a smooth function $f$ on $S$ such that $\Delta f + Ae^f = B$ with $\int_S BdV=\beta-\deg(\mathcal{L})$. But from the work of Kazdan and Warner \cite{K-W}, the existence and uniqueness of $f$ holds as soon as $A>0$ and
$\int_S BdV >0$.
\item The case of Bradlow pairs \cite{Bra2} studied in \cite{GP}. If $(E,\phi)$
is a pair (i.e. $E$ is a holomorphic vector bundle and $\phi\in H^0(M,E)$)
on $(M,\omega)$ and if $F$ is given by extension on $X=M \times \mathbb{P}^1$ $$0 \rightarrow p^*E \rightarrow F \rightarrow q^*\mathcal{O}(2) \rightarrow 0$$
then $(E,\phi)$ is $\lambda$-stable in the sense of Bradlow if and only if
$F$ is Mumford-stable with respect to the polarization associated to $$p^*\omega + \frac{2V}{(r(E)+1)\lambda-\deg(E)}q^*\omega_{FS}.$$ Of course, there is a natural identification between the pair $(E,\phi)$ and the triple $(E,\mathcal{O},\phi)$. Nevertheless, except in the case when $rk(E)=1$, it is not possible in general to relate in a simple way Bradlow's equation \cite{Bra2} to  equation (\ref{QuivVortex}). Instead, one obtains a new equation, called \textit{almost} Vortex equation
$$\sqrt{-1}\Lambda F_{h}+e^{-u} \phi \otimes \phi^{*_h}=\lambda Id_E,$$
where $u$ depends on the choice of a trivialisation on the structure sheaf. Even if this metric $h$ and the solution of Bradlow's equation are not related by a conformal change, they define the same point in the moduli space of solutions \cite{B-GP}.
\end{itemize}

Now, suppose that we consider a metric (or a family of metrics) solution to one of the previous equations (\ref{QuivVortex}),(\ref{chaine3}),(\ref{Higgs}),(\ref{VW}). Thanks to Theorem \ref{centralAC-GP}, we know that these solutions are related to a $K$-invariant $\boldsymbol{\tau}$-Hermite-Einstein metric for a filtration on $X=M\times G/P$.
From the equivariant version of the Approximation theorem (Theorem \ref{Approxequiv}), we see that these solutions can be approximated by \lq algebraic metrics', i.e metrics coming from a G.I.T construction as a zero of certain moment maps in finite dimensional setting. Finally we get,

\begin{theorem}  \label{final2}
Let  $(\mathscr{E},\boldsymbol{\phi})$ an irreducible $(\mathcal{Q},\mathcal{R})$-bundle on a smooth projective manifold $M$ such that for each $v\in\mathcal{Q}_v$ with $\mathcal{E}_v$  non trivial, there exists a smooth hermitian metrics $h_v\in Met(\mathcal{E}_v)$ satisfying  
the quiver Vortex equation (\ref{QuivVortex}). Then, up to some renormalizations by conformal changes, every metric $h_v$ is the limit in $C^\infty$ topology of a sequence of algebraic metrics. In particular, the solutions of coupled Vortex equations, critical Hitchin's self-duality equations over a curve, special Vafa-Witten equations can be approximated by algebraic metrics. For an irreducible stable pair, the sequence of algebraic metrics  obtained converges to the solution of an almost Vortex equation.
\end{theorem}

\section{Appendix}

\subsection{Endomorphism $\Pi^{\mathscr{F},\boldsymbol{\tau}}_h$}

In this part, we gather some elementary and technical results which are used in Section \ref{section3}.

Let $\mathscr{F}$ be a holomorphic filtration of length $m$, $h$ a hermitian metric on $\mathcal{F}$
and $\pi_{h,i}^{\mathscr{F}}$ the $h$-orthogonal projection onto the bundle $\mathcal{F}_i \subset \mathcal{F}$. 
For two smooth hermitian metrics $h_1$ and $h_2$ on $\mathcal{F}$, we know that they are related by the existence of an endomorphism $\eta$ such that $$h_1(X,Y)=h_2(\eta X,Y)$$ and $\eta$ is hermitian with respect to $h_2$ and definite positive. In particular, it is well known that $F_{h_1}=F_{h_2}+\overline{\partial}( \eta^{-1} \partial_{h_2} \eta)$.  
\begin{notation}
We set by $h \cdot \eta$ the metric $h(\eta \cdot, \cdot)$ for $\eta \in End(\mathcal{F})$, hermitian with respect to $h$.
\end{notation}

\begin{lemma} \label{annex0}
For all $(m-1)$-tuple $\{\tau_1,..,\tau_{m-1}\}$ of real numbers and all $h$-hermitian endomorphism $\boldsymbol{\eta} \in End(\mathcal{F})$  such that $h'=h\cdot \boldsymbol{\eta}$, we have
\begin{equation*}
d\left(\tsum_{i=1}^{m-1} \tau_i \pi_{h',i}^{\mathscr{F}}\right)=
\tsum_{i=1}^{m-1} \tau_i \pi_{h,i}^{\mathscr{F}}  d\boldsymbol{\eta} \left(Id-\pi_{h,i}^{\mathscr{F}} \right).
\end{equation*}
\end{lemma}

\begin{proof}
First of all, we can restrict to one of the factors $\pi_{h,i}^{\mathscr{F}}$, 
$h$-orthogonal projection onto the subbundle $\mathcal{F}_i$. 
Let $t \mapsto \pi_{i}(t)$ be a one parameter family of projections
onto the bundle $\mathcal{F}_i$ such that  $\pi_{i}(0)=
\pi_{h,i}^{\mathscr{F}}$. Since we have the relations
\begin{eqnarray*}
\pi_{i}(t)\pi_{i}(t)=\pi_{i}(t), \hspace{2cm} \pi_{i}(0)\pi_{i}(t)=\pi_{i}(t), 
\end{eqnarray*}
we obtain that
$$\pi_{i}(0)\pi_{i}(0)'=\pi_{i}(0)',$$ 
i.e. $\Im(\pi_{i}(0)') \subset \mathcal{F}_i$. From another side,
$\pi_{i}(0)'\pi_{i}(0)=\pi_{i}(0)'-\pi_{i}(0)\pi_{i}(0)'=0$
and therefore $\ker(\pi_{i}(0)')\supset \mathcal{F}_i$. Thus,
$$\pi_{i}(0)'=\pi_{i}(0)\pi_{i}(0)'(Id-\pi_{i}(0)),$$
and the space of solutions of this equation is $Hom(\mathcal{F}_i^{\perp},\mathcal{F}_i)$.
We notice that the differential is necessarily 
$U(\mathcal{F}_i^{\perp})\times U(\mathcal{F}_i)$ invariant. Finally, we can apply Schur   lemma since $U(\mathcal{F}_i^{\perp})\times U(\mathcal{F}_i)$ acts irreducibly. Then, up to a multiplicative constant, the differential is given by
$$X \mapsto \pi_{i}(0)X(Id-\pi_{i}(0)).$$
Set $h_t=h \cdot (Id+ \eta_t)$ avec $\eta_0=0$ and 
choose a $h$-orthonormal basis $(e^0_j)_{j=1,..,r}$ for which the first $r_i$ vectors generate $\mathcal{F}_i$. Then the new $h_t$-orthonormal basis $(e_j^t)_{j=1,..,r}$ is given by
\begin{equation}
R(e_j^t)=(e_i^0), \label{base}
\end{equation}
where $R$ is the unique upper triangular matrix with positive diagonal coefficients that satisfies the relation  $R R^{*_{h_t}}=Id+\eta_t$.
Now, by differentiating (\ref{base}) at $t=0$, we have
for $j\leq r_i$,
$$ \left(de_j^{t}\right)_{t=0}= -\frac{1}{2}d\left(\eta_0\right)_{jj}e_j^0 - \tsum_{k<j}d\left(\eta_0\right)_{jk}e_k^0.$$
Thus,  differentiating $h_t$ at $t=0$, we have
$$d\left({\tsum}_{j=1}^{r_i} e_j^t\otimes {e_j^t}^{*_{h_t}}\right)_{t=0}=
\sum_{j=1}^{r_i} e_j^0\otimes {e_j^0}^{*_{h_0}}d\eta_0-\tsum_{j=1}^{r_i} e_j^0\otimes {e_j^0}^{*_{h_0}}d\eta_0\tsum_{k=1}^{r_i} e_k^0\otimes {e_k^0}^{*_{h_0}},$$
which allows us to conclude.\qed
\end{proof}

Now, by direct application of the previous lemma we have,
\begin{lemma} \label{annex1}
For all $(m-1)$-tuple $\{\tau_1,...\tau_{m-1}\}$ of real numbers and all
hermitian endomorphism $\boldsymbol{\eta} \in End(\mathcal{F})$, we have
\begin{equation*}
{\tsum}_{i=1}^{m-1} \tau_i \pi_{h\cdot (Id + \boldsymbol{\eta}),i}^{\mathscr{F}}={\tsum}_{i=1}^{m-1} \tau_i \pi_{h,i}^{\mathscr{F}}  
+\Pi^{\mathscr{F},\boldsymbol{\tau}}_h(\boldsymbol{\eta}) +\mathbf{O}(\boldsymbol{\eta}^2)
\end{equation*}
where we have fixed the endomorphism,
$$\Pi^{\mathscr{F},\boldsymbol{\tau}}_h: \boldsymbol{\eta} \mapsto \tsum_{i=1}^{m-1} \tau_i \pi_{h,i}^{\mathscr{F}}\boldsymbol{\eta}\left(Id-\pi_{h,i}^{\mathscr{F}}\right).$$
Here $\mathbf{O}(\boldsymbol{\eta}^2)$ represents an hermitian endomorphism such that its   Hilbert-Schmidt norm can be bounded by $O(\Vert \boldsymbol{\eta} \Vert^2_{C^0})$.
\end{lemma}

\subsection{Resolution of a certain elliptic equation}
We will need the following classical K\"ahler identities.
\begin{lemma}
Let $E$ be a hermitian holomorphic vector bundle $E$ over a K\"ahler manifold and $F_E$ be its Chern curvature. We have the commuting identities: 
\begin{eqnarray*}
[\Lambda,\bar{\partial}]=-\sqrt{-1}\partial^*, \hspace{1cm} [\Lambda,{\partial}]=\sqrt{-1}\bar{\partial}^*, \hspace{1cm}
\Delta_{\bar{\partial}}=\Delta_{{\partial}}+[\sqrt{-1}F_E,\Lambda].
\end{eqnarray*}
\end{lemma}

In all the following, we will assume that the $\{\tau_1,...,\tau_m \}$ are non negative.

\begin{lemma} \label{Br36e}
Let $\mathscr{F}$ be a simple holomorphic filtration over $M$ and let $\Psi:End(\mathcal{F})\rightarrow End(\mathcal{F})$ be a positive self-adjoint operator of order zero. Then, for all hermitian metric $h \in Met(\mathcal{F})$, it is always possible to find a smooth solution, which preserves the filtration, of the following elliptic system:
$$ \Lambda_{\omega} \overline{\partial}\partial \mathsf{Q}^{\prime } + \Psi(\mathsf{Q}') =\mathsf{Q}$$
for all smooth endomorphism $\mathsf{Q}$ such that $\mathsf{Q}(\mathcal{F}_i)\subset \mathcal{F}_i$ and $\int_M \tr(\mathsf{Q}) dV= 0$.\\
 Moreover, if $\mathscr{F}$ is a holomorphic filtration such that there exists $\mathcal{F}$ a conformally $\boldsymbol{\tau}$-Hermite-Einstein metric $h$, and if one has fixed   
 \begin{equation} 
 \Psi:U \mapsto
  \Pi^{\mathscr{F},\boldsymbol{\tau}}_h(U),
  \label{operateur}
 \end{equation}
 then $\mathsf{Q}$ is self-adjoint if and only if $\mathsf{Q}'$ is self-adjoint.
\end{lemma}

\begin{proof}
We need to see that the operator  $\Lambda_{\omega} \overline{\partial}\partial + \Psi$  has trivial kernel. Recall that this operator is elliptic (of order 2) positive and self-adjoint since the $\tau_i$ are positive. Let us consider the kernel of this operator:  $\left( \partial ^{\ast }\partial \right) U=0$ implies
$ \left\vert \partial U\right\vert _{h}=0 $, i.e. since $\mathscr{F}$ is simple, $U=\gamma Id$ with $\gamma$ constant on $M$. If $Id \in \ker \Psi$ then, by Fredholm alternative, the elliptic system admits a solution if
 $\langle Id, Q \rangle=\int_M \tr(\mathsf{Q}) dV= 0$. The uniqueness is obvious once it has been assumed that $\int_M \tr(\mathsf{Q}') dV=0$. 
 If $Id \notin \ker \Psi$, then the system admits a unique solution.
 The defined operator by (\ref{operateur}) is self-adjoint and positive. Moreover, again with K\"ahler identities,
\begin{eqnarray*}
\sqrt{-1}\Lambda\bar{\partial}\partial {\mathsf{Q}^{\prime }}^*&=&\Delta_{\partial}{\mathsf{Q}^{\prime }}^*,\\
&=& - \left(\Delta_{\bar{\partial}}{\mathsf{Q}^{\prime }}\right)^* ,\\
&=& \left(\Delta_{\partial}{\mathsf{Q}^{\prime }}-[\sqrt{-1}\Lambda F_h,{\mathsf{Q}^{\prime }}]\right)^*, \\
&=& \left(\Delta_{\partial}{\mathsf{Q}^{\prime }}-[{\tsum}_{i}\pi^{\mathscr{F}}_{h,i},{\mathsf{Q}^{\prime }}]\right)^*,
\end{eqnarray*}
and we get
\begin{equation*}
\left(\sqrt{-1}\Lambda \left( \overline{\partial }\partial \mathsf{Q}^{\prime }\right)
+ \tsum_{i}\pi^{\mathscr{F}}_{h,i} \mathsf{Q}^{\prime }
\right)^* = \sqrt{-1}\Lambda \left( \overline{\partial }\partial {\mathsf{Q}^{\prime }}^*\right)  +   \tsum_{i}\pi^{\mathscr{F}}_{h,i} {\mathsf{Q}^{\prime }}^*.
\end{equation*}
Hence,
\begin{equation*}
\left(\sqrt{-1}\Lambda \left( \overline{\partial }\partial \mathsf{Q}^{\prime }\right)
+ \Psi(\mathsf{Q}^{\prime })
\right)^* = \sqrt{-1}\Lambda \left( \overline{\partial }\partial {\mathsf{Q}^{\prime }}^*\right)  +   \Psi({\mathsf{Q}^{\prime }}^*),
\end{equation*}
and by uniqueness of the solution, we have that ${\mathsf{Q}^{\prime }}$ is hermitian and definite positive if and only if it is the case for ${\mathsf{Q}}$.\qed
\end{proof}

\begin{acknowledgement}
The author is particularly grateful to Philippe Eyssidieux for introducing him to the subject. He deeply thanks Richard Thomas for enlightening conversations and Eva Miranda for her constant encouragement and support.
\end{acknowledgement}

\nocite{Zh}
\bibliography{biblio-matann}

\end{document}